\documentclass[11pt,a4paper]{amsart}

\usepackage{a4wide} 
\usepackage[utf8]{inputenc}

\usepackage{amsthm, amsmath, amssymb, amsfonts}
\usepackage{amsbsy,amscd,amsopn,amsxtra,amstext}

\usepackage{mathtools}

\usepackage{dsfont}
\usepackage{latexsym}

\usepackage{hyperref}
\usepackage{cleveref}
\usepackage{xcolor}

\usepackage{verbatim}
\usepackage{graphicx}

\usepackage{caption}
\usepackage{subcaption}
\usepackage{algorithm}
\usepackage{algorithmic}
\usepackage{xpatch}
\usepackage{lineno}
\usepackage{enumitem}

\usepackage[a4paper, total={6in, 8in}]{geometry}

\usepackage[colorinlistoftodos]{todonotes}   
\usepackage{makecell}

\numberwithin{equation}{section}
\newtheorem{theorem}{Theorem}[section]

\newtheorem{lemma}[theorem]{Lemma}

\theoremstyle{definition}

\newcommand{\rset}{\mathbb{R}}
\newcommand{\cset}{\mathbb{C}}
\newcommand{\Rd}{\,\mathrm{d}}
\newcommand{\bA}{\bar{A}} %
\newcommand{\by}{\bar{y}} %
\newcommand{\TR}{\mathrm{Tr}\,} %
\newcommand{\CSP}{\cset^{d\times d}}

\def\R{\mathbb{R}}
\def\E{\mathbb{E}}

\def\SEP{\,|\,}
\def\EXP#1{e^{#1}}
\def\OPER#1{\widehat{#1}}
\def\OPERW#1{\widehat{#1}}
\def\PERIOD{\,.}
\def\COMMA{\,,}
\def\BIGO{\mathcal{O}}
\def\Id{\mathrm{I}}
\def\IU{\mathrm{i}}
\def\NORMFAC{\left(\frac{\sqrt M}{2\pi}\right)}

\def\HSPACE{\mathcal{H}}
\def\SSPACE{\mathcal{S}}
\def\FT{\mathcal{F}}

\def\MP{{\# }}
\def\XYH{\tfrac{1}{2}(x+y)}
\def\XXH#1#2{\tfrac{1}{2}({#1}+{#2})}

\def\NABLAR{\nabla'}
\def\DT#1{\dot{#1}}
\def\BARYT{\Upsilon}

\def\YT{Y}
\def\WBT{W^\beta_t}
\def\HALF{\tfrac{1}{2}}

\begin{document}

\ifx\MMAN\undefined

%
%

\title[Canonical mean-field molecular dynamics]{Canonical mean-field molecular dynamics derived from quantum mechanics}

\author[X. Huang]{Xin Huang}
\address{Institutionen f\"or Matematik, Kungl. Tekniska H\"ogskolan, 100 44 Stockholm, Sweden}
\email{xinhuang@kth.se}

\author[P. Plech\'a\v{c}] {Petr Plech\'a\v{c}}
\address{Department of Mathematical Sciences, University of Delaware, Newark, DE 19716, USA}
\email{plechac@udel.edu}

\author[M. Sandberg]{Mattias Sandberg}
\address{Institutionen f\"or Matematik, Kungl. Tekniska H\"ogskolan, 100 44 Stockholm, Sweden}
\email{msandb@kth.se}

\author[A. Szepessy]{Anders Szepessy}
\address{Institutionen f\"or Matematik, Kungl. Tekniska H\"ogskolan, 100 44 Stockholm, Sweden}
\email{szepessy@kth.se}

\else

%
%

\title{Canonical mean-field molecular dynamics derived from quantum mechanics}

\author{Xin Huang}
\address{%
Institutionen f\"or Matematik, Kungl. Tekniska H\"ogskolan, 100 44 Stockholm, Sweden;
\email{xinhuang@kth.se}
}

\author{Petr Plech\'a\v{c}}
\address{%
Department of Mathematical Sciences, Univ. of Delaware, Newark, DE 19716, USA;
\email{plechac@udel.edu}
}

\author{Mattias Sandberg}
\address{%
Institutionen f\"or Matematik, Kungl. Tekniska H\"ogskolan, 100 44 Stockholm, Sweden; 
\email{msandb@kth.se}
}

\author{Anders Szepessy}
\address{%
Institutionen f\"or Matematik, Kungl. Tekniska H\"ogskolan, 100 44 Stockholm, Sweden; 
\email{szepessy@kth.se}
}
\fi
%

\begin{abstract}
Canonical quantum correlation observables can be approximated by 
classical molecular dynamics. 
In the case of low temperature  the ab initio 
molecular dynamics potential 
energy is based on the ground state electron eigenvalue 
problem and the 
accuracy has been proven to be $\BIGO (M^{-1})$, provided the 
first electron eigenvalue gap is sufficiently large compared to the given 
temperature and $M$ is the ratio of nuclei and electron masses. 
For higher temperature eigenvalues corresponding to excited 
electron states are 
required to obtain $\BIGO (M^{-1})$ accuracy and the derivations 
assume that all electron eigenvalues are separated, which for 
instance excludes conical intersections. 
This work studies a mean-field molecular dynamics approximation
where the mean-field Hamiltonian for the nuclei is the partial trace 
$h:=\TR(H e^{-\beta H})/\TR(e^{-\beta H})$ with respect to the electron 
degrees of freedom 
and $H$ is the Weyl symbol corresponding to a quantum many body Hamiltonian 
$\OPERW{H}$. It is proved that the mean-field molecular dynamics approximates  
canonical quantum correlation observables with accuracy
$\BIGO (M^{-1}+ t\epsilon^2)$, for correlation time $t$ where $\epsilon^2$ 
is related 
to the variance of mean value approximation $h$. 
Furthermore, the proof 
derives a precise asymptotic representation of the Weyl symbol of the Gibbs density operator using a path integral formulation. Numerical experiments on a model problem with one nuclei and two electron states show that the mean-field dynamics 
has similar or better accuracy than
standard molecular dynamics based on the ground state electron eigenvalue. 
\end{abstract}

\maketitle

\section{Classical approximation of canonical quantum observables}

\subsection{Introduction to the approximations}\label{subsec_intro_approx}
We study approximation of {\it quantum time-correlation observables} at the quantum canonical ensemble for a system consisting of nuclei (slow degrees of freedom) 
and electrons (fast degrees of freedom) at the inverse temperature $\beta = 1/(k_B T)$, 
where $k_B$ is the Boltzmann constant and $T>0$ is the temperature.  We work in 
Hartree atomic units in which the reduced Planck constant $\hbar=1$, the electron charge $e=1$, the Bohr radius $a_0=1$ and the electron mass $m_e=1$. Thus the semiclassical parameter in the subsequent analysis is 
given by the ratio of nucleus and electron masses $M$. 
For example, in the case of a proton--electron system the ratio is $M = m_p/m_e \approx 1836$.

The full quantum system is described by the Hamiltonian operator which includes
the kinetic energy of nuclei and the electronic kinetic energy operator 
together with the operator describing interaction between electrons, with coordinates $x_e$, and nuclei, with coordinates $x$,
$$
-\frac{1}{2M} \Delta_{x} - \frac{1}{2} \Delta_{x_e} +\OPER{V}_e(x,x_e)\PERIOD
$$
In this work, in the spirit of Born-Oppenheimer (adiabatic) approximation, 
we replace the time evolution of 
electrons by the Schr\"odinger electron eigenvalue problem.  
We represent the 
electronic kinetic energy operator and the interaction operator, $\OPER{H}_e = -\frac{1}{2}\Delta_{x_e}+\OPER{V}_e(x,x_e)$ as a matrix-valued potential $\OPER{V}(x)$ obtained
by a representation of the operator 
$-\tfrac{1}{2} \Delta_{x_e} +\OPER{V}_e(x,x_e)$ on a finite-dimensional ($d$-dimensional) subspace
of suitable normalized electronic eigenfunctions
$\{\phi_k\}_{k=1}^d$, as $V(x)_{k\ell}=\langle\phi_k,\OPER{H}_e(x,\cdot)\phi_\ell\rangle$, described precisely in Section~\ref{Sec2}.
Hence we work with
the Hamiltonian operator 
\begin{equation}\label{ham_defin}
\OPER{H} = -\frac{1}{2M} \Delta\otimes \Id + {V}(x)\,. 
\end{equation}
The first term $-\tfrac{1}{2M}\Delta_x\otimes{\Id}$
represents the kinetic energy of the nuclei where
$\Id$ is the $d\times d$ identity matrix. 
The second term, ${V}(x)$, is the matrix-valued potential approximation to $\OPER{H}_e$
and does not depend on $M$.
We assume that this finite-dimensional approximation of
the electronic operator results in a Hermitian matrix-valued smooth confining potential 
$V:\rset^N\to \rset^{d\times d}$ that depends on the positions $
x_i\in\R^3$ of nuclei $i=1,2,\dots, N'$, where we set $N=3N'$. 
For the sake of simplicity, we assume that the nuclei have the same mass; in the case with different nuclei masses $M$ becomes a diagonal matrix, which can be transformed to the formulation  \eqref{ham_defin} with the same mass by the change of coordinates $M_1^{1/2}\bar x=M^{1/2}x$.

The large nuclei/electron mass ratio $M\gg 1$ is the basis of semiclassical analysis 
and implies a separation of time scales, for which nuclei represent slow and electrons 
much faster degrees of freedom. 
The Weyl quantization takes this scale separation into account. In particular for the
Hamiltonian operator $\OPER{H}$ the corresponding matrix valued Weyl symbol becomes
$H(x,p)= \tfrac{1}{2}|p|^2\Id + V(x)$ for the nuclei phase-space points 
$(x,p)\in \R^{N}\times\R^{N}$, as described more precisely in Section~\ref{Sec2}.

In order to study correspondence between the quantum time-correlation function and its 
classical counter part we work in Heisenberg representation for the time-dependent quantum
observables given by self-adjoint operators $\OPER{A}_t$ and $\OPER{B}_t$. We employ 
the Weyl quantization to link the quantum dynamics given by 
the Heisenberg equation to classical  Hamiltonian equations of motions on the phase space $(x,p)\in\rset^N\times\rset^
N$ of
nuclei positions and momenta and to  averaging with respect to a suitable canonical Gibbs
distribution on the phase space. 

More precisely,
given a quantum system defined by the Hamiltonian $\OPER{H}$ acting on wave functions in 
$L^2(\rset^{N},\cset^d)\equiv [L^2(\rset^N)]^d$ we denote $\OPER{\rho}=e^{-\beta \OPERW{H}}$ the density operator for a quantum Hamiltonian operator $\OPER{H}$ at the inverse temperature 
$\beta>0$ and
consider quantum correlation observables based on the normalized trace 
\begin{equation}\label{qm_trace}
\frac{\TR(\OPER{A}_t \OPER{B}_0  e^{-\beta \OPERW{H}})}{\TR( e^{-\beta \OPERW{H}})}\,, 
\end{equation}
and the symmetrized version
\begin{equation}\label{qms}
 \mathfrak{T}_{\mathrm{qm}}(t):= 
 \frac{\TR\big(\frac{1}{2}(\OPER{A}_t \OPER{B}_0 +\OPER{B}_0\OPER{A}_t) e^{-\beta \OPERW{H}}\big) 
 }{\TR( e^{-\beta \OPERW{H}})}\,,
 \end{equation}
for quantum observables $\OPER{A}_t=e^{{\mathrm{i}}t M^{1/2} \OPER H}\OPER A_0 e^{-{\mathrm{i}}t M^{1/2} \OPER H}$ and $\OPER{B}_0$ at times $t$ and $0$, respectively. 
That is, the quantum observables solve the Heisenberg-von Neumann equation ${\mathrm{ d}}\OPER{A}_t/{\mathrm{d}}t={\mathrm{i}}M^{1/2}[\OPER H,\OPER A_t]$, where $[\OPER H,\OPER A_t]=\OPER H\OPER A_t-\OPER A_t\OPER H$ is the commutator. Here $M^{-1/2}$ plays the role of the Planck constant $\hbar$.
With this time scale the nuclei move a distance of order one in unit time in the classical setting. 
The aim is to derive a mean-field molecular dynamics  approximation 
\begin{equation}\label{md_mf}
\mathfrak{T}_{\mathrm{md}}(t):=\frac{\int_{\rset^{2N}} A_0\big(z_t(z_0)\big) B_0(z_0) 
\TR( e^{-\beta H(z_0)})\Rd z_0}{\int_{\rset^{2N}}  \TR( e^{-\beta H(z_0)})\Rd z_0}
\end{equation}
of the correlation function \eqref{qms},
where $A_0$ and $B_0$ are the Weyl symbols for the initial quantum observables 
and $z_t:=(x_t,p_t)$ solves the Hamiltonian system
\begin{equation}\label{h}
\begin{split}
\DT{x}_t =&\;\;\;\;\nabla_p h(x_t,p_t)\,,\\
\DT{p}_t =&-\nabla_x h(x_t,p_t)\,,\\
\end{split}
\end{equation}
with the initial data $z_0:=(x_0,p_0)\in\rset^{2N}$ of nuclei positions and momenta. 
The trace $\TR(e^{-\beta H(z_0)})$ is over the space $\rset^{d\times d}$ of $d\times d$ matrices which can be also viewed as the trace operator with respect to electron degrees 
of freedom under the finite dimensional approximation of the electronic Hamiltonian.
On the other hand $\TR( e^{-\beta\OPER{H}})$ represents the trace acting on the space 
of trace operators on $L^2(\rset^N,\cset^d)$ which we can view as the trace with respect to both nuclei and electron degrees of freedom.

\medskip

Two main questions arise: (a) how should the mean field Hamiltonian approximation $h:\rset^{2N}\to \rset$ be chosen,
and (b) how small is the corresponding estimate for the approximation error $\mathfrak{T}_{\mathrm{qm}}-\mathfrak{T}_{\mathrm{md}}$\,?

Assume $\Psi:\rset^N\to \mathbb C^{d\times d}$ is a differentiable orthogonal matrix.
Based on certain regularity assumptions on $A_0$, $B_0$, $V$ and $\Psi$, we prove in 
Theorem~\ref{thm_mf} that  for the mean-field Hamiltonian $h:\rset^{N}\times \rset^N\to\rset$ 
defined by
\begin{equation}\label{mf-def}
h(z):=\frac{\TR (H(z)e^{-\beta H(z)})}{\TR (e^{-\beta H(z)})}\,,
\end{equation}
and the symbols $A_0$ and $B_0$, which are independent of the electron coordinates,
we have
 \begin{equation}\label{trace_est}
 |\mathfrak{T}_{\mathrm{qm}}(t)-\mathfrak{T}_{\mathrm{md}}(t)|=\BIGO (M^{-1}+t\epsilon_1^2+t^2\epsilon_2^2)\,,
 \end{equation}
 where the parameters $\epsilon_j^2$  are the variances 
\begin{equation}\label{epsilon_2_square_model}
\begin{split}
\epsilon_1^2 &= \frac{\|\TR\big((H-h)^2 e^{-\beta H}\big)\|_{L^1(\rset^{2N})}}{
\|\TR\big(e^{-\beta H}\big)\|_{L^1(\rset^{2N})}}\,,\\
\epsilon_2^2 &= \frac{\|\TR\big(\nabla(H_\Psi-h)\cdot\nabla(H_\Psi-h)e^{-\beta H_\Psi}\big)\|_{L^1(\rset^{2N})}}{
\|\TR\big(e^{-\beta H_\Psi}\big)\|_{L^1(\rset^{2N})}}\,,
\end{split}
\end{equation}
with the definition $H_\Psi(x,p):=\frac{|p|^2}{2}{\rm I} + \Psi^*(x)V(x)\Psi(x)$ using the Hermitian transpose $\Psi^*(x)$\,.

We note that the mean-field Hamiltonian can be written
\begin{equation}\label{h_simple}
\begin{split}
h(x,p)&=\frac{|p|^2}{2}+\frac{\TR(V(x)e^{-\beta V(x)})}{\TR(e^{-\beta V(x)})}\\
&= \frac{|p|^2}{2} + \frac{\sum_{i=0}^{d-1}\lambda_i(x)e^{-\beta\lambda_i(x)}}{\sum_{i=0}^{d-1} e^{-\beta\lambda_i(x)}}=:
\frac{|p|^2}{2}+\lambda_*(x)\,,
\end{split}
\end{equation}
where $\lambda_i(x),\ i=0,\ldots, d-1,$ are the eigenvalues of $V(x)$, and
$\lambda_*:\rset^{N}\to\rset$ is the obtained mean-field potential. 
Therefore the mean-field $h$ is independent of the large mass ratio parameter $M$, so that the dynamics  \eqref{h} is independent of $M$ and consequently the nuclei move a distance of order one in unit time. We see that the mean-field $h=\TR(He^{-\beta H})/\TR(e^{-\beta H})$ is the mean value with respect to the Gibbs density. 
The error term  $\epsilon_1^2$ can be written as the corresponding normed variance
\begin{equation}\label{epsilon_1_square_model}
\epsilon_1^2=  \frac{\|\sum_{i=0}^{d-1}(\lambda_i-\lambda_*)^2e^{-\beta\lambda_i}\|_{L^1(\rset^{N})}}{\|\sum_{i=0}^{d-1} e^{-\beta\lambda_i}\|_{L^1(\rset^N)}}
\end{equation}
and at points $x$ where the eigenvalues $\lambda_i(x)$ are separated
a suitable choice of
$\Psi(x)$ 
is the matrix of eigenvectors to $V(x)$ which implies
\begin{equation}\label{eps_22}
\TR\big(\nabla(H_\Psi-h)\cdot\nabla(H_\Psi-h)e^{-\beta H_\Psi}\big)(x,p)
=\sum_{i=0}^{d-1} \big(
\nabla(\lambda_i(x)-\lambda_*(x))\cdot
\nabla (\lambda_i(x)-\lambda_*(x))
e^{-\beta\lambda_i(x)-\beta|p|^2/2}\big)\,.
\end{equation}
On the other hand to have sets with coinciding eigenvalues is generic in dimension two and higher, see \cite{teller}, and there the matrix $\Psi(x)$ is in general not differentiable. Therefore \eqref{eps_22} will typically not hold everywhere. 
Section \ref{sec_numerics} presents numerical experiments where also the size of the error terms are analyzed
in different settings.

In Section~\ref{Sec2} we review necessary background on Weyl calculus 
and state the main theoretical result, namely the error estimate \eqref{trace_est}, as 
Theorem~\ref{thm_mf}, 
together with the ideas of its proof.
Sections~\ref{sec_thm} and \ref{sec_lemmas} present the proof of the theorem. 
Section~\ref{sec_numerics} presents numerical experiments on the approximation error $\mathfrak{T}_{\mathrm{qm}}(t)-\mathfrak{T}_{\mathrm{md}}$ discussed in the next Section \ref{subsec_numerics}.

\subsection{Numerical comparisons}\label{subsec_numerics}
In Section~\ref{sec_numerics} we present some numerical examples with varying settings of $t,1/M, \epsilon^2_i$, related to different potentials $V$
with the purpose to study the following questions:
Is the estimate \eqref{trace_est} sharp or does the error in practise behave differently with respect to $t$, $1/M$ and $\epsilon_i^2$ ? Is the main contribution to the error coming  from approximation of the the 
matrix-valued potential by a scalar potential in the quantum setting or from the classical approximation of quantum dynamics based on scalar potentials ? Can the mean-field dynamics improve approximation compared to using molecular dynamics based on the ground state eigenvalue $\lambda_0$ instead of $\lambda_*$ ?

\begin{figure}[!htb]
    \begin{center}
  \includegraphics[height=6cm]{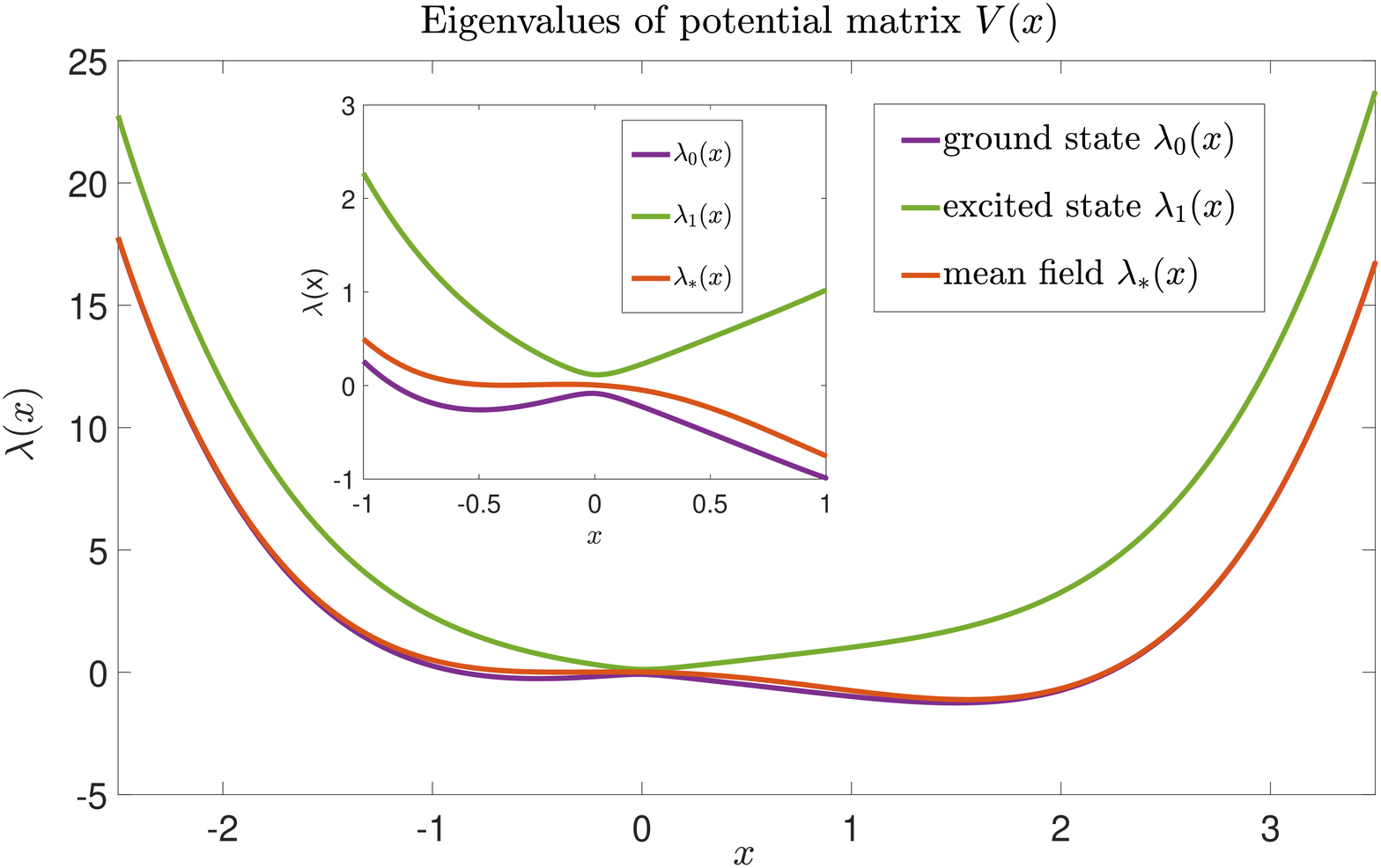}
  \caption{The eigenvalue functions $\lambda_0(x)$ and $\lambda_1(x)$ of the potential matrix $V(x)$, and the corresponding mean-field potential $\lambda_\ast(x)$, for the parameters $c=1$, $\delta=0.1$ (Case E).}\label{fig:eigval_plot_A51}
  \end{center}
\end{figure}

Theorem~\ref{thm_mf}
does not give precise answers to these questions. The aim of 
this section is to
provide some insight from several numerical experiments on a model problem, chosen to avoid the computational difficulties for realistic systems with many particles. Therefore we use one nuclei in dimension one, $N=1$, and two electron states, defined by the Hamiltonian

\begin{equation}\label{H_hat_oper1}
    \widehat{H}=-\frac{1}{2M}\Delta \otimes {\rm I}+V(x),
\end{equation}
where $\Id$ is the $2\times 2$ identity matrix, $V:\mathbb{R}\to\mathbb{R}^{2\times2}$ given by
\begin{equation}\label{pot_mat_V_def1}
V(x)=
\frac{1}{4}(x-\frac{1}{2})^4 {\rm I}+
c\begin{bmatrix}
x & \delta \\
\delta & -x\\
\end{bmatrix},\quad (x,c,\delta)\in\mathbb{R}\times\rset\times\rset\,, 
\end{equation}
with the two eigenvalues 
\begin{equation} \label{two_eigval_new1}
\lambda_0(x)=\frac{1}{4}\,(x-\frac{1}{2})^4-c\,\sqrt{x^2+\delta^2}\,,\quad \lambda_1(x)=\frac{1}{4}\,(x-\frac{1}{2})^4+c\,\sqrt{x^2+\delta^2}\,,
\end{equation}
plotted 
 in Figure~\ref{fig:eigval_plot_A51} (Case E in Table~\ref{table_case_A_to_E_new_written_caption}).

Section~\ref{sec_numerics} presents numerical results comparing quantum mechanics to the three different numerical approximations based on:
the ground state potential $\lambda_0$, mean-field potential $\lambda_*$ and excited state dynamics. The excited state molecular dynamics studied in \cite{KPSS} uses several paths related to different electron eigenvalues and is defined by
\begin{equation}\label{excited_state_formula1}
\mathfrak{T}_{\mathrm{es}}(\tau):=
\sum_{j=0}^{d-1}\int_{\mathbb{R}^{2N}}  A_0(z_\tau^j(z_0)) B_0(z_0) 
\frac{e^{-\beta(\frac{|p_0|^2}{2}+\lambda_j(x_0))}}{
\sum_{k=0}^{d-1}\int_{\mathbb{R}^{2N}} e^{-\beta(\frac{|p|^2}{2}
+\lambda_k(x))}
\mathrm{d}x {\mathrm{d}}p}\mathrm{d}z_0,
\end{equation}
where $z_\tau^j=(x_\tau^j,p_\tau^j)$, $j=0,\ldots, d-1$ solves the Hamiltonian dynamics for state $j$
\[ 
\begin{aligned}
\dot{x}_\tau^j&=p_\tau^j\,,\\
\dot{p}_\tau^j&=-\nabla\lambda_j(x_\tau^j)\,,
\end{aligned}
\]
with the initial condition $z_0^j=(x_0,p_0)=z_0$.

%

\subsubsection{Equilibrium observables}
For equilibrium observables, i.e., where $\tau=0$, the mean-field and excited state
dynamics are equal since then
\begin{equation}\label{equilibrium_tau_0_compare}
\mathfrak{T}_{\mathrm{md}}(0)=
\mathfrak{T}_{\mathrm{es}}(0)= \int_{\rset^{2N}} A_0(z_0)B_0(z_0) 
\frac{\TR(e^{-\beta H(z_0) }){\mathrm{d}}z_0}{\int_{\rset^{2N}}\TR(e^{-\beta H(z) }){\mathrm{d}}z}\,.
\end{equation}
Numerical results on equilibrium observables show that
mean-field and excited state molecular dynamics are more accurate than molecular dynamics based only on the ground state. In the case of correlation observables with $\tau>0$, mean-field and excited state molecular dynamics give in general different approximations.


\subsubsection{Correlation observables}
Observations on quantum dynamics for a system in interaction with a heat bath at thermal equilibrium can be approximated by correlations \eqref{qm_trace}
in the canonical ensemble, cf. \cite{feynman,morandi,ford_kac2,HS}. 
For instance, the classical observable for the diffusion constant 
\[
\frac{1}{6\tau}\frac{3}{N}\sum_{k=1}^{N/3} |x_k(\tau)-x_k(0)|^2=\frac{1}{2N\tau}\big(
|x({\tau})|^2 +|x(0)|^2-2x(\tau)\cdot x(0)\big)
\]
includes the time-correlation $x(\tau)\cdot x(0)$.
Hence the corresponding quantum correlation \eqref{qm_trace} would for this case use
$\OPER{A}_\tau=\widehat x_\tau
\rm I$ and $\OPER{B}_0=\widehat x_0\rm I$, and
\[
\widehat x_\tau\cdot \widehat x_0=\sum_{n=1}^{N/3}\sum_{j=1}^3\EXP{\IU\tau M^{1/2}\widehat H}\widehat x_{n_j} \EXP{-\IU\tau M^{1/2}\widehat H}  \widehat x_{n_j}\PERIOD
\]

The numerical results in Section~\ref{sec_numerics} for time-dependent observables are mainly
based  on the momentum auto-correlation \[\frac{\TR\big((\widehat{p}_\tau \cdot\widehat{p}_0 +\widehat{p}_0\cdot\widehat{p}_\tau) e^{-\beta \widehat H}\big) }{2\TR( e^{-\beta \widehat H})}\,,\]
which is related to the diffusion constant $D$ by the Green-Kubo formula \cite{Hansen}
\[
D=
\frac{\int_0^\infty \int_{\rset^{2N}}p_s(x_0,p_0)\cdot p_0 \TR\big(e^{-\beta H(x_0,p_0)}\big){\mathrm{d}}x_0{\mathrm{d}}p_0{\mathrm{d}}s}{
\int_{\rset^{2N}} \TR\big(e^{-\beta H(x_0,p_0)}\big){\mathrm{d}}x_0{\mathrm{d}}p_0}\,,
\]
since the velocity is equal to the momentum in our case with unit particle mass.

\medskip

The different numerical experiments in Table \ref{table_case_A_to_E_new_written_caption} are chosen by varying the parameters such that all three, two, one or no molecular dynamics approximate well.
In  \hyperlink{case_A}{Case A}, with low temperature and large eigenvalue gap, all three molecular dynamics approximate the quantum observable with similar small error
and also the error terms $1/M$, $\epsilon_1^2$ and $\epsilon_2^2$ are small.
%
%
In \hyperlink{case_B1}{Case B}, 
with small difference of the eigenvalues (i.e. $c$ is small), the mean-field and excited states dynamics is more accurate than ground state dynamics and the error terms $1/M$, $\epsilon_1^2$ and $\epsilon_2^2$ are still small. The result is similar in  \hyperlink{case_C}{Case C}, with an avoided crossing (i.e. $\delta$ is small) and small difference of the eigenvalues. 
In  \hyperlink{case_D}{Case D}, with high temperature and larger difference of the eigenvalues, only the excited state dynamics provides accurate approximations to the quantum observables and the error terms $\epsilon_1^2$ and $\epsilon_2^2$ are large. 
Finally in  \hyperlink{case_E}{Case E}, when the difference of the eigenvalues are sufficiently large and we have an avoided crossing, molecular dynamics is accurate only for short correlation time $\tau$ and the error terms $\epsilon_1^2$ and $\epsilon_2^2$ are large.

\begin{table}[h]
\begin{center}
\begin{tabular}{ |c|c|c|c|c|c|c|c|l| }
\hline
case & $\beta$ & $c$ & $\delta$ & $q_1$ & $\epsilon_1^2$ & $\epsilon_2^2$ & plot & features\\
 \hline
\hyperlink{case_A}{A} & 3.3 & 1 & 1 & 0.0002 & 0.001 & 0.00013 & Fig.\ref{fig:eigval_plot_case_A} & \makecell[l]{ high $\beta$,  large $\gamma_\lambda$, large gap} \\ 
 \hline
 \hyperlink{case_B1}{B} & 1 & 0.1 & 1 & 0.43 & 0.019 & 0.004 & Fig.\ref{fig:eigval_plot_case_B} & \makecell[l]{medium $\beta$, small $\gamma_\lambda$, small gap}
 \\  
 \hline
 \hyperlink{case_C}{C} & 1 & 0.1 & 0.01 & 0.46 & 0.009 & 0.010 & Fig.\ref{fig:eigval_plot_case_C} & \makecell[l]{medium $\beta$, small $\gamma_\lambda$, smallest gap}\\  
 \hline
 \hyperlink{case_D}{D} & 0.28 & 1 & 1 & 0.30 & 2.011 & 0.416 & Fig.\ref{fig:eigval_plot_case_D} & \makecell[l]{low $\beta$, large $\gamma_\lambda$, large gap} \\ 
 \hline
 \hyperlink{case_E}{E} & 1 & 1 & 0.1 & 0.16 & 0.290 & 0.503 & Fig.\ref{fig:eigval_plot_A51} & \makecell[l]{medium $\beta$, large $\gamma_\lambda$, small gap} \\  
 \hline
\end{tabular} 
\end{center}
\caption{A summary of different parameter settings in each test case, where $q_1$ denotes the probability of the electronic excited state, with the precise definition in \eqref{mu_cl_weighted_qj}. Low value of $\beta$ means a high system temperature, the parameter $c$ measures the difference between the two eigenvalues 
$\gamma_\lambda:= 
\int_\R |\lambda_1 - \lambda_0|\, e^{-\beta\lambda_0}\,{\mathrm{d}}x/\int_\R e^{-\beta\lambda_0}\,{\mathrm{d}}x$
and the parameters $c$ and $\delta$ determine the eigenvalue gap. The value of $\epsilon_2^2$ is computed following \eqref{eps_22} with $\Psi(x)$ the matrix of eigenvectors to $V(x)$.}
\label{table_case_A_to_E_new_written_caption}
\end{table}

Figures \ref{fig:err_qm_scalar_qm_case_A6} and \ref{fig:err_case_A5} 
show that in  \hyperlink{case_D}{Case D} and \hyperlink{case_E}{Case E}, where mean-field and ground state approximations are not accurate,
the approximation error is dominated by the part corresponding to
replacing the matrix valued potential by a scalar potential in the quantum formulation and not the part of the error resulting from classical approximation of quantum mechanics for scalar potentials.
More precise conclusions relating the error terms $1/M$, $t\epsilon_1^2$
and $t^2\epsilon_2^2$ in \eqref{trace_est} to the numerical experiments are in Section \ref{sec_numerics}. The numerical experiments here also show that the mean-field dynamics has similar or better accuracy compared to ground state molecular dynamics. It would be interesting to do this comparison for realistic problems. 

\subsection{Relation to previous work}
Classical approximation of canonical quantum correlation observables have been derived  with $\BIGO(M^{-1})$ accuracy for any temperature, see, e.g., \cite{Teufel,KPSS}. Computationally this accuracy requires to solve classical molecular dynamics paths related to several electron eigenvalues, while the mean-field dynamics has the advantage to use only one classical path at the price of loosing accuracy over long time. 


Classical limits of canonical quantum observables were first studied by Wigner, \cite{wigner}. His proof introduces 
the ``Wigner'' function for scalar Schr\"odinger equations and uses an expansion in the Planck constant $\hbar$ to relate equilibrium quantum observables to the corresponding classical  Gibbs phase space averages. 

To derive classical limits in the case of matrix or operator-valued Schr\"odinger equations previous works, see \cite{Teufel}, diagonalize the electron eigenvalue problem, which then excludes settings where the eigenvalues coincide at certain points due to the inherent loss of regularity at such points. The mean-field formulation presented here avoids diagonalization of the electron eigenvalue problem at points with low eigenvalue regularity.

The classical limit with a scalar potential $V$, e.g. the electron ground state eigenvalue, has been studied by three different methods:
\begin{itemize}
\item[(1)] Solutions to the von Neumann quantum Liouville equation, for the density operator, are shown to converge to the solution of the classical Liouville equation, using the Wigner function \cite{Lions-Paul}, also under low regularity of the potential, cf. \cite{Figalli-Paul}. These results use the Wigner function and compactness arguments, which do not provide a convergence rate. 

\item[(2)] In the second method, used in our work,
the two main mathematical tools are a generalized version of Weyl's law and quantization properties, as described by semiclassical analysis, e.g., in \cite{zworski} and \cite{martinez_book}. The generalized Weyl's law links the trace for canonical quantum observables to classical phase space integrals, related to Wigner's work. The quantization properties compare the quantum and classical dynamics and provide convergence rates. Our study differs 
from previous works that also used similar tools, e.g., \cite{Teufel, robert}.
The standard method to bound remainder terms in semiclassical expansions, based on the Planck constant $\hbar$, use the Calder\'on-Vaillancourt theorem to estimate operator norms. Such approach yields error bounds with constants depending on the $L^1(\rset^{2N})$-norms of Fourier transforms of symbols and the potential. 
The $L^1(\rset^{2N})$-norm of a Fourier transformed function
can be bounded by the $L^1(\rset^{2N})$ norm 
of derivatives of order $N$ of the function.
Therefore
the obtained constants in $\BIGO(\hbar^\alpha)$ error estimates are large in high dimension $N$.
In our work here we apply dominated convergence to obtain error estimates based on low regularity of symbols and the potential, while
Fourier transforms in $L^1(\rset^{2N})$ are only required to be finite and do not enter in the final error estimates.
\item[(3)] The third alternative in \cite{Golse-Paul} provides a new method that also avoids the large constants in the Calder\'on-Vaillancourt theorem in high dimensions by using convergence with respect to a generalized  Wasserstein distance and different weak topologies.   
\end{itemize}

A computational bottleneck in ab initio molecular dynamics simulations of canonical correlation observables is in solving electron eigenvalue problems at each time step. 
An alternative to approximate quantum observables is to use path integral Monte Carlo formulations in order to evaluate Hamiltonian exponentials. The Hamiltonian exponentials come in two forms: oscillatory integrals in time $t\in\rset$, 
based on $e^{\IU t \OPERW{H}}$ for the dynamics, and integrals for Gibbs function $e^{-\beta\OPERW{H}}$ that decay with increasing inverse temperature $\beta\in(0,\infty)$.
The high variance related to the oscillatory integrand $e^{\IU t \OPERW{H}}$ means that standard computational path integral formulations for molecular dynamics are applied only to the statistics based on the partition function 
$\TR (e^{-\beta\OPERW{H}})$ while the
dynamics is approximated classically. 

Two popular path integral methods are centroid molecular dynamics and
ring-polymer molecular dynamics, see, e.g.,  \cite{ring-poly} or \cite{tuckerman_compare}. 
In these methods the discretized path integral is interpreted as a classical Hamiltonian with a particle/bead for each
degree of freedom in the discretized path integral.  
For the centroid version the dynamics is 
based on the average of the particle/bead positions, i.e., the centroid, with forces related to a free energy potential for the partition function thereby forming a mean-field approximation.
It is related to the mean-field approximation \eqref{md_mf}  and \eqref{mf-def} 
but differs in that in our work  the forces are based on the mean Hamiltonian, 
for the partial trace over electron degrees of freedom, instead of on the partition function for the discretized path integral with respect to both nuclei and electron degrees of freedom, centered at the centroid. 
In ring-polymer molecular dynamics classical kinetic energy is added for each bead forming a Hamiltonian
with harmonic oscillators in addition to the original potential energy. Consequently the phase-space is related
to coupled ring polymers, one for each original particle. 

There is so far no convergence proof for centroid nor ring-polymer molecular dynamics.
Therefore it would be interesting to further study their relation to the mean-field model we analyse here.
The mean-field formulation \eqref{md_mf} and \eqref{mf-def} can also offer an alternative
to standard eigenvalue solutions by using a path integral formulation of the partial trace 
over the electron degrees of freedom, in the case of sufficiently large temperature 
avoiding the fermion sign problem,  see, e.g., \cite{dornheim}.
Another difference to previous work is that the convergence proof here derives a 
precise asymptotic representation of the Weyl symbol for
the Gibbs density operator using a path integral formulation, providing an example
that  using path integrals for the Gibbs density can also result in simplification 
of the theory.

\section{The main result and background from Weyl calculus}\label{Sec2}

In this section we state the main theorem and review necessary tools from semiclassical analysis and functional integration.

\medskip
To relate the quantum and classical observables
we employ Weyl calculus for matrix-valued symbols. 
First,
we introduce functional spaces that we use in the sequel:

(i) the Schwartz space of matrix-valued functions on the phase space
\[
\SSPACE := \big\{ A\in C^\infty(\rset^N\times\rset^N,\cset^{d\times d})\, \SEP 
   \sup_{z\in\R^N\times\R^N}|z^\gamma
   \partial^\alpha_z A_{ij}(z)| < \infty\,,\mbox{ for all indices}\big\}
\]
where we denote a point in the phase space $z=(x,p)$
and for a multi-index of non-negative integers $\alpha=(\alpha_1,\dots,\alpha_{2N})$, we have the partial derivatives $\partial^{\alpha}_z = \partial^{\alpha_1}_{z_1}\dots\partial^{\alpha_{2N}}_{z_{2N}}$ of the order  $|\alpha|=\sum_i\alpha_i$ and similarly we have $z^\gamma = z_1^{\gamma_1}\dots z_{2N}^{\gamma_{2N}}$ for the multi-index $\gamma$. 
For a matrix-valued symbol $A$ we also use the notation 
$A'(z)$ for the tensor $(A'(z))^m_{ij}:= 
\partial_{z_m} A_{ij}(z)$ and $A''(z)$ for 
the 4th-order tensor  $(A''(z))^{mn}_{ij}:= \partial^2_{z_m z_n} A_{ij}(z)$. The dual space of tempered distributions is denoted $\SSPACE'$. 

(ii) the space of $L^2$ vector-valued wave functions
\[
\HSPACE := L^2(\R^N,\cset^d)\equiv [L^2(\R^N)]^d\,.
\]

We define the Weyl quantization operator of a matrix-valued symbol $A\in\SSPACE$ 
as the mapping $A\mapsto \OPER{A}$
that assigns to the symbol $A$ the linear operator 
$\OPER{A}:\HSPACE \to \HSPACE$ 
defined  for all Schwartz functions $\phi(x)$ by
\begin{equation}\label{Weyl:definition}
\begin{split}
  \OPER{A} \phi(x)=&
  \int_{\rset^{N}} \NORMFAC^N\int_{\rset^{N}} \EXP{\IU M^{1/2}(x-y)\cdot p} 
   A(\XYH,p)  \Rd p\, \phi(y) \Rd y \\
   =:&\int_{\R^N} K_A(x,y) \phi(y)\,\Rd y\,,
\end{split}
\end{equation}
and extended to all wave functions $\phi\in \HSPACE$ by density.
The expression \eqref{Weyl:definition} shows that the kernel $K_A$ on $\mathcal H$ is the Fourier transform in the second argument of the symbol $A(x,p)$ and consequently
the Weyl quantization is well defined for
symbols in $\SSPACE'$, the space of tempered distributions. 

For example, the symbol 
\[
H(x,p):=\tfrac{1}{2}|p|^2 \Id + V(x)
\] 
yields the Hamiltonian operator 
\[
\OPER{H} = -\frac{1}{2M}\Delta \otimes \Id + V(x)\,.
\]

We formulate the main result as the following theorem estimating the mean-field approximation.
\begin{theorem}\label{thm_mf}
Let $\Psi:\rset^N\to \mathbb C^{d\times d}$ be a differentiable mapping into orthogonal matrices and define $V_\Psi:=\Psi^*V\Psi$,  for the Hermitian potential $V:\rset^{N}\to \cset^{d\times d}$. Assume that 
the components of the Hessian $V_\Psi''$ 
are in the Schwartz class 
and the scalar symbols $a_0:\rset^{2N}\to \mathbb C$ and $b_0:\rset^{2N}\to \mathbb C$
are infinitely  differentiable and compactly supported. 
Furthermore, suppose that  there is a constant $k$ such that $V+k\Id$ 
is positive definite everywhere, $\TR(e^{-\beta\OPER{H}})$ 
is finite, 
and there is a constant $C$ such that 
\[
\begin{split}
\|\TR e^{-\beta V}\|_{L^1(\rset^N)} + \|\TR( V^2e^{-\beta V})\|_{L^1(\rset^N)} +
\|\TR\big( (\sum_{m=1}^N\sum_{|\alpha|\le 3} \partial^\alpha_{x_m} V_\Psi \partial^\alpha_{x_m} V_\Psi)e^{-\beta V_\Psi}\big)\|_{L^1(\rset^N)}&\le C\,,\\
\sum_{|\alpha|\le 3}\big(\|\partial^\alpha_z a_0\|_{L^\infty(\rset^{2N})}
+\|\partial^\alpha_z b_0\|_{L^\infty(\rset^{2N})}\big)&\le C\,,\\
\max_i\sum_{|\alpha|\le 3}\|\partial_x^\alpha \partial_{x_i}h\|_{L^\infty(\rset^{2N})} &\le C\,,\\
\int_{\rset^{2N}}\TR\big(|\sum_{n,m}p_n\partial^2_{x_n x_m}V_\Psi(x)p_m|\,e^{-\beta (\frac{|p|^2}{2}{\rm I} + V_\Psi(x))}\big){\mathrm{d}}x{\mathrm{d}}p &\le C\,,\\
\int_{\rset^{2N}}\TR\big((\sum_{n}p_n\partial^2_{x_n}V_\Psi(x))^2\,e^{-\beta (\frac{|p|^2}{2}{\rm I}+V_\Psi(x))}\big){\mathrm{d}}x{\mathrm{d}}p &\le C\,,
\end{split}
\]
then there are constants $c, \bar M, \bar t$, depending on $C$ and $\beta$, such that the quantum canonical observable \eqref{qms}, with $A_0=a_0\Id$ and $B_0=b_0\Id$,
can be approximated by the mean-field molecular dynamics \eqref{md_mf}-\eqref{h} with the error
 \begin{equation}\label{mf_error}|\mathfrak{T}_{\mathrm{qm}}(t)-\mathfrak{T}_{\mathrm{md}}(t)|\le c(M^{-1}+t\epsilon_1^2+ t^2\epsilon_2^2)\,,
 \end{equation}
 for $M\ge \bar M $ and $0\le t\le \bar t$, where
 \[
 \begin{split}
\epsilon_1^2 &= \|\TR\big((H-h)^2 e^{-\beta H}\big)\|_{L^1(\rset^{2N})}/
\|\TR\big(e^{-\beta H}\big)\|_{L^1(\rset^{2N})}\,,\\
\epsilon_2^2 &= \|\TR\big(\nabla(H_\Psi-h)\cdot\nabla(H_\Psi-h)e^{-\beta H_\Psi}\big)\|_{L^1(\rset^{2N})}/
\|\TR\big(e^{-\beta H}\big)\|_{L^1(\rset^{2N})}\,, \quad H_\Psi:=\frac{|p|^2}{2}{\rm I} + V_\Psi\,.
 \end{split}
 \]
\end{theorem}


\subsection{Overview and background to the proof}
This section provides background and motivation to the proof of the theorem in three subsections.
The first subsection reviews application of Weyl calculus for the dynamics, the
second one is on a generalized form of Weyl's law in order to relate the
quantum trace to phase space integrals
and the third subsection introduces path integrals and their application in the context
of our result.

\subsubsection{Weyl calculus and dynamics}
This section first introduces the central relation between commutators and corresponding Poisson brackets for the classical limit of dynamics. Given two smooth functions $v(x,p)$, $w(x,p)$ on the phase space we define the Poisson bracket
\begin{equation}\label{eqn:Poissonbra}
\begin{split}
\{v,w\}&:=\nabla_p v(x,p)\cdot\nabla_x w(x,p)-\nabla_x v(x,p)\cdot\nabla_p w(x,p)\\
&=\big(\nabla_{p'},-\nabla_{x'})\cdot(\nabla_{x},\nabla_{p})\big)v(x',p')w(x,p)\big|_{(x,p)=(x',p')}\,.
\end{split}
\end{equation}
We denote the gradient operator in the variable $z = (x,p)$ as
$\nabla_z = (\nabla_x,\nabla_p)$ and $\NABLAR_z = (\nabla_p,-\nabla_x)$, 
hence the Poisson bracket is expressed as
\begin{equation}\label{eqn:Poissonbrabis}
\{v,w\} = \big(\NABLAR_{z'}\cdot\nabla_{z}\big)v(z')w(z)\big|_{z=z'}\,.
\end{equation}
For two operators $\OPER{C}$, $\OPER{D}$ on the space $\mathcal{H}$ we define their commutator
\[
[\OPER{C},\OPER{D}] = \OPER{C}\OPER{D} - \OPER{D}\OPER{C}\,.
\]

To relate the quantum and classical dynamics for particular observables with symbols of the type $a_0\Id$ treated in Theorem~\ref{thm_mf}, 
we assume that the classical Hamiltonian
flow $z_t(z_0):=\big(x_t(z_0),p_t(z_0)\big)$ solves the Hamiltonian system
\[
\begin{split}
\DT{x}_t &=p_t=\nabla_p h(x_t,p_t)\,,\\
\DT{p}_t &=-\nabla_x h(x_t,p_t)\,,\\
\end{split}
\]
with the initial data $(x_0,p_0)=z_0\in\rset^{2N}$.
Given a scalar  Schwartz function 
$a_0$ we define a smooth function 
on the flow $z_t(z_0)$ 
\begin{equation}\label{eq:at}
a_t(z_0):=a_0\big(z_t(z_0)\big)\,,
\end{equation}
and we have that $a_t(z_0)$ satisfies
\begin{equation}\label{ce}
\partial_t a_t(z_0)=\{h(z_0),a_t(z_0)\}
\end{equation}
since the direct calculation gives
\[
\begin{split}
\partial_t a_t(z_0)=\partial_s a_0\big(z_t(z_s)\big)\big|_{s=0}
&=\DT{z}_0\cdot\nabla_{z_0}a_0\big(z_t(z_0)\big)=\{h(z_0),a_t(z_0)\}\,.
\end{split}
\]
The corresponding quantum evolution of the observable $\OPER{A}$, for $t\in \rset$, 
is defined by the Heisenberg-von Neumann equation
\begin{equation}\label{qme}
\begin{split}
\partial_t \OPER{A}_t&= {\IU}M^{1/2} [\OPERW{H},\OPER{A}_t]\,,
\end{split}
\end{equation}
which implies the representation 
\[
\OPER{A}_t=e^{{\IU}t M^{1/2}\OPER{H}}\OPER{A}_0e^{-{\IU}t M^{1/2}\OPER{H}}\,.
\]
The basic property in Weyl calculus that links the quantum evolution \eqref{qme} 
to the classical \eqref{ce} is the relation
\begin{equation}\label{qe-ce}
{\IU}M^{1/2} [\OPERW{H},\OPER{a}_t] = \{H,a_t\}^{\widehat{}} + \OPER{r}_t^a
\end{equation}
where the remainder symbol $r^a_t$ is small. In Lemma~\ref{lemma1} we show that
\begin{equation}\label{Mra}
\lim_{M\to\infty} Mr^a_s=\frac{1}{12}(\nabla_{z_0}\cdot\nabla_{z_0'}')^3H(z_0)a_0\big(z_s(z_0')\big)\big|_{z_0=z_0'}\,.
\end{equation}

The main obstacle to establish the classical limit for dynamics based on the matrix-valued operator  Hamiltonian $\OPER{H}$ is that matrix symbols do not commute, i.e. $[H,A_t]\ne 0$, which implies 
the additional larger remainder 
${\IU}M^{1/2} [{H},{A}_t]$ 
in \eqref{qe-ce}. 
Therefore the usual semiclassical analysis perform approximate diagonalization of $ [\OPERW{H},\OPER{A}_t]$, see \cite{Teufel,KPSS}. Diagonalization of $V$ introduces  eigenvectors that are not smooth everywhere unless the eigenvalues are separated, due to the inherent loss of regularity for eigenvectors corresponding to coinciding eigenvalues, see \cite{teller}.
To have a conical intersection point with coinciding eigenvalues is  generic in dimension two, and in higher dimensions the intersection is typically a co-dimension two set, see \cite{teller}.  The non smooth diagonalization has been so far
a difficult obstacle to handle with the tools of Weyl calculus. 
The aim in our work  is therefore to avoid diagonalization  everywhere
by analyzing  a mean-field approximation differently. 

In order to apply \eqref{qe-ce} we use Duhamel's principle, see \cite{evans}:
the inhomogeneous linear equation in the variable $A_t-a_t\equiv A_t-a_t\Id$, where $A_t$ satisfies the evolution \eqref{qme} and $a_t$ is defined by \eqref{eq:at},
\[
\begin{split}
\partial_t (\OPER{A}_t-\OPER{a_t})&= {\IU}M^{1/2} [\OPERW{H},\OPER{A}_t]- \{h,a_t\}^{\widehat{}}\\
&={\IU}M^{1/2} [\OPERW{H},\OPER{A}_t-\OPER{a}_t] 
+{\IU}M^{1/2} [\OPERW{H},\OPER{a}_t] - \{h,a_t\}^{\widehat{}}\\
\end{split}
\]
can be solved by
integrating solutions to the homogeneous  problem with respect to the inhomogeneity
\begin{equation}\label{duhamel_1}
\OPER{A}_t -\OPER{a}_t = 
\int_0^t  e^{{\IU}M^{1/2}(t-s) \OPERW{H}}\big({\IU}M^{1/2}
[\OPERW{H},\OPER{ a}_s]- \{ h, a_s\}^{\widehat{}}\ \big) e^{{-\rm i}M^{1/2}(t-s) \OPERW{H}}
\Rd s\,.\\
\end{equation}

The quantum statistics has a similar remainder term to \eqref{Mra}, namely the difference $\rho-e^{-\beta H}$
of the Weyl symbol $\rho$ for the quantum Gibbs density operator $\OPER{\rho}=e^{-\beta\OPER{H}}$ and the classical Gibbs density $e^{-\beta H}$.
To characterize asymptotic behaviour as $M\to\infty$ of this difference we employ 
representation of the symbol $\rho$ based on
Feynman-Kac path integral formulation, as presented in Section~\ref{sec:FK}.

\medskip
\subsubsection{Generalized Weyl law}
To link the quantum trace to a classical phase space integral we will use a generalized form
of Weyl's law, see \cite{zworski,Teufel}.
The semiclassical analysis is based on the fact that
the 
$\HSPACE$-trace of a Weyl operator, with a $d\times d$ 
matrix-valued symbol,  is equal to
the phase-space average of its symbol trace. 
Indeed  we have by the definition of the integral kernel in \eqref{Weyl:definition} for $A\in\SSPACE$
\begin{equation}\label{trace1}
\TR{\OPER{A}} =
\int_{\mathbb R^{N}} \TR K_A(x,x) \Rd x = \NORMFAC^N\int_{\mathbb R^{N}}
\int_{\mathbb R^{N}}
\TR  A(x,p)  \Rd p\, \Rd x\PERIOD
\end{equation}
In fact also the composition of two Weyl operators is determined by the phase-space average as follows, see \cite{Teufel}.

\begin{lemma}\label{lemma_composition}
The composition of  two Weyl operators $\OPER A$ and $\OPER B$, with $A\in\SSPACE$ and $B\in\SSPACE$ satisfies
\[
\TR(\OPER{A}\OPER{B})=\NORMFAC^N\int_{\mathbb R^{2N}}  \TR\big(A(x,p)B(x,p)\big) \Rd x\Rd p\COMMA
\]
where $A(x,p)B(x,p)$ is the matrix product of the two $d\times d$ matrices $A(x,p)$ and $B(x,p)$.
\end{lemma}
\begin{proof}
The well-known proof is a straight forward evaluation of the integrals involved in the composition of two kernels and it is given here 
for completeness.
The kernel of the composition is
\[
\begin{split}
   K_{AB}(x,y) &=
    \NORMFAC^{2N}\int_{\mathbb R^{3N}}
      A(\XXH{x}{x'},p)B(\XXH{x'}{y},p')\\
      &\qquad\times \EXP{iM^{1/2}\big((x-x')\cdot p + (x'-y)\cdot p'\big)} \Rd p' \Rd p \Rd x'
\end{split}
\]
so that the trace of the composition becomes
\[
\begin{split}
    &\TR(\OPER{A}\OPER{B}) =\int_{\mathbb R^N}\TR K_{AB}(x,x) \Rd x\\
    &=\NORMFAC^{2N}\int_{\mathbb R^{4N}}
       \TR\big(A(\XXH{x}{x'},p)B(\XXH{x'}{x},p')\big) \\
&\qquad \times \EXP{iM^{1/2}\big((x-x')\cdot p + (x'-x)\cdot p'\big)} \Rd p' \Rd p \Rd x' \Rd x\\
&=\NORMFAC^{2N}\int_{\mathbb R^{4N}}
\TR\big(A(y,p)B(y,p')\big) \EXP{iM^{1/2}y'\cdot (p-p')} \Rd p' \Rd p \Rd y' \Rd y\\
&=\NORMFAC^{N}\int_{\mathbb R^{2N}}
\TR\big(A(y,p)B(y,p)\big)  \Rd p  \Rd y\COMMA\\
\end{split}
\]
using the change of variables $(y,y')=\big((x+x')/2,x-x'\big)$, which verifies the claim.
\end{proof}
The composition of three operators does not have a corresponding phase-space representation. We will instead
 use the composition operator $\MP$ (Moyal product), defined by
$\OPER{A}\OPER{B}=\widehat{A\MP B}$, to reduce the number of Weyl quantizations to two, e.g., as  $\OPER{A}\OPER{B}\OPER{C}\OPER{ D}=\widehat{A\MP B}\,\widehat{C\MP  D}$. 
More precisely, the Moyal product of two symbols has the representation
\begin{equation}\label{composition_bis}
A\MP B = e^{\frac{{\IU}}{2M^{1/2}} (\nabla_{x'}\cdot\nabla_{p}- \nabla_{x}\cdot\nabla_{p'})}
A(x,p) B(x',p')\big|_{(x,p)=(x',p')}\,.
\end{equation}
For general background we refer the reader to \cite{zworski} or \cite{martinez_book}.

The isometry between  Weyl operators with the Hilbert-Schmidt  inner product,  $\TR({\OPER{A}^*\OPER{B}})$, and the corresponding $L^2(\rset^N\times\rset^N,\mathbb C^{d\times d})$ symbols obtained by  Lemma~\ref{lemma_composition} shows how to extend from symbols in $\SSPACE$ to $L^2(\rset^N\times\rset^N,\mathbb C^{d\times d})$ by density of $\SSPACE$
in $L^2(\rset^N\times\rset^N,\mathbb C^{d\times d})$, see \cite{Teufel}. We will use the Hilbert-Schmidt norms 
$\|\OPER{A}\|_{\mathcal{HS}}^2=\TR(\OPER{A}^*\OPER{A})=\TR(\OPER{A}^2)$ and $\|\TR A^2\|_{L^2(\rset^{2N})}$
to estimate Weyl operators and Weyl symbols, respectively.

We show in Lemma~\ref{lemma1} that having the Weyl symbols and $V''$  in the Schwartz class
imply that dominated convergence can be applied in the phase space integrals 
obtained from the generalized form of Weyl's law.

\medskip
\subsubsection{Feynman-Kac path integrals}\label{sec:FK}
In order to analyse the symbol $\rho$ for the Gibbs density operator $\OPER{\rho}=e^{-\beta \OPERW{H}}$ we will use path integrals (in so called imaginary time), as in \cite{feynman}, \cite{amour} and \cite{papanicolaou}, based on Feynman-Kac formula applied to the kernel of the Gibbs density operator and its corresponding Weyl quantization. 
We start with the kernel representation
\[
(e^{-\beta \OPERW{H}}\phi)(x') = \int_{\rset^N} K_\rho(x',y')\phi(y')\Rd y'\,,
\]
obtained from that $(e^{-\beta \OPERW{H}}\phi)(x')=:u(x',\beta)$ solves
the parabolic partial differential equation
\[
\partial_\beta u(\cdot,\beta)+\OPERW{H}u(\cdot,\beta)=0\,,\quad \beta>0\,, \ u(\cdot,0)=\phi\,.
\]
To motivate the construction of the path integral representation 
we first identify the Weyl symbol  of the density operator in the case of a scalar potential $V$, i.e., the case $d=1$. In order to emphasize that we consider the scalar case, we denote this Weyl symbol $\rho_s(x,p)$, and thereby the
associated  kernel is denoted $K_{\rho_s}$. 
Direct application of
the Feynman-Kac formula, see Theorem 7.6 in \cite{karatzas}, implies that the kernel can be written as the expected value
\begin{equation}\label{4.2}
K_{\rho_s}(x',y') = \E\Big[e^{-\int_0^\beta V(\omega'_r)\Rd r} \delta(\omega'_\beta-y')\SEP \omega'_0=x'\Big]\,,
\end{equation}
where $\omega'_t$ solves the stochastic differential equation
\[
\Rd \omega'_t=M^{-1/2}\Rd W_t\,,\;\;\;\;\;\omega'_0 = x'\,,
\]
with the standard Wiener process $W_t$ in $\rset^N$ and the delta measure $\delta(\omega'_\beta-y')$ concentrated at the point $y'$. 

We recall the definition \eqref{Weyl:definition} of Weyl quantization for the 
(scalar) symbol $\rho_s(x,p)$
\[
\OPER{\rho}_s\phi(x) = \int_{\rset^{N}}\underbrace{\NORMFAC^N\int_{\rset^{N}}
e^{{\IU}M^{1/2}(x- y)\cdot p}
 \rho_s(\HALF(x+y), p) \Rd p}_{=K_{\rho_s}(x,y)}\, \phi(y)\Rd y\,,
\]
from which we obtain the expression for the symbol of an operator
associated with the kernel $K_{\rho_s}(x,y)$
\begin{equation}\label{4.1}
\rho_s(x,p)=\int_{\rset^{N}} e^{-{\IU}M^{1/2} y\cdot p} K_{\rho_s}(x+\frac{y}{2}, x-\frac{y}{2}) \Rd y\,.
\end{equation}
Using the substitution $x'=x+\tfrac{y}{2}$, and $y'=x-\frac{y}{2}$, i.e.,
\[
y=x'-y' \quad \mbox{and} \quad x=\frac{x'+y'}{2}\,,
\]
and letting
\[
\omega'_t=x+\frac{y}{2}-\omega_t
\]
where 
\begin{equation}\label{eqn:sde}
\Rd \omega_t=M^{-1/2}\Rd  W_t\,,\;\;\;\;\;\;
\mbox{$\omega_0=0$,}
\end{equation}
imply by combining \eqref{4.1}, \eqref{4.2} and the transformed path \eqref{eqn:sde}
\[
\begin{split}
\rho_s(x,p)&=\int_{\rset^N}e^{-M^{1/2}{\IU} y\cdot p}
\mathbb E[ e^{-\int_0^\beta V(x+\frac{y}{2}-\omega_r){\mathrm{d}}r} \delta(y-\omega_\beta)\, |\, \omega_0=0]{\mathrm{d}}y\\
&=\E[e^{-{\IU}M^{1/2}\omega_\beta\cdot p}
e^{-\int_0^\beta V(x+\frac{1}{2}\omega_\beta-\omega_r)\Rd r}\, |\, \omega_0=0]\,.
\end{split}
\]
We have obtained the path integral representation of the symbol corresponding
to the Gibbs density operator $e^{-\beta\OPER{H}}$ in the case of the scalar potential $V$
\begin{equation}\label{4.3}
\begin{split}
\rho_s(x,p)
&=\E[e^{-{\IU} W_\beta\cdot p}
e^{-\int_0^\beta V(x+M^{-1/2}(\frac{1}{2}W_\beta-W_r))\Rd r}] \,,
\end{split}
\end{equation}
as derived in \cite{amour}.

\medskip

We can proceed along similar lines in the case of the {\it matrix-valued} potential $V$ studied here. The Feynman-Kac formula has been derived for operator-valued potentials $V$, see \cite{papanicolaou}, and
in the case where the potential $V(x)$ is a general  matrix, the exponential 
\[
e^{-\int_0^\beta V(x+M^{-1/2}(\frac{1}{2}W_\beta-W_r))\Rd r}\,,
\] 
in \eqref{4.3}
for the time evolution, is replaced by the
corresponding matrix-valued process  
$\BARYT_\beta^+ \in \R^{d\times d}$ which solves
\begin{equation}\label{4.4}
\DT{\BARYT}_t^+= -\BARYT_t^+ \underbrace{V\big(x+M^{-1/2}(\tfrac{1}{2}W_\beta-W_t)\big)}_{=:V^+_t}\,, \;\;\;\;\mbox{$t\in (0,\beta)$ and $\BARYT_0^+={\Id}$.}
\end{equation}
We will use the notation ${W^\beta_t}:=\tfrac{1}{2}W_\beta-W_t$.
The steps in the derivation of  \eqref{4.3} then imply that the symbol in the case of a matrix-valued potential can  again be expressed by
\begin{equation*}\label{matFeynman-Kac1}
\rho(x,p)=\E[e^{-{\IU}W_\beta\cdot p} \, \BARYT^+_\beta(W)\, |\, W_0=0]\,.
\end{equation*}
To estimate $\rho-e^{-\beta {H}}$ we use
the symmetry property that
the Weyl symbol, $\rho$, for the Gibbs density operator $\OPER{\rho}=e^{-\beta\OPER{H}}$
is a Hermitian matrix.
Indeed, we have $e^{-\beta \OPERW{H}}$ represented as an $L^2$-integral operator with the kernel $K_\rho$ and since \[\OPERW{H}=-\frac{1}{2M}\Delta\otimes{\Id}+V(x)\]
is real and Hermitian  also 
$e^{-\beta\OPERW{H}}=\sum_{n=0}^\infty (-\beta \OPERW{H})^n/n!$ is real and Hermitian.
Therefore the Weyl symbol corresponding to the Gibbs density operator $\OPER{\rho}$ satisfies
\begin{equation}\label{eq:intrep}
\rho(x,p) = \int_{\rset^N}e^{-{\IU} M^{1/2}y\cdot p}  K_\rho(x+\frac{y}{2}, x-\frac{y}{2}) \Rd y
= \int_{\rset^N}e^{-{\IU} M^{1/2}y\cdot p}  K_\rho(x-\frac{y}{2}, x+\frac{y}{2}) \Rd y\,.\\
\end{equation}
Either of these integral representations show that $\rho$ is Hermitian. 
The same steps as above leading to \eqref{4.3} can by \eqref{eq:intrep} be applied to the change of variables $x'= x-y/2$ and $y'=x+y/2$, which implies that we  also have
\begin{equation*}\label{matFeynman-Kac2}
\rho(x,p)=\E[e^{-{\IU}W_\beta\cdot p} \, \BARYT_\beta^-(W)\, |\, W_0=0]\,,
\end{equation*}
where
\begin{equation}\label{4.4-}
\DT{\BARYT}_t^-= -\BARYT_t^- \underbrace{V\big(x-M^{-1/2}(\tfrac{1}{2}W_\beta-W_t)}_{=:V^-_t}\big)\,, \;\;\;\;\mbox{$t\in (0,\beta)$ and $\BARYT_0^-={\Id}$.}
\end{equation}
Therefore we have the symmetrized representation
\begin{equation}\label{matFeynman-Kac}
\rho(x,p)= \frac{1}{2}\E[e^{-{\IU}W_\beta\cdot p} \, \big(\BARYT_\beta^+(W) +\BARYT_\beta^-(W)\big)\, |\, W_0=0]\,.
\end{equation}
 
Using the path-integral representation of the symbol $\rho(x,p)$ we 
prove in Section~\ref{sec_path-proof}
\begin{lemma}\label{path_lemma}
Assume that the bounds in Theorem \ref{thm_mf} hold. Then
\begin{equation}\label{symb_estim}
\begin{split}
\lim_{M\to\infty} &M\big(\rho(x,p)-{e^{-\beta H}}(x,p)\big) \\
&= e^{-\beta |p|^2/2}\Big(
\int_0^\beta(\frac{\beta}{2}-t)^2 e^{-tV(x)}p^TV''(x)p \,e^{-(\beta-t)V(x)}\Rd t\\
&\quad+ \int_0^\beta\int_0^t (\frac{\beta}{2}-s)(\frac{\beta}{2}-t) e^{-sV(x)}V'(x)pe^{-(t-s)V(x)} V'(x)p
e^{-(\beta-t)V(x)}{\mathrm{d}}s{\mathrm{d}}t\Big)\,,
\end{split}
\end{equation}
and
\begin{equation}\label{symb_estim_bis}
\begin{split}
\|\TR\rho^2\|_{L^\infty(\rset^{2N})}&=\BIGO (1)\,,\\
\|\TR\rho^2\|_{L^2(\rset^{2N})}&<\infty\,.
\end{split}
\end{equation}
\end{lemma}

\section{Proof of Theorem \ref{thm_mf}}\label{sec_thm}
The section proves Theorem \ref{thm_mf} based on three steps:
\begin{itemize}
\item[{ Step 1.}] use Duhamel's principle recursively to analyse the dynamics based on $H$,
\item[{ Step 2.}] use estimates of remainders for Weyl compositions and the Weyl Gibbs density to analyse the statistics at $t=0$,
\item[{ Step 3.}] repeat Step 1 and Step 2 with $H$ replaced by $\bar H:= \Psi^* \# H \# \Psi$ to approximately  diagonalize $H$. 
\end{itemize}
\begin{proof}[Proof of the theorem]

{\bf Step 1.}
Lemma \ref{lemma1} shows that
the commutator has the representation
\begin{equation}\label{Hy}
{\IU} M^{1/2} [\OPERW{H},\OPER{a}_s]=(\{H,a_s\}+r^a_s)^{\widehat{}} \,, 
\end{equation}
where the remainder $r^a$ vanishes
as $M\to\infty$.
By Duhamel's principle we have by \eqref{ce} and \eqref{qme} as in \eqref{duhamel_1}
\begin{equation}\label{first_step}
\begin{split}
&\TR\big((\OPER{A}_t\OPER{B}_0e^{-\beta\OPERW{H}} -\OPER{a}_t\OPER{B}_0e^{-\beta\OPERW{H}}) \\
&=\int_0^t \TR\Big( e^{{\IU}M^{1/2}(t-s) \OPERW{H}}\big({\IU}M^{1/2}
[\OPERW{H},\OPER{ a}_s]- \{ h, a_s\}^{\widehat{}}\ \big) e^{{-\rm i}M^{1/2}(t-s) \OPERW{H}}\OPER{B}_0 e^{-\beta\OPERW{H}}\Big)\Rd s\\
&=\int_0^t \TR\Big( e^{{\IU}M^{1/2}(t-s) \OPERW{H}}\underbrace{\big(
\{ H,{ a}_s\}^{\widehat{}}+\hat{r}^a - \{ h, a_s\}^{\widehat{}} \ \big)}_{=:\widehat{Da_s}} e^{{-\rm i}M^{1/2}(t-s) \OPERW{H}}\OPER{B}_0 e^{-\beta\OPERW{H}}\Big)\Rd s\\
\end{split}
\end{equation}
and the cyclic invariance of the trace together with 
$[e^{{\IU}M^{1/2}(t-s) \OPERW{H}},e^{-\beta \OPERW{H}}]=0$
imply
\begin{equation}\label{duhamel_2}
\begin{split}
&\TR\big((\widehat{A}_t\OPER{B}_0e^{-\beta\OPERW{H}} -\OPER{a}_t\OPER{B}_0e^{-\beta\OPERW{H}}) \\
     &=\int_0^t \TR\Big( \widehat{Da_s} e^{{-\rm i}M^{1/2}(t-s) \OPERW{H}}\OPER{B}_0 e^{{\IU}M^{1/2}(t-s) \OPERW{H}}
     e^{-\beta\OPERW{H}}\Big)\Rd s\\
&=\int_0^t \TR\Big( \widehat{Da_s} \OPER{B}_{s-t} e^{-\beta\OPERW{H}}\Big)\Rd s\,.\\
\end{split}
\end{equation}
The right hand side can again be estimated by  applying Duhamel's principle \eqref{duhamel_2}, now to $\OPER{B}_{s-t}$ and $b_{s-t}$,
as follows
\begin{equation}\label{second_step}
\begin{split}
&\int_0^t \TR(\widehat{Da_s} \OPER{B}_{s-t} e^{-\beta\OPERW{H}})\Rd s
-\underbrace{\int_0^t \TR(\widehat{Da_s} \OPER{b}_{s-t} e^{-\beta\OPERW{H}})\Rd s}_{=:T_0} \\
&=
\int_0^t\int_0^{s-t} \TR(e^{-\beta\OPERW{H}}\widehat{Da_s} e^{{\IU}M^{1/2}(s-t-\tau) \OPERW{H}}\widehat{D b_\tau}
e^{-{\IU}M^{1/2}(s-t-\tau) \OPERW{H}}\Rd \tau\Rd s
=:T_3\,.
\end{split}
\end{equation}
We have by Cauchy's inequality and the cyclic invariance of the trace
\begin{equation}\label{forth_step}
\begin{split}
|T_3| &\le \int_0^t\int_0^{t-s}\Big(\TR\big((\widehat{Da_s} )^*\widehat{Da_s} e^{-\beta\OPERW{H}}\big)
\times \\ &\quad\quad\times
\TR\big((e^{{\IU}M^{1/2}(s-t-\tau) \OPERW{H}}\widehat{Db_\tau}e^{-\beta\OPERW{H}/2}e^{-{\IU}M^{1/2}(s-t-\tau) \OPERW{H}})^*
\times \\ &\quad\quad\times
(e^{{\IU}M^{1/2}(s-t-\tau) \OPERW{H}}\widehat{Db_\tau}e^{-\beta\OPERW{H}/2}
e^{-{\IU}M^{1/2}(s-t-\tau) \OPERW{H}})\big)\Big)^{1/2} \Rd \tau\Rd s\\
&=\int_0^t\int_0^{t-s}\Big(\TR\big((\widehat{Da_s} )^*\widehat{Da_s} e^{-\beta\OPERW{H}} \big)
\TR\big((\widehat{Db_\tau})^*\widehat{Db_\tau}e^{-\beta\OPERW{H}}\big)\Big)^{1/2}\Rd \tau\Rd s\,.\\
\end{split}
\end{equation}
The following two lemmas, proved in Section~\ref{sec_lemmas}, estimate the remainder terms $T_0$ and $T_3$
\begin{lemma}[Mean-field approximation]\label{lemma3}
Assume that the bounds in Theorem \ref{thm_mf} hold, then
\begin{equation}\label{3.1}
\frac{|T_0|}{\TR (e^{-\beta \OPERW{H}})}=\BIGO \big(t\epsilon_1^2+tM^{-1}\big) \,.\\
\end{equation}
\end{lemma}

\begin{lemma}[Composition analysis]\label{lemma1} Assume that the bounds in Theorem \ref{thm_mf} hold and that $c$ and $d$ are in the Schwartz space $\SSPACE$, then 
\[
\begin{split}
c\MP d  &= cd + \mathfrak{r}_{cd}\,,\\
\lim_{M\to\infty} M^{1/2}\mathfrak{r}_{cd} &=\frac{{\mathrm{i}}}{2}(\nabla_{z}\cdot\nabla_{z'}')c(z)d(z')\big|_{z=z'}\,,
\end{split}
\]
and if $c$ and $d$ are scalar valued
\[
\begin{split}
\frac{1}{2}(c\MP d + d\MP c) &=
    cd + r_{cd}\,,\\
\lim_{M\to\infty} Mr_{cd} &=\frac{1}{8}(\nabla_{z}\cdot\nabla_{z'}')^2c(z)d(z')\big|_{z=z'}\,,\\
\end{split}
\]
and if $c$  is scalar valued
\[
\begin{split}
&{\IU}M^{1/2}(H\MP c-c\MP H)= \{H,c\}
+{r}_c\,,\\
&{r}_c =\frac{1}{8M}\int_0^1\cos(\frac{\bar s}{2M^{1/2}}\nabla_{z}\cdot\nabla_{z'}')
(\nabla_{z}\cdot\nabla_{z'}')^3H(z)c(z')\big|_{z=z'}(1-\bar s)^2\Rd \bar s\,,\\
&\lim_{M\to\infty} M{r}_c=\frac{1}{24}(\nabla_{z}\cdot\nabla_{z'}')^3H(z)c(z')\big|_{z=z'}\,,\\
\end{split}
\]
where the limits hold in $L^1(\rset^{2N})$ and $L^\infty(\rset^{2N})$.
Furthermore the function $a_t:\rset^{2N}\to \mathbb C$, defined in \eqref{eq:at}, is
in the Schwartz class and there holds
\[
\frac{|T_3|}{\TR (e^{-\beta \OPERW{H}})}=\BIGO \big(t^2(\epsilon_2^2+M^{-1})\big) \,.
\]
\end{lemma}

The combination of \eqref{first_step}-\eqref{forth_step}, Lemmas \ref{lemma3} and \ref{lemma1} imply
\begin{equation}\label{dyn_pre}
\frac{|\TR\big((\OPER{A}_t\OPER{B}_0e^{-\beta\OPERW{H}} -\OPER{a}_t\OPER{B}_0e^{-\beta\OPERW{H}})|}{\TR (e^{-\beta \OPERW{H}})}=\mathcal O(tM^{-1}+t\epsilon_1^2+t^2\epsilon_2^2)
\,.
\end{equation}
Similarly the symmetrized
difference has the same bound
\begin{equation}\label{dynamic_step}
\begin{split}
&\TR\big((\OPER{A}_t\OPER{B}_0+\OPER{B}_0\OPER{A}_t)e^{-\beta\OPER{H}}\big)-\TR\big((\OPER{a}_t\OPER{B}_0+\OPER{B}_0\OPER{a}_t)e^{-\beta\OPER{H}}\big)\\
&=
\TR\big((\OPER{A}_t\OPER{b}_0+\OPER{b}_0\OPER{A}_t)e^{-\beta\OPER{H}}\big)-\TR\big((\OPER{a}_t\OPER{b}_0+\OPER{b}_0\OPER{a}_t)e^{-\beta\OPER{H}}\big)
=\mathcal O(tM^{-1}+t\epsilon_1^2+t^2\epsilon_2^2)\TR (e^{-\beta \OPERW{H}})\,,
\end{split}
\end{equation}
obtained by interchanging the role of $\OPER{A}$ and $\OPER{B}$ in \eqref{dyn_pre}.

{\bf Step 2.} Here we estimate the second term in the left hand side of \eqref{dynamic_step}, which can be split into
\begin{equation}\label{stat_1}
\TR\big((\OPER{a}_t\OPER{b}_0+\OPER{b}_0\OPER{a}_t)e^{-\beta\OPER{H}}\big)=\TR\big((\OPER{a}_t\OPER{b}_0+\OPER{b}_0\OPER{a}_t)
\OPER{e^{-\beta H}} \big)+
\TR\big((\OPER{a}_t\OPER{b}_0+\OPER{b}_0 \OPER{a}_t)(e^{-\beta\OPER{H}}-\OPER{e^{-\beta H}})\big)\,.
\end{equation}
The first term in the right hand side in \eqref{stat_1} has by
Lemma \ref{lemma_composition} the classical molecular dynamics approximation
\begin{equation}\label{md-y0}
\begin{split}
    &\TR\big(\frac{\OPER{a}_t\OPER{b}_0+\OPER{b}_0\OPER{a}_t}{2}\widehat{{e^{-\beta H}}}\big)= \NORMFAC^N \int_{\mathbb R^{2N}}  \TR\big(\frac{(a_t\MP b_0+b_0\MP a_t)(x,p)}{2} {e^{-\beta H(x,p)}}\big) \Rd x\Rd p\\
 &= \NORMFAC^N \int_{\mathbb R^{2N}}  \TR\Big(a_0\big(z_t(z_0)\big)b_0(z_0) {e^{-\beta H(z_0)}}\Big)+\TR\big(r_{ab}(z_0){e^{-\beta H(z_0)}}\big) \Rd z_0\,,\\
\end{split}
\end{equation}
where by Lemma \ref{lemma1}
\begin{equation}\label{Mr_ab}
\begin{split}
&\lim_{M\to\infty} M\int_{\mathbb R^{2N}}\TR\big(r_{ab}(z_0){e^{-\beta H(z_0)}}\big)\Rd z_0\\
&=\frac{1}{16}\int_{\mathbb R^{2N}}
\TR\big(e^{-\beta H(z_0)}\big)(\nabla_{z_0}\cdot\nabla_{z_0'}')^2
\Big(a_0\big(z_t(z_0)\big)b_0(z_0')\Big)\big|_{z_0=z_0'}\Rd z_0\,.
\end{split}
\end{equation}
It remains to estimate the second term in the right hand side of
\eqref{stat_1}. We have
\[
\begin{split}
&\TR\big(\frac{\OPER{a}_t\OPER{b}_0+\OPER{b}_0\OPER{a}_t}{2}(\widehat{{e^{-\beta H}}}-\widehat\rho)\big)
=\NORMFAC^N\int_{\rset^{2N}}\TR\big( \frac{a_t\MP b_0+b_0\MP a_t}{2}({e^{-\beta H}}-\rho)\big)\Rd z
\end{split}
\]
and 
Lemmas~\ref{path_lemma} and \ref{lemma1} imply
\begin{equation}\label{Mrho}
\begin{split}
&\lim_{M\to\infty} \Big(M\int_{\rset^{2N}}\TR\big( \frac{a_t\MP b_0+b_0\MP a_t}{2}({e^{-\beta H}}-\rho)\big)\Rd z\Big)\\
&= \int_{\rset^{2N}}\TR\Big( 2^{-1}({a_tb_0+b_0a_t})(x,p)
e^{-\beta |p|^2/2}\times\\
&\quad\times\int_0^\beta(\frac{\beta}{2}-r)^2 e^{-rV(x)}p^TV''(x)p e^{-(\beta-r)V(x)}\Rd r\Big)
\Rd x\Rd p\\
&\quad + \int_{\rset^{2N}}\TR\Big( 2^{-1}({a_tb_0+b_0a_t})(x,p)
e^{-\beta |p|^2/2}\times\\
&\qquad\times
\int_0^\beta
 \int_0^r (\frac{\beta}{2}-s)(\frac{\beta}{2}-r) e^{-sV(x)}V'(x)pe^{-(r-s)V(x)} V'(x)p
e^{-(\beta-r)V(x)}{\mathrm{d}}s{\mathrm{d}}r\Big)\Rd x\Rd p\\
&= \frac{\beta^3}{24}\int_{\rset^{2N}} ({a_tb_0+b_0a_t})(x,p)
e^{-\beta |p|^2/2} \TR\big( p^TV''(x)p e^{-\beta V(x)}\big)\Rd x\Rd p+\BIGO(1)\\
&=\BIGO (1)\,,
\end{split}
\end{equation}
where the second last equality follows by interchanging the order of
the trace and the integration with respect to $\beta$ and using the cyclic invariance of the trace. The first equality is obtained by splitting the integral as 
\[
\int_{\rset^{2N}}\ldots =\int_{\mathbb L_m^c}\ldots\ + \
\int_{\mathbb L_m}\ldots \] 
over a compact set $\mathbb L_m^c:=\{z\in\rset^{2N}\ :\ h(z)\le m\}$ and its complement $\mathbb L_m:=\{z\in\rset^{2N}\ :\ h(z)>m\}$ and 
using that $\rho$ is uniformly bounded in $L^\infty(\rset^{2N})$:  the second integral  is zero for $m$ sufficiently large,
as verified in \eqref{alfa_int} below, and in the compact set we apply dominated convergence. 

{\it Verification of $\int_{\mathbb L_m}|a_t(z_0)|{\mathrm{d}}z_0=0$.} The integration with respect to the initial data measure $\Rd  z_0$ can be replaced by integration with respect
to $\Rd  z_t$ since the phase space volume is preserved, i.e. the Jacobian determinant 
\[
|{\rm det}\big(\frac{\partial z_0}{\partial z_t}\big)| =1
\]
 is constant for all time by Liouville's theorem. 
 We first verify that
 \begin{equation}\label{h_infty}
 h(z)\to\infty \ \mbox{as $|z|\to\infty$}\,.
 \end{equation}
 By assumption $\int_{\rset^{N}}\TR e^{-\beta V(x)}\Rd x<\infty$. Therefore the smallest eigenvalue $\lambda_0(x)$
 of $V(x)$ also satisfies $\int_{\rset^{N}} e^{-\beta \lambda_0(x)}\Rd x<\infty$
 and as $h(x,p)\ge |p|^2/2+\lambda_0(x)$ we have $\int_{\rset^{2N}}e^{-\beta h(x,p)}\Rd x\Rd p<\infty$,
 which combined with the assumption $\|\nabla_xh(x,p)\|_{L^\infty(\rset^{2N})}\le C$ establishes \eqref{h_infty}.
 By using the two properties $h\big(z_t(z_0)\big)=h(z_0)$ and $h(z)\to\infty$ as $|z|\to\infty$
 together with the compact support of $a_0$ we obtain
\begin{equation}\label{alfa_int}
\begin{split}
\int_{\mathbb L_m} |
a_{0}\big(z_t(x_0,p_0)\big)| \Rd  z_0
&=\int_{\mathbb L_m} |
a_0\big(z_t(x_0,p_0)\big)|
|{\rm det}\big(\frac{\partial z_0}{
\partial z_t}\big)|\Rd  z_t\\
&=\int_{\mathbb L_m} |
a_{0}(z)|
\Rd  z = 0\ \mbox{as $m\to\infty$}\PERIOD\\
\end{split}
\end{equation}

In conclusion we have 
\[
\frac{1}{2}\TR\big((\OPER{a}_t\OPER{b}_0+\OPER{b}_0\OPER{a}_t)e^{-\beta\OPER{H}}\big)=\NORMFAC^N \Big(\int_{\mathbb R^{2N}}  \TR\big(a_0(z_t(z_0))b_0(z_0) {e^{-\beta H}}(z_0)\big)
+\mathcal O(M^{-1})\Big)\,,
\]
which combined with \eqref{dynamic_step} proves the theorem for $\Psi={\rm I}$.

{\bf Step 3.} 
To improve the error estimate $\epsilon_2^2$ we study the transformed
Hamiltonian operator \[
\widehat{\bar H}:=\Psi^* \widehat H \Psi\]
where
$\Psi:\rset^{N}\rightarrow \CSP$ and $\Psi(x)$ is any twice differentiable orthogonal matrix  with the Hermitian transpose
 $\Psi^*(x)$.
This Step 3 has the three substeps:
\begin{itemize}
\item[Step 3.1] study the dynamics under $\bar H$, 
\item[Step 3.2] analyse $\bar H=\Psi^*\# H\#\Psi$,
\item[Step 3.3] modify Steps 1 and 2.
\end{itemize}

{\it Step 3.1.} Let $\alpha$ be any complex number and define for $t\in\mathbb R$ the exponential 
\[
\OPERW{\by_t}:=\Psi^* \EXP{t\alpha\hat H}\Psi.
\]
Differentiation shows that
\begin{equation*}\label{eq:y_t_diffeq}
\partial_t \OPERW{\by_t} = \alpha  \Psi^*\hat H\Psi \Psi^* \EXP{t\alpha\hat H}\Psi
= \alpha  \Psi^*\hat H\Psi\OPERW{\by_t}
\end{equation*}
and consequently 
\begin{equation*}\label{y_t_exp}
\OPERW{\by_t} =\EXP{t\alpha \Psi^*\hat H\Psi}=e^{t\alpha \widehat{\bar H}}\PERIOD
\end{equation*}
Therefore the transformed variable
\begin{equation}\label{tilde_a_t}
\OPERW{\bA(t,z)} :=  \Psi^* \OPER{ A}_t\Psi \COMMA\quad t \in\rset\,,
\end{equation}
evolves with the dynamics $\widehat{\bar H}$
\[
\begin{split}
\widehat{\bar A_t} &=  \Psi^* \OPER{ A}_t\Psi\\
&= \Psi^* e^{{\mathrm{i}} tM^{1/2}\widehat{ H}}\OPER{A_0}e^{-{\mathrm{i}} tM^{1/2}\widehat{H}}\hat \Psi\\
&=  e^{{\mathrm{i}} tM^{1/2}\widehat{ \bar H}}\OPER{\bar A_0}e^{-{\mathrm{i}} tM^{1/2}\widehat{\bar H}}\,.\\
\end{split}
\]
and the cyclic property of the trace implies that the quantum observable satisfies
\begin{equation}\label{tr_abe}
\TR(\OPER{A_t}\OPER{B_0}e^{-\beta \OPER{H}})=
\TR(\underbrace{\Psi\Psi^*}_{={\rm I}}\OPER{A_t}\OPER{B_0}e^{-\beta \OPER{H}})
=\TR(\Psi^*\OPER{A_t}\Psi\Psi^*\OPER{B_0}\Psi\Psi^*e^{-\beta \OPER{H}}\Psi)
=\TR(\OPER{\bar A_t}\OPER{\bar B_0}e^{-\beta \OPER{\bar H}})\,,
\end{equation}
where the initial symbols are given by
\[
\begin{split}
\bar A_0  &= \Psi^* \# a_0\# \Psi\,,\\
\bar B_0  &= \Psi^* \# b_0\# \Psi\,.\\
\end{split}
\]

{\it Step 3.2.} We have 
\begin{equation}\label{bar_H}
\bar H(x,p)= \Psi^*\MP  H\MP \Psi(x,p)
=\Psi^*(x)H(x,p)\Psi(x) + \frac{1}{4M}\nabla\Psi^*(x)\cdot\nabla\Psi(x)\,,
\end{equation}
as derived in \cite{KPSS} from 
the composition in \eqref{composition_bis}:
 \[
 \begin{split}
& \Psi^*\MP  H\MP \Psi(x,p)\\
 &=\Psi^*(x)\MP  \big(H(x,p)\Psi(x) +\frac{\IU M^{-1/2}}{2} p\cdot \nabla\Psi(x) -\frac{M^{-1}}{4}\Delta\Psi(x)\big)\\
 &=\Psi^*(x)H(x,p)\Psi(x)+\frac{\IU M^{-1/2}}{2}p\cdot \nabla\Psi^*(x)\Psi(x)-\frac{M^{-1}}{4}\Delta\Psi^*(x)\Psi(x)\\
&\quad +\frac{\IU}{2M^{1/2}} \Psi^*(x) p\cdot \nabla\Psi(x)
 -\frac{1}{4M}\nabla\Psi^*(x)\cdot\nabla\Psi(x)-\frac{1}{4M}\Psi^*(x)\Delta\Psi(x)\\
  \end{split}
 \]
where by the orthogonality $\Psi^*\Psi=\Id$ 
 \[
 \begin{split}
 \Psi^*\MP  H\MP \Psi(x,p)&=\Psi^*(x)H(x,p)\Psi(x)+\frac{\IU}{2M^{1/2}}p\cdot\nabla\big(\Psi^*(x)\Psi(x)\big)\\
&\qquad  -\frac{1}{4M}\Big(\Delta \big(\Psi^*(x)\Psi(x)\big) -\nabla\Psi^*(x)\cdot\nabla\Psi(x)\Big)\\
  &=\Psi^*(x)H(x,p)\Psi(x) + \frac{1}{4M}\nabla\Psi^*(x)\cdot\nabla\Psi(x)\PERIOD
 \end{split}
 \]
 Let $\Psi(x)$ be the orthogonal matrix composed of the eigenvectors to $V(x)$, then the matrix 
 \[
 \Psi^*(x)H(x,p)\Psi(x)= \frac{|p|^2}{2}{\rm I} + \Lambda(x)
 \]
 is diagonal, with the eigenvectors $\lambda_i(x)$ of $V(x)$ forming the diagonal $d\times d$ matrix $\Lambda(x)$. 
 The non diagonal part $ \frac{1}{4M}\nabla\Psi^*(x)\cdot\nabla\Psi(x)$ of $\bar H(x,p)$ is small  if $\Psi(x)$ is differentiable everywhere. If the eigenvectors to $V$ are not differentiable in a point $x_*$, we may use a regularized  version of $\Psi$ in a neighbourhood of $x_*$ to form an approximate diagonalization of $H$.
 

{\it Step 3.3.} The derivation in Step 1 can now be repeated with $H$ replaced by $\bar H$ and $A,B$ by $\bar A,\bar B$.
Duhamel's principle \eqref{duhamel_1} implies
\begin{equation*}\label{duhamel_22}
\OPER{\bar A_t} -\OPER{a_t}  -(\OPER{\bar A_0}- \OPER{a_0})= 
\int_0^t  e^{{\IU}M^{1/2}(t-s) \OPERW{\bar H}}\big({\IU}M^{1/2}
[\OPERW{\bar H},\OPER{ a}_s]- \{ h, a_s\}^{\widehat{}}\ \big) e^{{-\rm i}M^{1/2}(t-s) \OPERW{\bar H}}
\Rd s\,,\\
\end{equation*}
and we obtain by \eqref{tr_abe} as in \eqref{dyn_pre}
\[
\begin{split}
\TR\big((\widehat{ A}_t\OPER{ B}_0e^{-\beta\OPERW{ H}} 
-\OPER{a}_t\OPER{\bar B}_0e^{-\beta\OPERW{\bar H}}) &=
\TR\big((\widehat{\bar A}_t\OPER{\bar B}_0e^{-\beta\OPERW{\bar H}} -\OPER{a}_t\OPER{\bar B}_0e^{-\beta\OPERW{\bar H}})  \\
&=T_0+T_3 +\TR\big((\widehat{\bar A_0}-\OPER{a_0})\OPER{\bar B}_0e^{-\beta\OPERW{\bar H}}\big)
+ \int_0^t \TR\big(\widehat{Da_s}(\widehat{\bar B_0}-\widehat{b_0})e^{-\beta\widehat{\bar H}} \big){\mathrm{d}}s\,,\\
\end{split}
\]
where $H$ is replaced by $\bar H$ in $\widehat{Da}$, $\widehat{Db}$, $T_0$ and $T_3$, so that
\[
\begin{split}
T_0 &=\int_0^t \TR(\underbrace{\widehat{Da_s}}_{=:\{\bar H-h,a_s)\}^{\widehat{} \, +\hat r^a}} \OPER{b}_{s-t} e^{-\beta\OPERW{\bar H}})\Rd s\,,\\
T_3 &= \int_0^t\int_0^{s-t} \TR(e^{-\beta\OPERW{\bar H}}\widehat{Da_s} e^{{\IU}M^{1/2}(s-t-\tau) \OPERW{\bar H}}\widehat{D b_\tau}
e^{-{\IU}M^{1/2}(s-t-\tau) \OPERW{\bar H}}\Rd \tau\Rd s\,.
\end{split}
\]

The two terms $T_0$ and $T_3$ have the bounds in Lemmas \ref{lemma3} and \ref{lemma1} with $H$ replaced by $\bar H$.
It remains to show that the initial errors  satisfy
\begin{equation}\label{initial_bar}
\frac{\TR\big((\widehat{\bar A_0}-\OPER{a_0})\OPER{\bar B}_0e^{-\beta\OPERW{\bar H}}\big)}{
\TR(e^{-\beta\OPERW{\bar H}})}
+ \frac{\int_0^t \TR\big(\widehat{Da_s}(\widehat{\bar B_0}-\widehat{b_0})
e^{-\beta\widehat{\bar H}} \big){\mathrm{d}}s}{\TR(e^{-\beta\OPERW{\bar H}})}=\mathcal O(M^{-1})\,.
\end{equation}
 To prove \eqref{initial_bar} let $\widehat{\bar \rho}= e^{-\beta \widehat{\bar H}}$, then by the composition \eqref{composition_bis} and Lemma  \ref{lemma_composition} we obtain
\[
\begin{split}
\TR\big((\widehat{\bar A_0}-\OPER{a_0})\OPER{\bar B}_0e^{-\beta\OPERW{\bar H}}\big)
&= \TR\Big(\big((\bar A_0 -a_0)\# \bar B_0\big)^{\widehat{}} \ \widehat{\bar\rho}\Big)\\
&= \big(\frac{\sqrt{M}}{2\pi}\big)^N\int_{\rset^{2N}} \TR\Big(\big((\bar A_0 -a_0)\# \bar B_0\big)
\,  \bar \rho\Big) \, {\mathrm{d}}x{\mathrm{d}}p\,.
\end{split}
\]
To estimate the initial error $\bar A_0-a_0$ we use  Lemma \ref{lemma1}
\[
\begin{split}
a_0\# \Psi(x,p) &= a_0(x,p)\Psi(x) + \frac{{\mathrm{i}}}{2\sqrt M} \nabla_p a_0(x,p)\cdot \nabla_x \Psi(x)\\
&\quad - \frac{1}{4M} \int_0^1 e^{-\frac{{\mathrm{i}} s}{2} M^{-1/2} \nabla_x\cdot \nabla_p} (\nabla_x\cdot\nabla_p)^2 a_0(x',p)\Psi(x)(1-s){\mathrm{d}}s \big|_{x'=x}\,,
\end{split}
\]
which by the composition \eqref{composition_bis} implies
\[
\begin{split}
\bar A_0(x,p) &= \Psi^*\# a_0\#\Psi(x,p)\\
&= \Psi^*(x) (a_0\# \Psi)(x,p) + \frac{{\mathrm{i}}}{2\sqrt M} \nabla_x\Psi^*(x)\cdot \nabla_p (a_0\#\Psi)(x,p)\\
&\quad -\frac{1}{4M} \int_0^1 e^{-\frac{{\mathrm{i}} s}{2} M^{-1/2} \nabla_x\cdot \nabla_p} (\nabla_x\cdot\nabla_p)^2 \Psi^*(x) (a_0\#\Psi)(x',p)(1-s){\mathrm{d}}s \big|_{x'=x}\\
&= a_0(x,p)\underbrace{\Psi^*(x)\Psi(x)}_{={\rm I}} + \frac{{\mathrm{i}}}{2\sqrt M}  \nabla_p a_0(x,p) \cdot 
\underbrace{\nabla_x\big(\Psi^*(x)\Psi(x)\big)}_{=0}\\
&\quad - \frac{1}{4M} \Psi^*(x) \int_0^1 e^{-\frac{{\mathrm{i}} s}{2} M^{-1/2} \nabla_x\cdot \nabla_p} (\nabla_x\cdot\nabla_p)^2 \big(a_0(x',p)\Psi(x)\big)(1-s){\mathrm{d}}s \big|_{x'=x}\\
&\quad -\frac{1}{4M} (\nabla_{x'}\cdot\nabla_p)(\nabla_x\cdot\nabla_p) \big(\Psi^*(x')a_0(x'',p)\Psi(x)\big)\big|_{x''=x'=x}\\
&\quad +\frac{{\mathrm{i}}}{8M^{3/2}} \int_0^1 e^{-\frac{{\mathrm{i}} s}{2} M^{-1/2} \nabla_x\cdot \nabla_p} (\nabla_x\cdot\nabla_p)^2 \Psi^*(x) \big(a_0(x',p)\Psi(x)\big)(1-s){\mathrm{d}}s \big|_{x'=x}\\
&- \frac{1}{4M} \int_0^1 e^{-\frac{{\mathrm{i}} s}{2} M^{-1/2} \nabla_x\cdot \nabla_p} (\nabla_x\cdot\nabla_p)^2 
\Psi^*(x) (a_0\#\Psi)(x',p)(1-s){\mathrm{d}}s \big|_{x'=x}\,.\\
\end{split}
\]
Here the orthogonality $\Psi^*\Psi={\rm I}$ implies $\nabla(\Psi^*\Psi)=0$. We obtain as in \eqref{Mr_ab} the limit
\[
\begin{split}
&\lim_{M\to\infty}  M\int_{\rset^{2N}} \TR\big((\bar A_0-a_0)\bar B_0 \bar \rho\big)\, {\mathrm{d}}x{\mathrm{d}}p\\
&= -\int_{\rset^{2N}} \TR\Big(\big[\frac{1}{8} \Psi^*(x) (\nabla_x\cdot\nabla_p)^2 a_0(x',p)\Psi(x)\big|_{x'=x}\\
&\quad +\frac{1}{4}(\nabla_{x'}\cdot\nabla_p)(\nabla_x\cdot \nabla_p)\Psi^*(x')a_0(x'',p)\Psi(x)\big|_{x''=x'=x}\\
&\quad +\frac{1}{8}(\nabla_x\cdot\nabla_p)^2 \Psi^*(x)a_0(x',p)\Psi(x')\big|_{x'=x}\big] b_0(x,p) e^{-\beta \Psi^*(x)H(x,p)\Psi(x)}\Big){\mathrm{d}}x{\mathrm{d}}p
\end{split}
\]
and similarly
\[
\begin{split}
\int_0^t \TR\big(\widehat{Da_s}(\widehat{\bar B_0}-\widehat{b_0})e^{-\beta \widehat{\bar H}}\big){\mathrm{d}} s
= \big(\frac{\sqrt M}{2\pi}\big)^N\int_{\rset^{2N}} \TR\Big(\big(Da_s\#(\bar B_0-b_0)\big) \bar \rho\Big) {\mathrm{d}} x {\mathrm{d}}p
\end{split}
\]
where
\[
\begin{split}
&\lim_{M\to\infty} \big(M\int_{\rset^{2N}} \TR\Big(\big(Da_s\#(\bar B_0-b_0)\big) \bar \rho\Big) {\mathrm{d}} x {\mathrm{d}}p\big)\\
&= -\int_{\rset^{2N}} \TR\Big(\{\bar H(x,p)-h(x,p),a_s(x,p)\}
\big[\frac{1}{8} \Psi^*(x) (\nabla_x\cdot\nabla_p)^2 b_0(x',p)\Psi(x)\big|_{x'=x}\\
&\quad +\frac{1}{4}(\nabla_{x'}\cdot\nabla_p)(\nabla_x\cdot \nabla_p)\Psi^*(x')b_0(x'',p)\Psi(x)\big|_{x''=x'=x}\\
&\quad +\frac{1}{8}(\nabla_x\cdot\nabla_p)^2 \Psi^*(x)b_0(x',p)\Psi(x')\big|_{x'=x}\big]
e^{-\beta \Psi^*(x)H(x,p)\Psi(x)}\Big){\mathrm{d}}x{\mathrm{d}}p\,.
\end{split}
\]
which proves \eqref{initial_bar}.

\end{proof}

\section{Proof of Lemmas}\label{sec_lemmas}
This section estimates the remainder terms $T_0$ and $T_3$ and the statistical error $\rho-e^{-\beta H}$.
The error term $T_3$ in Lemma \ref{lemma1} is due to remainders from classical approximation while $T_0$ in Lemma \ref{lemma3} is the main term regarding the mean-field approximation. The statistical error is estimated in Lemma \ref{path_lemma}.

\subsection{Proof Lemma \ref{lemma3}}
\begin{proof} We consider first the case $\Psi={\rm I}$, i.e. $\bar H=H$.
We have 
\begin{equation}\label{L3.1}
\begin{split}
\int_0^t \TR(\widehat{Da_s}\OPER{b}_{s-t}e^{-\beta\OPERW{H}})\Rd s
& = \int_0^t \TR\big((\{H,a_s\}^{\widehat{}} + \hat r_s^a -\{h,a_s\}^{\widehat{}}\ )\OPER{b}_{s-t}e^{-\beta\OPERW{H}} \big)\Rd s\\
&=  \int_0^t \TR\big(\{H-h,a_s\}^{\widehat{}} \ \OPER{b}_{s-t}e^{-\beta\OPERW{H}}\big)\Rd s
+ \int_0^t \TR\big(\hat r_s^a \OPER{b}_{s-t}e^{-\beta\OPERW{H}} \big)\Rd s\\
&=(\frac{M^{1/2}}{2\pi})^{N}\int_0^t \int_{\rset^{2N}}\TR\big((\{H-h,a_s\}\MP b_{s-t}) \rho\big)\Rd z\Rd s\\
&\quad +(\frac{M^{1/2}}{2\pi})^{N}\int_0^t \int_{\rset^{2N}}\TR\big((r_s^a\MP b_{s-t}) \rho\big)\Rd z\Rd s
\end{split}
\end{equation}
where the trace in the integrals over phase space $\rset^{2N}$ is with respect to $d\times d$ matrices. 
The inner integral in the second term of the right hand side in \eqref{L3.1} has by Lemmas \ref{path_lemma} and \ref{lemma1} the limit
\begin{equation}\label{3.2}
\begin{split}
&\lim_{M\to\infty}\Big(M\int_{\rset^{2N}}\TR\big((r_s^a \MP b_{s-t}) \rho)\big)\Rd z \Big)\\
&=
\frac{1}{24}\int_{\rset^{2N}}
\TR\Big((\nabla_{z_0'}\cdot\nabla_{z_0}')^3\big(H(z_0)a_0\big(z_s(z_0')\big)\big)\big|_{z_0=z_0'}
b_0\big(z_{s-t}(z_0)\big) e^{-\beta H(z_0)}\Big)\Rd z_0\,,
\end{split}
\end{equation}
using splitting of the phase space integral as in \eqref{Mrho}.
Similarly we have by Lemmas \ref{path_lemma} and \ref{lemma1} the limit
\[
\lim_{M\to\infty}\int_{\rset^{2N}} \TR\big((\{H-h,a_s\}\MP b_{s-t} )\rho\big)\Rd z
=\int_{\rset^{2N}} \TR(\{H-h,a_s\}b_{s-t}{e^{-\beta H}})\Rd z\,
\]
and integration by parts together with the mean-field definition \eqref{mf-def} simplify
 the first term in the right hand side of \eqref{L3.1} to
 \begin{equation}\label{3.3}
 \begin{split}
& \int_{\rset^{2N}}\TR(\{H-h,a_s\}b_{s-t}{e^{-\beta H}})\Rd z\\
& =\int_{\rset^{2N}}\TR\big(\NABLAR(H-h)\cdot\nabla a_s \ b_{s-t}e^{-\beta H}\big)\Rd z\\
&= -\int_{\rset^{2N}}\underbrace{\TR\big((H-h)  e^{-\beta H} \big)}_{=0}\NABLAR\cdot (b_{s-t}\nabla a_s)\Rd z\\
&\quad +\int_{\rset^{2N}}\TR\big((H-h) 
\int_0^\beta e^{-\tau H} \nabla a_s\cdot\NABLAR (H-h) e^{-(\beta-\tau) H} b_{s-t}{\mathrm{d}}\tau\big) \Rd z\\
 &\quad +\int_{\rset^{2N}}\underbrace{\TR\big((H-h)e^{-\beta H}\big)}_{=0} \nabla a_s\cdot \NABLAR h \ \beta b_{s-t}\Rd z\\
 &=\frac{1}{2}\int_{\rset^{2N}}\TR\Big( 
\int_0^\beta e^{-\tau H}\big( (H-h)\nabla a_s\cdot\NABLAR (H-h) \\
&\qquad +(\nabla a_s\cdot\NABLAR (H-h) (H-h))\big) e^{-(\beta-\tau) H}b_{s-t}{\mathrm{d}}\tau\Big) \Rd z\\
 &=\frac{\beta}{2}\int_{\rset^{2N}}\TR\big( \NABLAR(H-h)^2\cdot\nabla a_se^{-\beta H}b_{s-t}\big)\Rd z\\
 &=-\frac{\beta}{2}\int_{\rset^{2N}}\TR\big( (H-h)^2 \nabla a_s \cdot\NABLAR (e^{-\beta H}b_{s-t})\big)\Rd z\\
 &\quad -\frac{\beta}{2}\int_{\rset^{2N}}\TR\big( (H-h)^2 \underbrace{\NABLAR\cdot\nabla (a_s)}_{=0}  (e^{-\beta H}b_{s-t})\big)\Rd z\\
 &=-\frac{\beta}{2}\int_{\rset^{2N}}\TR\Big( (H-h)^2
 \nabla a_s \cdot (\NABLAR b_{s-t}e^{-\beta H}- b_{s-t}\int_0^\beta e^{-\tau H} \NABLAR H e^{-(\beta-\tau) H}{\mathrm{d}}\tau)\Big)\Rd z\\
 &=-\frac{\beta}{2}\int_{\rset^{2N}}\TR\big(e^{-\beta H} (H-h)^2 \nabla a_s \cdot 
 (\NABLAR b_{s-t}-\beta b_{s-t}\NABLAR H)\big)\Rd z\,.\\
 \end{split}
 \end{equation}
 The third equality uses that $(H-h)$ commutes with $e^{-\beta H}$ and the forth and the last equality is obtained by interchanging the order of the trace and the integral with respect to $\tau$ in combination with the cyclic invariance of the trace. Cauchy's inequality, the positive definiteness of $e^{-\beta H} (H-h)^2$
 and $H-h$ depending only on $x$ 
 imply that the right hand side has the bound
 \begin{equation}\label{Kx}
 \begin{split}
 &|\frac{\beta}{2}\int_{\rset^{2N}}\TR\big(e^{-\beta H} (H-h)^2 \nabla a_s \cdot 
 (\NABLAR b_{s-t}-\beta b_{s-t}\NABLAR H)\big)\Rd z|\\
 &\le \frac{\beta}{2}\int_{\rset^{2N}}\Big(\TR\big((e^{-\beta H} (H-h)^2)^2\big)
 \TR\big((\nabla a_s \cdot 
 (\NABLAR b_{s-t}-\beta b_{s-t}\NABLAR H))^2\big)\Big)^{1/2}\Rd z\\
 &\le\frac{\beta}{2}\int_{\rset^{2N}}\TR\big(e^{-\beta V} (H-h)^2\big)
 e^{-\beta|p|^2/2}\Big(\TR\big((\nabla a_s \cdot 
 (\NABLAR b_{s-t}-\beta b_{s-t}\NABLAR H))^2\big)\Big)^{1/2}\Rd z\\
 &\le K \int_{\rset^N} \TR\big(e^{-\beta V} (H-h)^2\big){\mathrm{d}}x
 \end{split}
 \end{equation}
 for a positive constant $K$.
The combination of \eqref{3.2}-\eqref{Kx} implies \eqref{3.1} and we note that the exponential form 
 $e^{-\beta H}$ of the Gibbs density was crucial to obtain \eqref{3.3}.
 
 In the case that $\Psi\ne {\rm I}$ we have by \eqref{bar_H} $\bar H=\Psi^*H\Psi + \frac{1}{4M} \nabla \Psi^*\cdot \nabla\Psi$
 and  the factor  $\TR\big((H-h)e^{-\beta H}\big)=0$ in \eqref{3.3}  is replaced by 
 $\TR\big((\bar H-h)e^{-\beta \bar H}\big)$ where
 \[
\lim_{M\to\infty} \Big(M \TR\big((\bar H-h)e^{-\beta \bar H}\big)\Big)
=
 \TR\big(\nabla\Psi^*\cdot\nabla\Psi e^{-\beta \Psi^* H\Psi}\big)\,.
 \]
 which as in \eqref{3.3} and \eqref{Kx} imply
 \[
 |T_0|\le \mathcal O(tM^{-1} + t\epsilon_1^2)\,.
 \]
\end{proof}

\subsection{Proof of Lemma \ref{path_lemma}}\label{sec_path-proof}
\begin{proof}
To estimate the difference of the symbols for the Gibbs density, $\rho-{e^{-\beta H}}$, we define
the solution operator 
\[
\dot{\YT_t}(x)=-\YT_t(x)V(x)\,,\quad\quad \YT_0={\Id}
\]
satisfying $\YT_{t}= e^{-tV}$, which implies ${e^{-\beta H(x,p)}}=e^{-\beta( |p|^2/2+V(x))}=
e^{-\beta |p|^2/2} \,\YT_{\beta}(x)$. We have by Duhamel's principle
\begin{equation}\label{yY}
\BARYT^\pm_\beta-\YT_\beta= \int_0^\beta \BARYT^\pm_t\,(V- V^\pm_t) e^{-(\beta-t)V}\Rd  t\,,
\end{equation}
which by dominated convergence implies
\begin{equation}\label{y+-}
\begin{split}
\lim_{M\to\infty} M^{1/2}(\BARYT^\pm_t-\YT_t) &=
\lim_{M\to\infty} M^{1/2}\int_0^t\BARYT_s^\pm(V-V_s^\pm)e^{-(t-s)V}{\mathrm{d}}s\\
&=\mp\int_0^te^{-sV(x)}V'(x)W^\beta_s e^{-(t-s)V(x)}{\mathrm{d}}s\,.\\
\end{split}
\end{equation}
The symmetrized relation \eqref{matFeynman-Kac}  yields
\[
\rho-{e^{-\beta H}}= \frac{1}{2}\E[e^{-{\IU}W_\beta\cdot p}\big((\BARYT_\beta^+-\YT_\beta)+ (\BARYT_\beta^--\YT_\beta)\big)]
\]
which by \eqref{yY}, \eqref{y+-} and dominated convergence establish
\begin{equation}\label{4.8}
\begin{split}
&\lim_{M\to\infty} \Big(M\big(\rho(x,p)-{e^{-\beta H(x,p)}}\big)\Big) \\
&=\frac{1}{2}\E[e^{-{\IU}W_\beta\cdot p}
\int_0^\beta \lim_{M\to\infty}\big(\BARYT^+_t M(V-V_t^+) +\BARYT^-_t M(V-V_t^-)\big)e^{-(\beta-t)V}\Rd  t]\\
&=\E\Big[e^{-{\IU}W_\beta\cdot p}\Big(\frac{1}{2}\int_0^\beta \lim_{M\to\infty}\big((\BARYT_t^+ - \BARYT_t^-) M^{1/2} \big)V'(x) W_t^\beta e^{-(\beta-t)V}{\mathrm{d}} t\\
&\quad+\frac{1}{4}\int_0^\beta \lim_{M\to\infty} \big(\BARYT_t^+(W^\beta_t)^T \int_0^1 V''(x+sM^{-1/2}W^\beta_t){\mathrm{d}}s W_t^\beta
e^{-(\beta-t)V}{\mathrm{d}}t\\
&\quad +
\frac{1}{4}\int_0^\beta \lim_{M\to\infty} \big(\BARYT_t^-(W^\beta_t)^T \int_0^1 V''(x-sM^{-1/2}W^\beta_t){\mathrm{d}}s W_t^\beta\big)e^{-(\beta-t)V}{\mathrm{d}}t\big)\Big)\Big]\\
&= \E\Big[e^{-{\IU}W_\beta\cdot p}\int_0^\beta\int_0^t e^{-sV}V'(x)W_s^\beta e^{-(t-s)V} 
V'(x)W_t^\beta e^{-(\beta-t)V}{\mathrm{d}}s{\mathrm{d}}t\Big]\\
&\quad +\frac{1}{2}\E\Big[e^{-{\IU}W_\beta\cdot p} \int_0^\beta e^{-tV} (W_t^\beta)^T V''(x)W_t^\beta e^{-(\beta-t)V}{\mathrm{d}}t\Big]\,.
\end{split}
\end{equation}
To determine the path integrals in the right hand side of \eqref{4.8} we 
make a partition of $[0,\beta]$ into time intervals 
$[t_j,t_{j+1})$, where $t_j=j\beta/J$ for $j=0,\ldots, J$  and corresponding Wiener increments
$\Delta W_j=W(t_{j+1})-W(t_j)$ and time steps $\Delta t:=\beta/J$,
to obtain
\[
W^\beta_{t_k}\equiv \frac{1}{2}W_\beta-W_{t_k}= \sum_{j=0}^{J-1}\frac{1}{2}\Delta W_j-\sum_{j<k}\Delta W_j
=\sum_{j=0}^{J-1}S_{j,k}\Delta W_j
\]
where
\[
S_{j,k}:=\left\{
\begin{array}{c}
\frac{1}{2}\quad  \mbox{if } j\ge k\,,\\
-\frac{1}{2}\quad  \mbox{if } j< k\,.\\
\end{array}\right.
\]
This partition implies
\begin{equation}\label{4.9}
\begin{split}
&-\E[e^{-{\IU}W_\beta\cdot p}
\int_0^\beta e^{-tV} (\WBT)^T
V''(x)\WBT \,e^{-(\beta-t)V}\big)\Rd  t]\\
&=- \lim_{J\to \infty}\E\Big[ e^{-{\IU}\sum_{m=0}^{J-1} \Delta W_m\cdot p}
\sum_{k=0}^{J-1}\sum_{j=0}^{J-1}\sum_{\ell=0}^{J-1} S_{j,k}S_{\ell,k}\Delta W_j^T\underbrace{(e^{-k\Delta t V}V''e^{-(J-k)\Delta t V})}_{=:V''_k}
\Delta W_\ell  \Delta t\Big]\\
&=-\lim_{J\to \infty}
\sum_{k=0}^{J-1}\sum_{j=0}^{J-1}\sum_{\ell=0}^{J-1} \Delta t S_{j,k} S_{\ell,k}\times\\
&\qquad\times \int_{\rset^{N J}}
e^{-{\IU}\sum_{m=0}^{J-1} \Delta W_m\cdot p}
\Delta W_j^T V''_k \Delta W_\ell 
\prod_{n=0}^{J-1}\frac{e^{-|\Delta W_n|^2/(2\Delta t)}}{(2\pi\Delta t)^{N/2}} \Rd (\Delta W_n)\\
&= \lim_{J\to \infty}e^{-J|p|^2\Delta t/2}
\sum_{k=0}^{J-1} (\frac{J}{2}-k)^2\Delta t^2
p^T V''_k p \,\Delta t\\
&=e^{-\beta |p|^2/2}\int_0^\beta(\frac{\beta}{2}-t)^2 e^{-tV(x)}p^T V''(x)p\, e^{-(\beta-t)V(x)}\Rd t=\BIGO (1)\,.
\end{split}
\end{equation}
Similarly we have
\[
\begin{split}
&-\E\Big[e^{-{\IU}W_\beta\cdot p}\int_0^\beta\int_0^t e^{-sV}V'(x)W_s^\beta e^{-(t-s)V} 
V'(x)W_t^\beta e^{-(\beta-t)V}{\mathrm{d}}s{\mathrm{d}}t\Big]\\
&=-\lim_{J\to\infty}
\sum_{k=0}^{J-1}\sum_{r=0}^{J-1}\sum_{j=0}^{J-1}\sum_{\ell=0}^{J-1} S_{j,r}S_{\ell,k}(\Delta t)^2
e^{-{\IU}\sum_{m=0}^{J-1} \Delta W_m\cdot p} 
\times\\
&\qquad\times
e^{-r\Delta t V} V'\Delta W_j e^{(k-r)\Delta t V}
V'\Delta W_\ell e^{-(J-k)\Delta t V}
\prod_{n=0}^{J-1}\frac{e^{-|\Delta W_n|^2/(2\Delta t)}}{(2\pi\Delta t)^{N/2}} \Rd (\Delta W_n)\\
&=\lim_{J\to\infty}
\sum_{k=0}^{J-1}\sum_{r=0}^{J-1}(\Delta t)^4 (\frac{J}{2}-r)(\frac{J}{2}-k) e^{-J\Delta t |p|^2/2}
e^{-r\Delta t V}V'pe^{-(k-r)\Delta tV} V'p e^{-(J-k)\Delta t V}\\
&=e^{-\beta |p|^2/2}\int_0^\beta\int_0^t (\frac{\beta}{2}-s)(\frac{\beta}{2}-t) e^{-sV(x)}V'(x)pe^{-(t-s)V(x)} V'(x)p
e^{-(\beta-t)V(x)}{\mathrm{d}}s{\mathrm{d}}t
\end{split}
\]
so that $(\rho-e^{-\beta H})=\mathcal O(M^{-1})$. 

The construction \eqref{4.4} implies 
\[
\frac{{\mathrm{d}}}{{\mathrm{d}}t}\TR\big((\BARYT_t^+)^2\big)=
\TR\Big(\BARYT_t^+\big(-2V(x+M^{-1/2}W^\beta_t)\big)\BARYT_t^+\Big)
\]
and by assumption there is a constant $k$ such that
$V+k\Id$ is positive definite everywhere. Therefore we have
\[
\frac{{\mathrm{d}}}{{\mathrm{d}}t}\TR\big((\BARYT_t^+)^2\big)\le k\TR\big((\BARYT_t^+)^2\big)
\]
which establishes $\TR\big((\BARYT_\beta^+)^2\big)\le e^{k\beta}$ and  shows that for independent Wiener processes $W$ and $W'$ we obtain by Cauchy's inequality 
\[
\begin{split}
\TR\rho^2&=\TR\big(\E[e^{-{\mathrm{i}} W_\beta\cdot p} \BARYT_\beta^+(W)]^*\E[e^{-{\mathrm{i}} W'_\beta\cdot p} \BARYT_\beta^+(W')]\big)\\
&= \E[e^{{\mathrm{i}} (W_\beta-W'_\beta)\cdot p}\TR\big(\BARYT_\beta^+(W)\BARYT_\beta^+(W')\big)]\\
&\le \E[ \Big(\TR\big((\BARYT_\beta^+(W))^2\big)\TR\big((\BARYT_\beta^+(W'))^2\big)\Big)^{1/2}]\\
&\le \E[\TR\Big(\big(\BARYT_\beta^+(W)\big)^2\Big)]\\
&\le e^{k\beta}\,.
\end{split}
\]
We also observe that
the generalized Weyl's law implies that $\rho$ is in $L^2(\rset^{2N}, \cset^{d\times d})$, namely 
\[
\begin{split}
(\frac{M^{1/2}}{2\pi})^N\int_{\rset^{2N}} \TR\big(\rho^2(z)\big){\mathrm{d}}z 
=\TR(\OPER{\rho}\OPER{\rho})= \TR(e^{-\beta\OPER{H}}e^{-\beta\OPER{H}})
=\TR(e^{-2\beta\OPER{H}})<\infty\,,
\end{split}
\]
which proves \eqref{symb_estim_bis}.
\end{proof}

\subsection{Proof Lemma \ref{lemma1}}
\begin{proof}

To estimate remainder terms we will use the composition operator for Weyl symbols
defined by $ \widehat{c\MP d} = \widehat c\widehat d$.
The composition operator has the representation
\begin{equation}\label{composition}
\begin{split}
c\MP d&= e^{\frac{{\IU}}{2M^{1/2}} (\nabla_{x'}\cdot\nabla_{p}- \nabla_{x}\cdot\nabla_{p'})}
c(x,p)d(x',p')\big|_{(x,p)=(x',p')}\\
&=e^{\frac{{\IU}}{2M^{1/2}} (\nabla_{x''}\cdot\nabla_{p'}- \nabla_{x'}\cdot\nabla_{p''})}
\underbrace{c(x+x',p+p')d(x+x'',p+p'')}_{=:f_{xp}(x',p',x'',p'')}\big|_{(x',p')=(x'',p'')=0}
\end{split}
\end{equation}
which can be written as an expansion using the Fourier transform $\FT$,  %
defined for $f:\rset^{4N}\to\rset$ by
\begin{equation}\label{fourier_def}
\FT \{f\}(\xi_{x'}, \xi_{p'},\xi_{x''}, \xi_{p''}):= %
\int_{\rset^{4N}}f(x',p',x'',p'')\EXP{-{\IU}(x'\cdot\xi_{x'}+p'\cdot\xi_{p'}+x''\cdot \xi_{x''} +p''\cdot \xi_{p''})} \Rd x'\Rd p'\Rd x''\Rd p''\,.
\end{equation}
 Its inverse transform implies
\[
\begin{split}
&\EXP{-\frac{\IU}{2M^{1/2}}(\nabla_{x''}\cdot\nabla_{p'}- \nabla_{x'}\cdot\nabla_{p''})}f(x',p',x'',p'')\Big|_{(x',p')=(x'',p'')=0}\\
&= %
(\frac{1}{2\pi})^{4N} 
\int_{\rset^{4N}} \FT  f(\xi_{x'},\xi_{p'},\xi_{x''},\xi_{p''}) \EXP{\frac{\IU}{2}M^{-1/2}
(\xi_{x''}\cdot\xi_{p'}-\xi_{x'}\cdot\xi_{p''})} \Rd\xi_{x'}\Rd\xi_{p'} \Rd\xi_{x''}\Rd\xi_{p''}\\
\end{split}
\]
and Taylor expansion of the exponential function  yields
\begin{equation}\label{exp_def}
\begin{split}
&\EXP{-\frac{\IU}{2M^{1/2}} (\nabla_{x''}\cdot\nabla_{p'}- \nabla_{x'}\cdot\nabla_{p''}) }f_{xp}(x',p',x'',p'')\Big|_{(x',p')=(x'',p'')=0}\\
&= (\frac{1}{2\pi})^{4N} 
\int_{\rset^{4N}} \FT  f_{xp}(\xi_{x'},\xi_{p'},\xi_{x''},\xi_{p''}) \EXP{\frac{\IU}{2}M^{-1/2}
(\xi_{x''}\cdot\xi_{p'}-\xi_{x'}\cdot\xi_{p''})} \Rd\xi_{x'}\Rd\xi_{p'} \Rd\xi_{x''}\Rd\xi_{p''}\\
%
&= %
(\frac{1}{2\pi})^{4N}
\int_{\rset^{2N}} \FT  f_{xp}(\xi_{x'},\xi_{p'},\xi_{x''},\xi_{p''}) \Big(\sum_{n=0}^m
(\frac{\IU(\xi_{x''}\cdot\xi_{p'}-\xi_{x'}\cdot\xi_{p''})}{2M^{1/2}})^n \frac{1}{n!} \\
&\quad +(\frac{\IU(\xi_{x''}\cdot\xi_{p'}-\xi_{x'}\cdot\xi_{p''})}{2M^{1/2}})^{m+1} \times\\
&\qquad \times \frac{1}{m!} 
\int_0^1(1-s)^m \EXP{\frac{\IU s}{2}M^{-1/2}(\xi_{x''}\cdot\xi_{p'}-\xi_{x'}\cdot\xi_{p''})} 
\Rd s\Big)\Rd\xi_{x'}\Rd\xi_{p'} \Rd\xi_{x''}\Rd\xi_{p''} \\
&= 
\sum_{n=0}^m\frac{1}{n!} (-\frac{\IU (\nabla_{x''}\cdot\nabla_{p'}- \nabla_{x'}\cdot\nabla_{p''})
}{2M^{1/2}})^n  f_{xp}(x',p',x'',p'')\Big|_{(x,'p')=(x'',p'')=0} \\
&\quad + (\frac{1}{2M^{1/2}})^{m+1}
\int_0^1  \EXP{-\frac{\IU s}{2}M^{-1/2} (\nabla_{x''}\cdot\nabla_{p'}- \nabla_{x'}\cdot\nabla_{p''})}
(-\IU  (\nabla_{x''}\cdot\nabla_{p'}- \nabla_{x'}\cdot\nabla_{p''}))^{m+1} \\
&\qquad\times f_{xp}(x',p',x'',p'')  \frac{(1-s)^m}{ m! }\Rd s\Big|_{(x,'p')=(x'',p'')=0}\PERIOD\\
\end{split}
\end{equation}
The pointwise limit of the remainder term can be estimated by dominated convergence
\[
\begin{split}
\lim_{M\to\infty}&\int_0^1  \EXP{-\frac{\IU s}{2}M^{-1/2} (\nabla_{x''}\cdot\nabla_{p'}- \nabla_{x'}\cdot\nabla_{p''})}
(-\IU  (\nabla_{x''}\cdot\nabla_{p'}- \nabla_{x'}\cdot\nabla_{p''}))^{m+1} \\
&\qquad\times f_{xp}(x',p',x'',p'')  \frac{(1-s)^m}{ m! }\Rd s\Big|_{(x,'p')=(x'',p'')=0}\\
&=(-\IU  (\nabla_{x''}\cdot\nabla_{p'}- \nabla_{x'}\cdot\nabla_{p''}))^{m+1}  f_{xp}(x',p',x'',p'')  \frac{1}{ (m+1)! }\Big|_{(x,'p')=(x'',p'')=0}\\
\end{split}
\]
provided $\int_{\rset^{4N}}| (\xi_{x''}\cdot\xi_{p'}-\xi_{x'}\cdot\xi_{p''})^{m+1} 
\FT f_{xp}(\xi)|\Rd \xi<\infty$. In addition we need convergence in $L^1(\rset^{2N})$,
as a function of $z=(x,p)$, to apply dominated convergence in the phase-space integrals. We have
\[
\FT f_{xp}(\xi_{x'},\xi_{p'},\xi_{x''},\xi_{p''})=\FT c(\xi_{x'},\xi_{p'})e^{{\IU}(x\cdot\xi_{x'}+p\cdot\xi_{p'})}
\FT d(\xi_{x''},\xi_{p''})e^{{\IU}(x\cdot\xi_{x''}+p\cdot\xi_{p''})}
\]
so that
\[
\begin{split}
&\int_{\rset^{4N}} \FT f_{xp}(\xi_{x'},\xi_{p'},\xi_{x''},\xi_{p''})\EXP{\frac{\IU s}{2}M^{-1/2}(\xi_{x''}\cdot\xi_{p'}-\xi_{x'}\cdot\xi_{p''})} \Rd\xi_{x'}\Rd\xi_{p'} \Rd\xi_{x''}\Rd\xi_{p''}\\
&= \int_{\rset^{4N}} (1-\Delta_{\xi_{x'}}-\Delta_{\xi_{p'} })^k 
(1-\Delta_{\xi_{x''}}-\Delta_{\xi_{p''} })^k  \Big(\FT c(\xi_{x'},\xi_{p'}) \FT d(\xi_{x''},\xi_{p''})\times\\
&\qquad\times \EXP{\frac{\IU s}{2}M^{-1/2}(\xi_{x''}\cdot\xi_{p'}-\xi_{x'}\cdot\xi_{p''})}\Big) \times\\
&\qquad\times 
\frac{e^{{\IU}(x\cdot\xi_{x''}+p\cdot\xi_{p''})}}{(1+|x''|^2+|p''|^2)^k}
\frac{e^{{\IU}(x\cdot\xi_{x'}+p\cdot\xi_{p'})}}{(1+|x'|^2+|p'|^2)^k}\Rd\xi_{x'}\Rd\xi_{p'} \Rd\xi_{x''}\Rd\xi_{p''}\,.
\end{split}
\]
Therefore we obtain remainder terms that are uniformly bounded in $L^1(\rset^{2N})$ provided
the Fourier transform $({\IU}\xi)^{\alpha}\FT c(\xi) $ of $\partial_z^\alpha c(z)$ satisfies
\begin{equation}\label{schwartz}
\int_{\rset^{2N}} |(1-\Delta_\xi)^{N+1}\big(\xi^{\alpha}\FT c(\xi)\big)|\Rd \xi<\infty
\end{equation}
and similarly for $d$
\[
\int_{\rset^{2N}} |(1-\Delta_\xi)^{N+1}\big(\xi^{\alpha}\FT d(\xi)\big)|\Rd \xi<\infty\,.
\]
We will apply the composition expansion to functions in the Schwartz class so that \eqref{schwartz} holds.

We conclude that Schwartz functions $c$ and $d$ satisfy
\begin{equation}\label{sharp2}
\begin{split}
&\lim_{M\to\infty}\big(M(\frac{ c\MP d+d\MP  c}{2}-cd)\big) \\
&= \lim_{M\to\infty}\Big( - \frac{1}{4}\int_0^1\cos(\frac{s}{2M^{1/2}}\nabla_{z}\cdot\nabla_{z'}')
(\nabla_{z}\cdot\nabla_{z'}')^2c(z)
d(z')\big|_{z=z'}(1-s)\Rd s\Big)\\
&= - \frac{1}{8}
(\nabla_{z}\cdot\nabla_{z'}')^2c(z)
d(z')\big|_{z=z'}\,,
\end{split}
\end{equation}
and 
\begin{equation}\label{sharp1}
\begin{split}
&\lim_{M\to\infty}( c\MP d)\\
&= \lim_{M\to\infty}\Big( cd - \frac{{\mathrm{i}}}{2M^{1/2}}\int_0^1e^{\frac{{\mathrm{i}} s}{2M^{1/2}}\nabla_{z}\cdot\nabla_{z'}'}
(\nabla_{z}\cdot\nabla_{z'}')c(z)
d(z')\big|_{z=z'}\Rd s\Big)\\
&= cd \,,
\end{split}
\end{equation}
as limits in $L^1(\rset^{2N})$ and in $L^\infty(\rset^{2N})$.
We also have
\[
\widehat{Da_s}={\IU}M^{1/2}[\OPERW{H},\OPERW{a_s}] -\{h,a_s\}^{\widehat{}}= \{H-h,a_s\}^{\widehat{}}\
+\hat r^a_s
\]
where
\begin{equation}\label{ras}
r^a_s =\frac{1}{8M}\int_0^1\cos(\frac{\bar s}{2M^{1/2}}\nabla_{z}\cdot\nabla_{z'}')
(\nabla_{z}\cdot\nabla_{z'}')^3H(z)a_s(z')\big|_{z=z'}(1-\bar s)^2\Rd \bar s
\end{equation}
so that
\[
\lim_{M\to\infty} Mr^a_s=\frac{1}{24}(\nabla_{z_0}\cdot\nabla_{z_0'}')^3H(z_0)a_0\big(z_s(z_0')\big)\big|_{z_0=z_0'}
\]
and the Poisson bracket takes the form
\[
\begin{split}
&\{H-h,a_s\}=\NABLAR_{z_0'}\cdot\nabla_{z_0}
(H-h)(z'_0)a_0\big(z_s(z_0)\big)\big|_{z_0=z_0'}\\
&=\big(\nabla_{p_0'}\cdot\nabla_{x_0}-\nabla_{x_0'}\cdot\nabla_{p_0}\big)
(H-h)(x'_0,p'_0)a_0\big(z_s(x_0,p_0)\big)\big|_{(x_0,p_0)=(x_0',p_0')}
\,.
\end{split}
\]
Based on this remainder estimate there holds
\[
\begin{split}
\TR\big((\widehat{Da_s})^2\widehat\rho\big)
&=\NORMFAC^N\int_{\rset^{2N}}\TR\Big( \big({\IU} M^{1/2}(H\MP a_s-a_s\MP H)
-\{h,a_s\}\big)
\\&\qquad
\MP  \big({\IU} M^{1/2}(H\MP a_s-a_s\MP H)-\{h,a_s\}\big)\ \rho\Big)\Rd z
\end{split}\]
and in the limit we obtain by \eqref{sharp1}, \eqref{ras} and Lemma \ref{path_lemma}
\[
\begin{split}
&\lim_{M\to\infty}\int_{\rset^{2N}}\TR\Big( \big({\IU} M^{1/2}(H\MP a_s-a_s\MP H)
-\{h,a_s\}\big)
\\&\qquad
\MP  \big({\IU} M^{1/2}(H\MP a_s-a_s\MP H)-\{h,a_s\}\big)\ \rho\Big)\Rd z\\
&=\int_{\rset^{2N}}\TR\Big(\big(\nabla_{z_0}a(z_s(z_0)\cdot \NABLAR_{z_0}(H-h)\big)^2e^{-\beta H(z_0)}\Big)\Rd z_0\,, 
\end{split}\]
using splitting of the phase space integral as in \eqref{Mrho}, which together with \eqref{sharp2}, \eqref{sharp1}, \eqref{ras} and Lemma \ref{a'} below proves the lemma.
\end{proof}

\begin{lemma}\label{a'}
Assume that the bounds in Theorem \ref{thm_mf} hold, then $a_t$ and $b_t$ are in the Schwartz class and
there is a constant $C'$, depending in $C$, such that
\begin{equation}\label{ab_bound}
\begin{split}
    \sum_{|\alpha|\le 3}\|\partial_z^\alpha a_t(z)\|_{L^\infty(\rset^{2N})}
    +\sum_{|\alpha|\le 3}\|\partial_z^\alpha b_t(z)\|_{L^\infty(\rset^{2N})} &\le C'\,.\\
\end{split}
\end{equation}
\end{lemma}

\begin{proof}
To estimate $\sum_{|\alpha|\le 3}\|\partial_{p_0}^\alpha a_0\big(z_t(z_0)\big)\|_{L^\infty(\rset^{2N})}$ we use the first order flow $\nabla_{z_0} z_t(z_0)=:z'(t)$, second order flow 
$z''_{,km}(t)=\partial_{z_k(0) z_m(0)} z(t)$ and third order flow $z'''(t)$, which
are solutions to the system
\[
\begin{split}
\dot z_i(t) &= (\nabla_{z_t}' h(z_t))_i=:f_i(z_t)\COMMA \\
z'_{i,k}(t) &= \Id_{ik}+ \int_0^t\sum_{k'} f'_{i,k'}(z_s)z'_{k',k}(s)\Rd s\COMMA \quad f'_{i,k'}(z):=\partial_{z_{k'}}f_i(z) \COMMA\\
z''_{i,km}(t) &=  \int_0^t \Big(\sum_{k'}  f'_{i,k'}(z_s)z''_{k',k m}(s) +  \sum_{k'm'} f''_{i,k'm'}(z_s)z'_{k',k}(s)z'_{m',m}(s)\Big) \Rd s\COMMA\\
&\qquad f''_{i,k'm'}(z):=\partial_{z_{k'}z_{m'}}f_i(z)\COMMA\\
 z'''_{i,kmn}(t) &=  \int_0^t \Big(\sum_{k'} f'_{i,k'}(z_s)z'''_{k',k mn}(s) 
+  \sum_{k'm'} f''_{i,k'm'}(z_s)z'_{k',k}(s)z''_{m',mn}(s) \\
&\qquad +\sum_{k'm'} f''_{i,k'm'}(z_s)z''_{k',kn}(s)z'_{m',m}(s)\\
&\qquad +\sum_{k'n'} f''_{i,k'n'}(z_s)z''_{k',k m}(s)z'_{n',n}(s)\\
&\qquad+\sum_{k'm'n'} f'''_{i,k'm'n'}(z_s)z'_{k',k}(s)z'_{m',m}(s)z'_{n',n}(s)\Big)\Rd s
\,.\\
\end{split}
\]
By summation and maximization over indices we obtain the integral inequalities 
\begin{equation}\label{z_flow_est}
\begin{split}
\max_{ik}
|z'_{i,k}(t)| &\le 1+ \int_0^t \sum_{k'} |f'_{i,k'}(z_s)| \max_{ik}|z'_{i,k}(s)|\Rd  s\COMMA\\
\max_{ik}\sum_{m}
|z''_{i,km}(t)| &\le \int_0^t \sum_{k'} |f'_{i,k'}(z_s)| \max_{ik}\sum_{m}
|z''_{i,km}(s)|\Rd  s\\
&\qquad + \int_0^t \sum_{k'm'} |f''_{i,k'm'}(z_s)| (\max_{ik}|z'_{i,k}(s)|)^2\Rd  s\COMMA\\
\max_{ik}\sum_{mn}
|z'''_{i,kmn}(t)| &\le \int_0^t \sum_{k'} |f'_{i,k'}(z_s)| \max_{ik}\sum_{mn}
|z'''_{i,kmn}(s)|\Rd  s\\
&\quad + \int_0^t \sum_{k'm'} |f''_{i,k'm'}(z_s)| \max_{ik}|z'_{i,k}(s)| \max_{ik}\sum_m
|z''_{i,km}(s)|\Rd  s\\
&\quad + \int_0^t \sum_{k'm'n'} |f'''_{i,k'm'n'}(z_s)| (\max_{ik}|z'_{i,k}(s)|)^3\Rd  s
\PERIOD\\
\end{split}
\end{equation}
The functions $\max_{ij}\sum_{|\alpha|\le 2} \partial_{z_0}^\alpha \partial_{z_j}z_i(t,z_0)$ can therefore be estimated as in \cite{gronwall} by Gronwall's inequality, which states: if there is a positive constant $K$ and  continuous positive functions $\gamma, u:[0,\infty)\rightarrow [0,\infty)$ such that  
\[u(t)\le K + \int_0^t \gamma(s)u(s)\Rd  s, \quad \mbox{ for } t>0\COMMA\]
 then 
 \[ u(t)\le K\EXP{\int_0^t \gamma(s)\Rd  s}, \quad \mbox{ for } t>0\PERIOD\]
Gronwall's inequality applied to \eqref{z_flow_est} implies
\begin{equation}\label{z_tz_0}
\max_{ij}\sum_{|\alpha|\le 2} \|\partial_{z_0}^\alpha \partial_{z_j}z_i(t,z_0)\|_{L^\infty(\rset^{2N})}=\BIGO(1)
\end{equation}
provided that
\begin{equation}\label{R7'}
\begin{split}
\max_i\sum_{|\alpha|\le 3}\|\partial^\alpha_x \partial_{x_i}h\|_{L^\infty(\rset^{N})} &= \BIGO(1)\PERIOD\\
\end{split}
\end{equation}

The flows $z',z'',z'''$ determine the derivatives of the scalar symbol $a_0\big(z_t(z_0)\big)=a_t(z_0)$
\begin{equation}\label{A_z_t}
\begin{split}
 \partial_{z_k}a_t &= \sum_{k'} a'_{0,k'}(z_t)z'_{k',k}(t)\COMMA\qquad a'_{0,k'}(z):=\partial_{z_{k'}}a_0(z),\\
\partial_{z_kz_m}a_t &=  \sum_{k'} a'_{0,k'}(z_t)z''_{k',k m}(t) +  \sum_{k'm'} a''_{0,k'm'}(z_t)z'_{k',k}(t)z'_{m',m}(t)\COMMA\\
\partial_{z_kz_mz_n}a_t&=  \sum_{k'} a'_{0,k'}(z_t)z'''_{k',k mn}(t) 
+  \sum_{k'm'} a''_{0,k'm'}(z_t)z'_{k',k}(t)z''_{m',mn}(t)\\
&\qquad +\sum_{k'm'} a''_{0,k'm'}(z_t)z''_{k',kn}(t)z'_{m',m}(t)\\
&\qquad +\sum_{k'n'} a''_{0,k'n'}(z_t)z''_{k',k m}(t)z'_{n',n}(t)\\
&\qquad+\sum_{k'm'n'} a'''_{0,k'm'n'}(z_t)z'_{k',k}(t)z'_{m',m}(t)z'_{n',n}(t)
\COMMA\\
\end{split}
\end{equation}
which together with \eqref{z_tz_0} proves \eqref{ab_bound}.

Similarly to verify that $a_t$ is a Schwartz function,  we first
extend both \eqref{z_flow_est} (for $\partial_{z_0}^\alpha z_t(z_0)$) 
and the representation in \eqref{A_z_t} to order $|\alpha|$ to obtain the bounds 
\begin{equation}\label{a_t_bound}
|\partial_z^\alpha a_t(z)|<\infty\,,
\end{equation}
for all indices $\alpha$.
Then define the compact set $\mathbb L_m^c=\{z\in\rset^{2N}\ :\ h(z)\le m\}$. The property $h(z)\to\infty$  as $|z|\to\infty$, which is verified below \eqref{h_infty}, implies that $L_m^c$ includes the support of $a_0$
for $m$ sufficiently large and the invariance $h\big(z_t(z_0)\big)=h(z_0)$, for all $z_0\in\rset^{2N}$, establishes by \eqref{a_t_bound}
\[
\sup_{z\in \rset^{2N}}|z^\gamma\partial_z^\alpha a_t(z)|
=\sup_{z\in L_m^c} |z^\gamma\partial_z^\alpha a_0\big(z_t(z)\big)|<\infty\,,
\]
for all indices $\gamma$ and $\alpha$.
We conclude that $a_t$ is a Schwartz function.
\end{proof}

The constant  in the right hand side of \eqref{z_tz_0} grows typically exponentially with respect to $t$, i.e.
\[
\max_{ij}\sum_{|\alpha|\le 2} \|\partial_{z_0}^\alpha \partial_{z_j}z_i(t,z_0)\|_{L^\infty(\rset^{2N})}
\le \EXP{ c't}\COMMA\]
 where $c'$ is the positive constant in the right hand side of \eqref{R7'}.
 We note that the assumptions on $h$ and $V$ are compatible with the assumption
to have $\int_{\rset^{2N}} \TR e^{-\beta H(z)}{\mathrm{d}}z$ bounded.

\section{Numerical experiments}\label{sec_numerics}
In this section we define a model problem that allows us to systematically
study approximation of canonical quantum correlation observables. 
The model problem is constructed so that we can accurately approximate the quantum dynamics, 
thereby avoiding the computational challenges to accurately approximate 
realistic quantum systems with many particles and excited electron states, 
cf. \cite{cances_nonlinear}. 
Subsection ~\ref{subsection_5_1} includes the approximations of equilibrium density function of position observable $\Hat{x}$, applying the mean-field molecular dynamics and the electron ground state molecular dynamics, respectively. In Subsection \ref{subsection_5_2}, we compare time-dependent correlation observables obtained from molecular dynamics evolving on the ground state,  
on the mean-field energy surface, and on a weighted average of all eigenstates (denoted the excited state dynamics below). 
In particular, we study whether the mean-field approximation can be more accurate than using only the ground state. 

In order to demonstrate the proposed mean-field molecular dynamics approximation, we devise the model problem as described by equations \eqref{H_hat_oper1} to \eqref{two_eigval_new1} in Subsection \ref{subsec_numerics},
where the difference between two eigenvalues $\lambda_1(x)-\lambda_0(x)=2c\,\sqrt{x^2+\delta^2}$ can be tuned by the two parameters $c$ and $\delta$. For $\delta$ small this model relates to the avoided crossing phenomenon in quantum chemistry
 where the two potential surfaces get almost intersected at a certain point, see \cite{teller}. The assumptions in Theorem \ref{thm_mf} are satisfied for positive $\delta$ but not for $\delta=0$ and therefore the approximation error is expected to vary with different $\delta$.

\subsection{Equilibrium observables}\label{subsection_5_1}
At the inverse temperature $\beta$, the quantum canonical ensemble average of a time-independent observable $\Hat{A}$ is obtained from the normalized trace
\begin{equation}\label{def_qm_average}
\frac{\mathrm{Tr}(e^{-\beta\Hat{H}}\Hat{A})}{\mathrm{Tr}(e^{-\beta\Hat{H}})}=\frac{\sum_n e^{-\beta E_n}\langle \Phi_n,\, \Hat{A}\Phi_n\rangle }{\sum_n e^{-\beta E_n}}\,, 
\end{equation}
where $(E_n,\Phi_n)_{n=1}^\infty$ are eigenvalues and the corresponding normalized eigenstates of the Hamiltonian operator $\widehat{H}$.
Consequently, for an observable depending only on the position $x$ with symbol $A(x)$, we have
\begin{equation}\label{mu_qm}
\frac{\mathrm{Tr}(e^{-\beta\Hat{H}}\Hat{A})}{\mathrm{Tr}(e^{-\beta\Hat{H}})}=\int_{\mathbb{R}^N}A(x)\,\mu_{\mathrm{qm}}(x) \,\mathrm{d}x,\ \mbox{ where }\ \mu_{\mathrm{qm}}(x)=\frac{\sum_n |\Phi_n(x)|^2e^{-\beta E_n}}{\sum_n e^{-\beta E_n}}\,.
\end{equation}
We apply a fourth-order finite difference scheme for the Laplacian operator in the Hamiltonian \eqref{H_hat_oper1} to approximate the equilibrium density $\mu_{\mathrm{qm}}(x)$ in \eqref{mu_qm}. The numerical implementation is explained with more details in Appendix~\ref{subsection_5_3_1}.

As an approximation to the quantum canonical ensemble average \eqref{def_qm_average}, we consider the normalized trace
$\TR(\widehat{e^{-\beta H}}\widehat{A})\,/\,{\TR(\widehat{e^{-\beta H}})}$
and apply the Lemma \ref{lemma_composition} and equation \eqref{equilibrium_tau_0_compare} to write the mean-field observable as
\[
\begin{split}
    \mathfrak{T}_{\mathrm{md}}(0)=\mathfrak{T}_{\mathrm{es}}(0)&=\frac{\TR(\widehat{e^{-\beta H}}\widehat{A})}{\TR(\widehat{e^{-\beta H}})}=\frac{\int_{\mathbb{R}^{2N}}\mathrm{Tr}(e^{-\beta H(x,p)}a(x,p)\mathrm{I})\,\mathrm{d}x\,\mathrm{d}p}{\int_{\mathbb{R}^{2N}}\mathrm{Tr}(e^{-\beta H(x',p')})\,\mathrm{d}x'\,\mathrm{d}p'}\\
&=\frac{\int_{\mathbb{R}^{2N}}\,a(x,p)(e^{-\beta\lambda_0(x)}+e^{-\beta\lambda_1(x)})\,e^{-\frac{\beta|p|^2}{2}}\,\mathrm{d}x\,\mathrm{d}p}{\int_{\mathbb{R}^{2N}}\,(e^{-\beta\lambda_0(x')}+e^{-\beta\lambda_1(x')})\,e^{-\frac{\beta|p'|^2}{2}}\,\mathrm{d}x'\,\mathrm{d}p'}\,,
\end{split}
\]
where $H(x,p)$ and $A(x,p)=a(x,p)\mathrm{I}$ with $a:\mathbb{R}^{2N}\to\mathbb{C}$ are the Weyl symbols corresponding to the operators $\Hat{H}$ and $\Hat{A}$, respectively. 
Specifically for an observable depending only on the position $x$, we obtain
\begin{equation}\label{mu_mf}
\frac{\TR(\widehat{e^{-\beta H}}\widehat{A})}{\TR(\widehat{e^{-\beta H}})}=\int_{\mathbb{R}^{N}}A(x)\,\mu_{\mathrm{mf}}(x)\,\mathrm{d}x,\\ \mbox{ where }\ \mu_{\mathrm{mf}}(x)=\frac{e^{-\beta\lambda_0(x)}+e^{-\beta\lambda_1(x)}}{\int_{\mathbb{R}^N}\,(e^{-\beta\lambda_0(x')}+e^{-\beta\lambda_1(x')})\,\mathrm{d}x'}\,. 
\end{equation}
The classical mean-field density $\mu_{\mathrm{mf}}$ in \eqref{mu_mf} can also be rewritten as a weighted average
\begin{equation}\label{mu_cl_weighted_qj}
\mu_{\mathrm{mf}}(x)=\sum_{j=0}^{1} \,q_j\,\frac{e^{-\beta\lambda_j(x)}}{\int_{\mathbb{R}}e^{-\beta\lambda_j(x')}\mathrm{d}x'}\,,\;\;\mbox{where }\;q_j=\frac{\int_{\mathbb{R}^N}\,e^{-\beta\lambda_j(x)}\mathrm{d}x}{\sum_{k=0}^{1}\int_{\mathbb{R}^N}\,e^{-\beta\lambda_k(x)}\mathrm{d}x}\,,\  j=0,1\,.
\end{equation}
The weights $q_0$ and $q_1$ can be interpreted as the probability for the system to be in the corresponding electron eigenstate $\lambda_0$ and $\lambda_1$, respectively obtained by integration with the corresponding Gibbs density.

We first plot the equilibrium quantum mechanics density $\mu_{\mathrm{qm}}$ using the formula \eqref{mu_qm}, and compare it with the classical mean-field density $\mu_\mathrm{mf}$ in \eqref{mu_mf}. They are also compared with the density based only on the ground state $\mu_{\mathrm{gs}}$, with the formula
\begin{equation}\label{mu_gs}
     \mu_{\mathrm{gs}}(x)=\frac{e^{-\beta\lambda_0(x)}}{\int_{\mathbb{R}^N}e^{-\beta\lambda_0(x)}\mathrm{d}x}\,,
\end{equation}
 which is obtained from the classical density formula \eqref{mu_cl_weighted_qj} by taking the probability for the excited state as $q_1=0$, and the probability for the ground state as $q_0=1$.
\begin{figure}[!htb]
    \begin{center}
  \includegraphics[height=5cm]{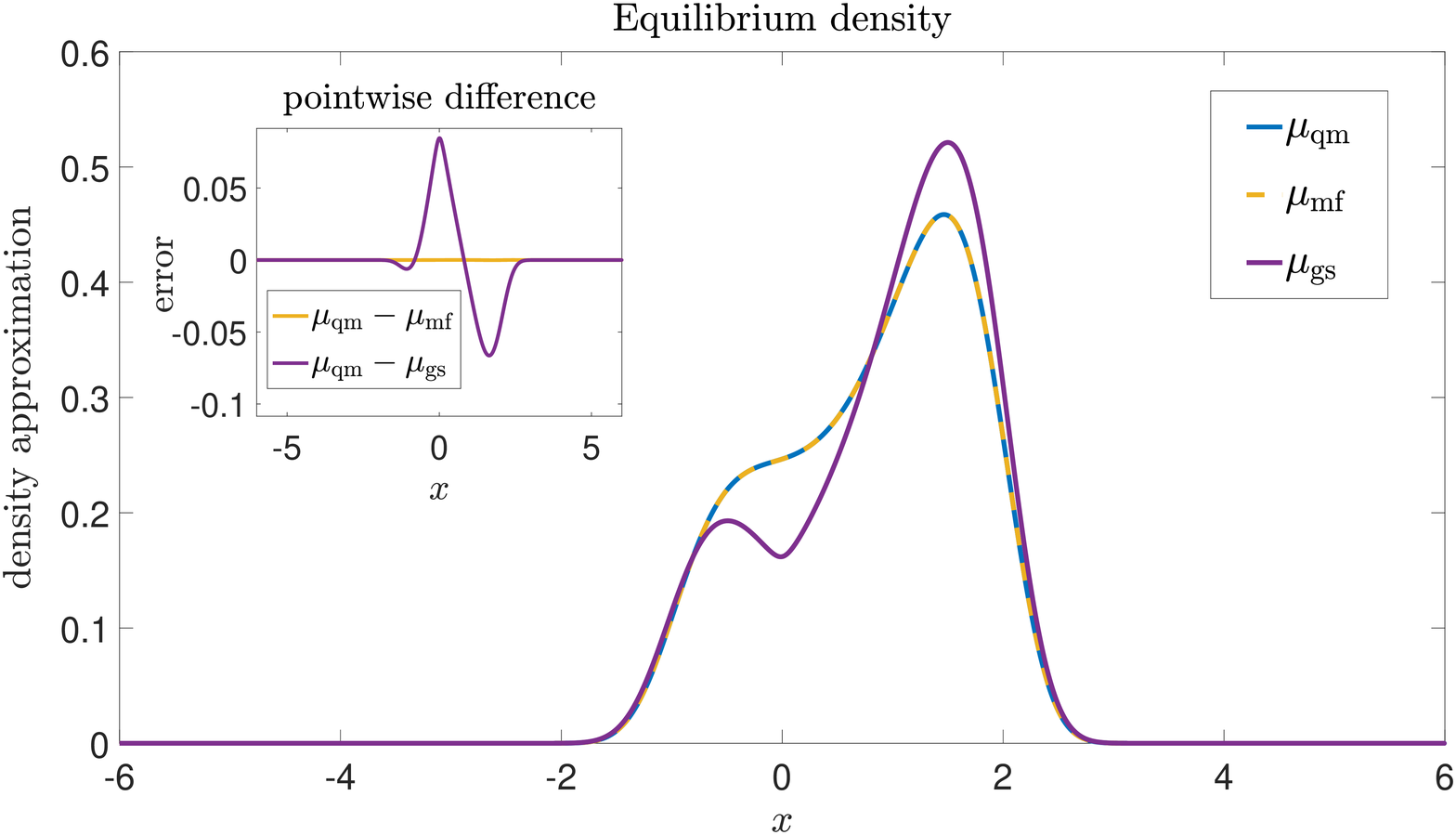}
  \caption{The density functions computed with the quantum mechanics formula \eqref{mu_qm}, the mean-field approximation formula \eqref{mu_mf}, and using only the ground state in the classical formula \eqref{mu_gs}, with mass ratio $M=1000$. The mean-field density (the dashed yellow line) is quite close to the quantum mechanics density curve (the solid blue line), implying a better accuracy than the ground state density (the solid violet line). }\label{fig:three-density-curves}
  \end{center}
\end{figure}

In Figure~\ref{fig:three-density-curves} the first reference density curve with quantum mechanics formula \eqref{mu_qm} is plotted in blue colour with a solid line. 
The density curve $\mu_{\mathrm{mf}}(x)$, obtained from the classical mean-field formula \eqref{mu_mf}, is plotted as the yellow dashed line and agrees well with the quantum mechanics density $\mu_{\mathrm{qm}}(x)$. Besides, the mean-field density $\mu_{\mathrm{mf}}(x)$ incurs much smaller error than the ground-state density $\mu_{\mathrm{gs}}(x)$ (the violet solid curve) in approximating the $\mu_{\mathrm{qm}}(x)$ density.
For Figure~\ref{fig:three-density-curves}, we use the parameters $M=1000$, $c=1$, and $\beta=1$ such that with the 
eigenvalue gap $\delta=0.1$, the system has a probability $q_1=0.16$ to be in the excited state.

\begin{figure}[!htb]
    \begin{center}
  \includegraphics[height=5cm]{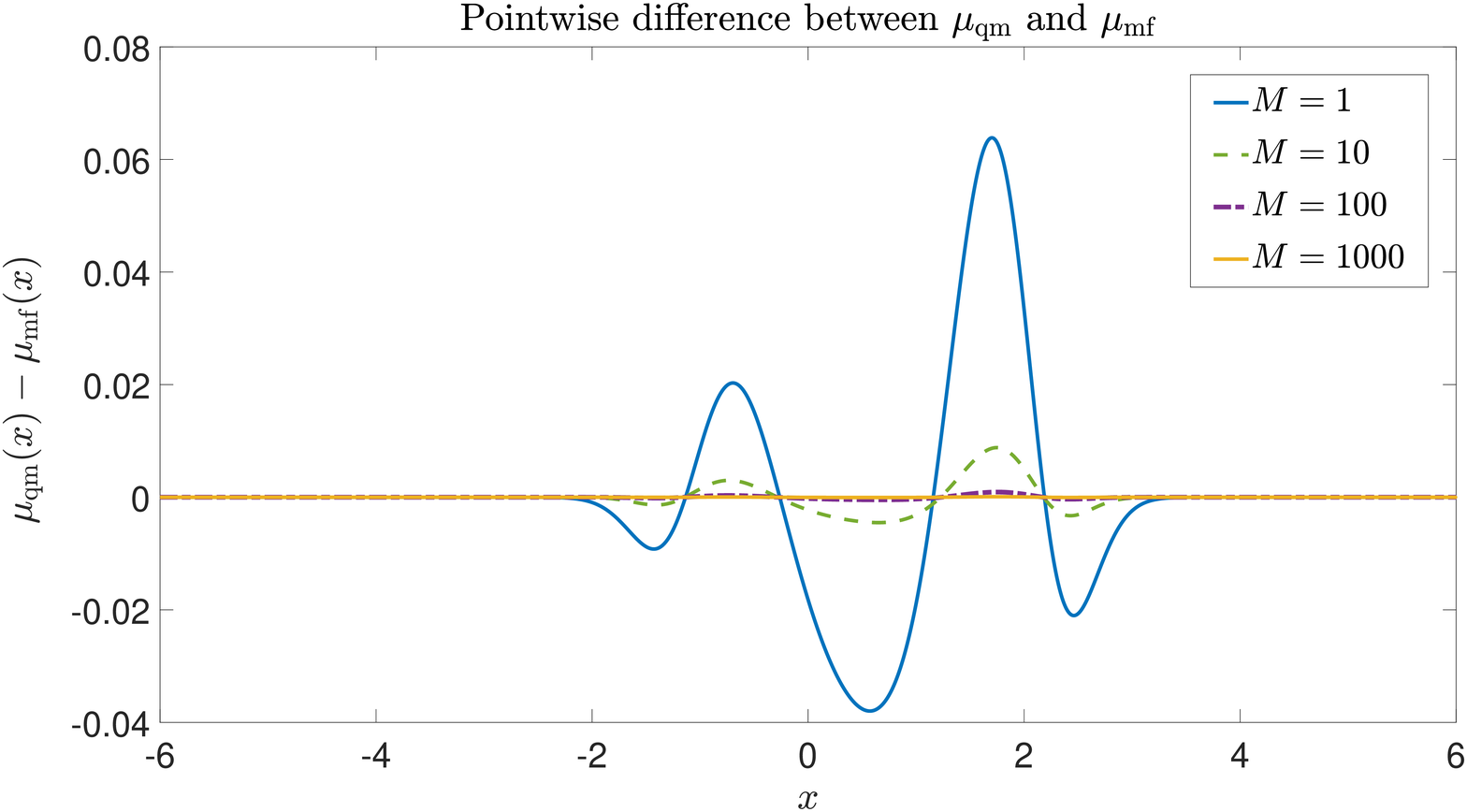}
  \caption{The point-wise difference between the quantum mechanics density $\mu_{\mathrm{qm}}$ and the mean-field density $\mu_{\mathrm{mf}}$ with inverse temperature $\beta=1$. The dashed violet curve with $M=100$ has so small an error that it is almost indiscernible from the solid yellow curve with $M=1000$. }\label{fig:density_error_with_M}
  \end{center}
\end{figure}
In Figure~\ref{fig:density_error_with_M} we depict a point-wise difference between the classical mean-field density $\mu_{\mathrm{mf}}(x)$ and the quantum mechanics 
density $\mu_{\mathrm{qm}}(x)$, 
for different values of the mass ratio $M$. 
The inverse temperature is still taken as $\beta=1$ for the eigenvalue gap $\delta=0.1$ and $c=1$, so that the probability for the excited state is kept as $q_1=0.16$. It is observed from Figure~\ref{fig:density_error_with_M} that as $M$ increases the error in the classical mean-field density approximation decreases. 

\begin{figure}[!htb]
  \begin{center}
  \includegraphics[height=5cm]{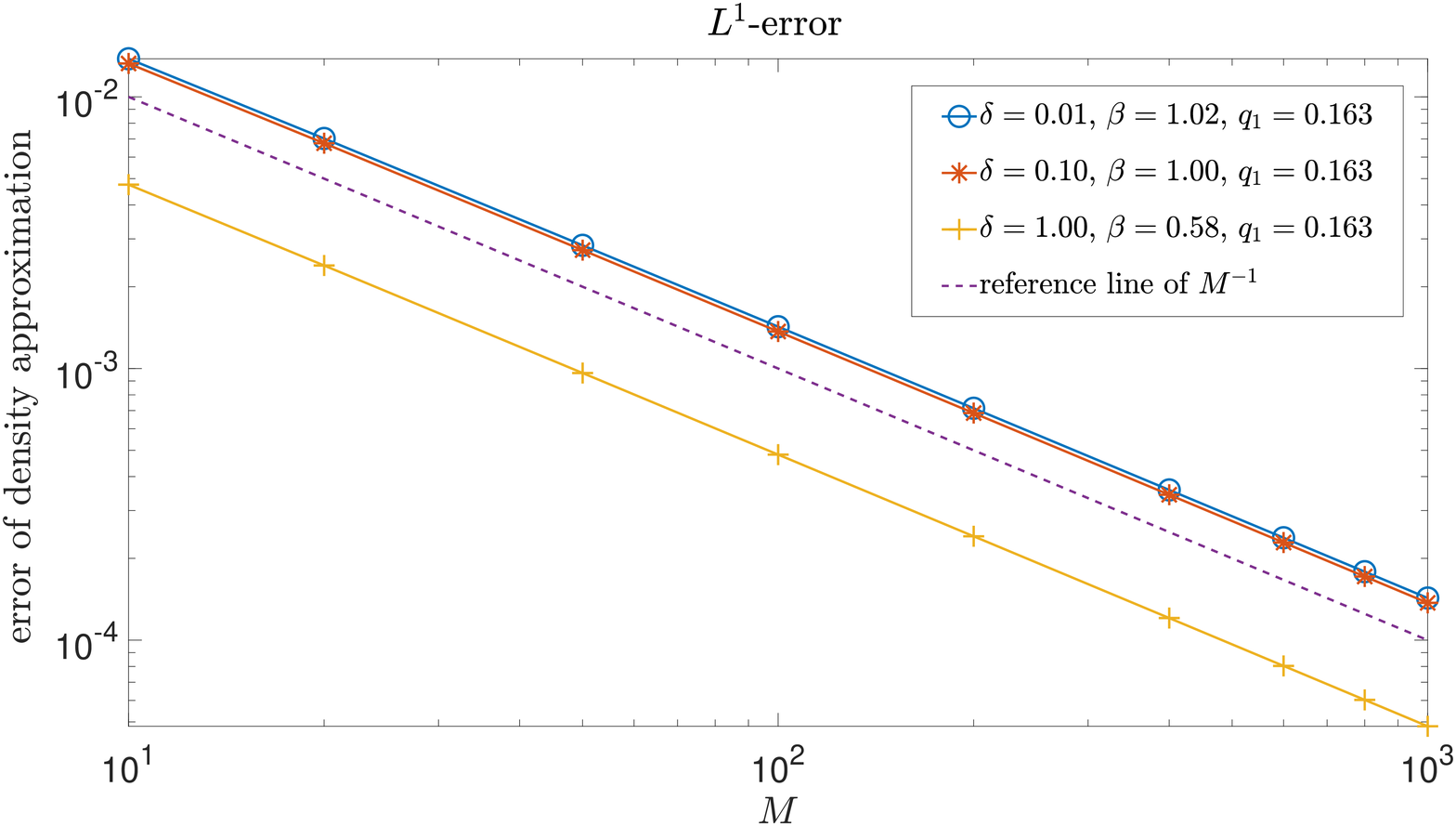}
  \caption{Dependence of the $L^1$-error between the quantum density $\mu_{\mathrm{qm}}$ and the classical mean-field density approximation $\mu_\mathrm{mf}$, shown in log-log scale.}\label{fig:density_error_with_M_and_T}
  \end{center}
\end{figure}

In order to study the dependence of the approximation error $\|\mu_{\mathrm{qm}}-\mu_{\mathrm{mf}}\|_{L^1}$ on the mass ratio $M$, we vary $M$ for three different inverse temperatures, with the corresponding eigenvalue gaps 
$\delta$ such that the probability to be in the excited state remains to be $q_1=0.16$. 
As seen from Figure~\ref{fig:density_error_with_M_and_T}, the $\mathcal{O}(M^{-1})$ dependence of the error in the equilibrium density using the classical mean-field 
approximation is in accordance with the theoretical result of Theorem~\ref{thm_mf}. 

\begin{figure}[!htb]
  \begin{center}
  \includegraphics[height=5cm]{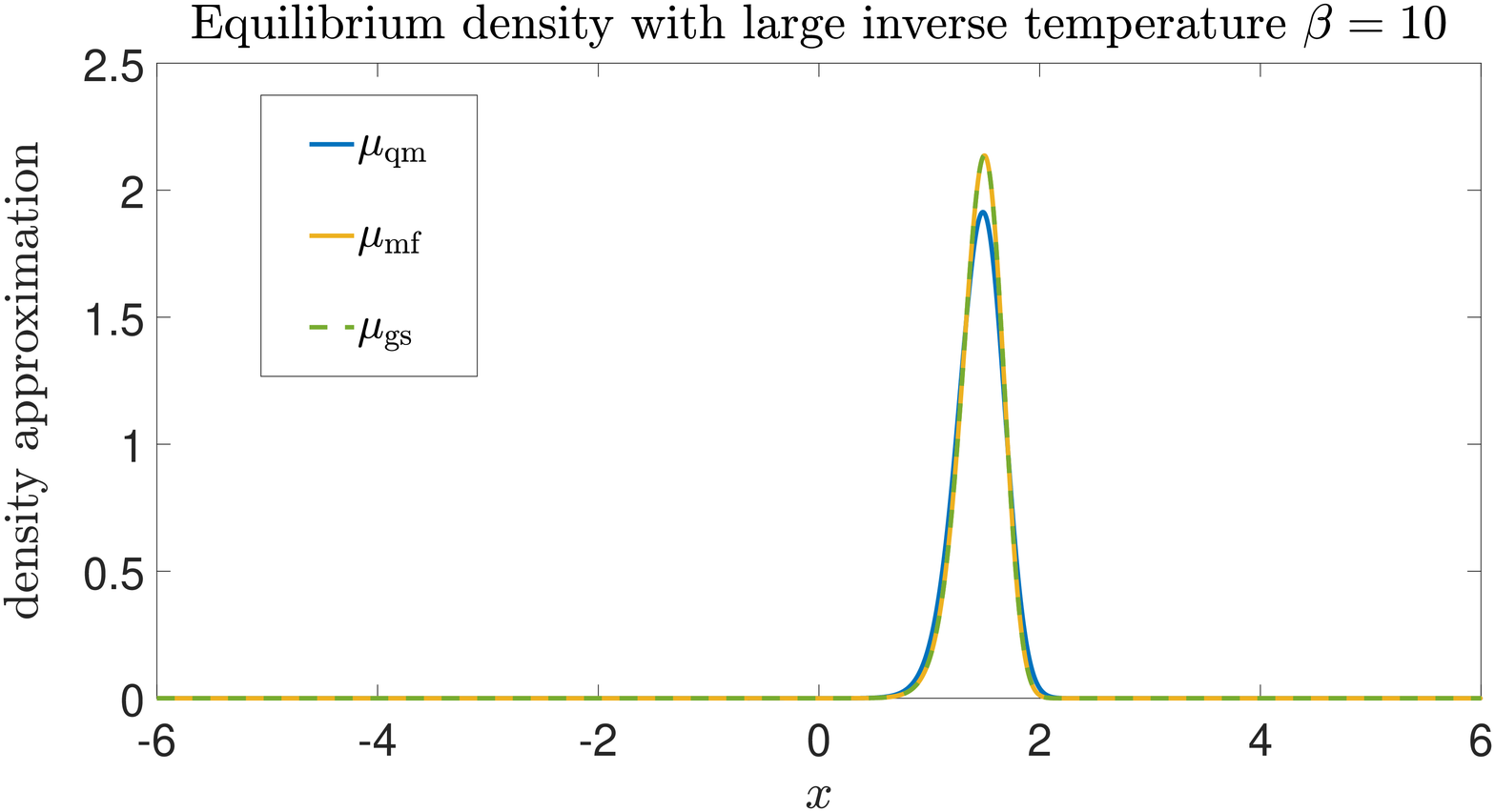}
  \caption{Equilibrium density $\mu_{\mathrm{qm}}$ with the classical mean-field and ground state approximations $\mu_{\mathrm{mf}}$ and $\mu_{\mathrm{gs}}$, with inverse temperature $\beta=10$, mass ratio $M=100$. The probability of electron excited state $q_1=7\times 10^{-7}$ is tiny.}\label{fig:density_large_beta}
  \end{center}
\end{figure}

Besides the $M$-dependence of the classical approximation, we also experiment with a relatively large inverse temperature $\beta = 10$ for mass ratio $M=100$, with parameters $c=1$, $\delta=0.1$. The quantum density $\mu_{\mathrm{qm}}$ together with its classical mean-field and ground state approximations $\mu_{\mathrm{mf}}$ and $\mu_{\mathrm{gs}}$ are plotted in the Figure~\ref{fig:density_large_beta}. The large value of $\beta$ implies a rather low temperature, which leads to a tiny probability for the electron excited state as $q_1=7\times 10^{-7}$. Consequently the density functions concentrate near the minimum of the ground state eigenvalue, and there is almost no difference between the mean-field and the ground state density curves.
\subsection{Time-correlated observables}\label{subsection_5_2}
\subsubsection{The model problem}\label{subsection_5_2_1}
We apply the mean-field molecular dynamics to approximate the auto-correlation function between the momentum observables  $\widehat{p}_0$ (at time $0$) and $\widehat{p}_\tau$ (at time $\tau$). 
In the Heisenberg representation the time evolution of the momentum 
observable is given by
\[ 
\widehat{p}_\tau:=e^{\mathrm{i}\tau\sqrt{M}\widehat{H}}\widehat{p}_0 e^{-\mathrm{i}\tau\sqrt{M}\widehat{H}}\,. 
\]
We study the two-eigenvalue model with the potential matrix $V(x)$ 
as defined by \eqref{pot_mat_V_def1} in Subsection~\ref{subsec_numerics}.
For computing the quantum correlation function, we approximate the initial position observable $\widehat{x}_0$ and the initial momentum observable $\widehat{p}_0$ by the matrices
\[ 
\widehat{x}_0\simeq
\left[\begin{smallmatrix}
x_0 &   &   &   &   &   & \mathbf{0} \\
  & x_0 &   &   &   &   &   \\
  &   & x_1 &   &   &   &   \\
  &   &   & x_1 &   &   &   \\
  &   &   &   & \ddots &   &   \\
  &   &   &   &   & x_K &   \\
\mathbf{0}  &   &   &   &   &   & x_K \\
\end{smallmatrix}\right]
=:X,\quad
\widehat{p}_0\simeq \mathrm{i}\sqrt{M}(H_dX-XH_d)=:P,
\]
respectively, where we discretize a sufficiently large computational domain $\Omega=[x_0,x_K]$ with uniform grid points $x_k=x_0 + k\Delta x$, for $k=0,1,\cdots,K$ and $\Delta x=\frac{|x_K-x_0|}{K}$. 
In the definition of the matrix $P$, the real symmetric matrix $H_d$ is of size $2(K+1)\times 2(K+1)$, corresponding to a fourth order finite difference approximation of the Hamiltonian operator $\widehat{H}$. 
More details about this approximation are provided in Section~\ref{subsection_5_3_1},
and the definition of matrix $H_d$ is given in \eqref{Hd-matrix}. 
The matrix $H_d$ generates the approximations
\[ 
e^{\mathrm{i}\tau\sqrt{M}\widehat{H}}\simeq e^{\mathrm{i}\tau\sqrt{M}H_d},\quad e^{-\beta\widehat{H}}\simeq e^{-\beta H_d}\,.  
\]

We apply the \texttt{eig} function of \texttt{Matlab} to obtain the eigenpairs $(e_n,\phi_n)$ of the $H_d$ matrix, and rearrange them to obtain the eigendecomposition
\[ 
Q:=[ \phi_1\quad \phi_2\quad\phi_3\quad\cdots\quad\phi_{2(K+1)}],\ \textup{ and }\ D:=\left[\begin{smallmatrix}
e_1 &   &   & \mathbf{0} \\
  & e_2 &   &   \\
  &   & \ddots &  \\
 \mathbf{0} &   &   & e_{2(K+1)} \\
\end{smallmatrix}\right]\,,
\]
such that $H_d$ can be diagonalized with the orthogonal matrix $Q$ as $H_d=Q\,D\,Q^\ast$,
and hence
\[ 
e^{\mathrm{i}\tau\sqrt{M}H_d}=Q\,e^{\mathrm{i}\tau\sqrt{M}D}\,Q^\ast,\ \textup{ and }\ e^{-\beta H_d}=Q\,e^{-\beta D}\,Q^\ast\,. 
\]
Thus the right-hand side of \eqref{qms} with $\widehat{A}_\tau=\widehat{p}_\tau$ and $\widehat{B}_0=\widehat{p}_0$ can be approximated by
\begin{equation}\label{T_qm_def}
    \begin{aligned}
    \mathfrak{T}_{\mathrm{qm}}(\tau)&:= 
    \frac{\TR\big((\widehat{p}_\tau\, \widehat{p}_0 +\widehat{p}_0\,\widehat{p}_\tau) \,e^{-\beta \widehat H}\big) }{2\,\TR( e^{-\beta \widehat H})}
    =\frac{\TR\big(\widehat{p}_\tau\,( \widehat{p}_0\,e^{-\beta\widehat{H}} + e^{-\beta\widehat{H}}\,\widehat{p}_0 ) \big)}{2\,\TR\big(e^{-\beta\widehat{H}}\big)}\\
    &\simeq \frac{\TR\big( Q\,e^{\mathrm{i}\tau\sqrt{M}D}\,Q^\ast\,P\,Q\,e^{-\mathrm{i}\tau\sqrt{M}D}\,Q^\ast\,(P\,Q\,e^{-\beta D}\,Q^\ast+Q\,e^{-\beta D}\,Q^\ast\,P) \big)}{2\,\TR\big( Q\,e^{-\beta D}\,Q^\ast \big)}\,,
    \end{aligned}
\end{equation}
where in the second equality we use the cyclic property of the trace.

By applying the mean-field molecular dynamics formula \eqref{md_mf} with momentum observables $\widehat{p}_0$ and $\widehat{p}_\tau$ in the model problem, 
we have the approximation for the time-correlation function as
\begin{equation}\label{mean-field_md_formula}
\mathfrak{T}_{\mathrm{mf}}(\tau)=\frac{\int_{\rset^{2N}} p_\tau p_0
\,\TR( e^{-\beta H(z_0)})\,{\mathrm{d}}z_0}{\int_{\rset^{2N}}  \TR( e^{-\beta H(z_0)})\,{\mathrm{d}}z_0}\,,
\end{equation}
where $z_t:=(x_t,p_t)$ solves the Hamiltonian system
\begin{equation}
\begin{aligned}
\dot{x}_t&=\nabla_p \,h(x_t,p_t)\,,\\
\dot{p}_t&=-\nabla_x \,h(x_t,p_t)\,,\\
\end{aligned}\label{mf_dynamics}
\end{equation}
with an initial state $z_0=(x_0,p_0)\in\mathbb{R}^{2}$. 
The Hamiltonian system \eqref{mf_dynamics} is solved numerically 
with the second-order velocity Verlet scheme, see \cite{md_book}. 
More details about this numerical implementation is in Appendix~\ref{subsection_5_3_2}.

We also apply the classical molecular dynamics formula for correlation functions introduced in \cite[Section 2.3.2]{KPSS}, which considers the contribution 
from the ground state and the excited states. 
For our specific example the {\it excited state dynamics} approximation 
of momentum correlation observable is given by
\begin{equation}\label{excited_state_formula}
\mathfrak{T}_{\mathrm{es}}(\tau):=\sum_{j=0}^1\int_{\mathbb{R}^2} q_j\,p_\tau^j(z_0)\,p_0^j(z_0)\,\frac{e^{-\beta(\frac{|p_0|^2}{2}+\lambda_j(x_0))}}{\int_{\mathbb{R}^2} e^{-\beta(\frac{|p|^2}{2}+\lambda_j(x))}\mathrm{d}z}\,\mathrm{d}z_0, \,,
\end{equation}
with the weights $q_0$ and $q_1$ as defined in \eqref{mu_cl_weighted_qj}, where $z_\tau^j=(x_\tau^j,p_\tau^j)$, $j=0,1$ solves the Hamiltonian dynamics
\[ 
\begin{aligned}
\dot{x}_\tau^j&=p_\tau^j\,,\\
\dot{p}_\tau^j&=-\nabla\lambda_j(x_\tau^j)\,,
\end{aligned}
\]
with the initial condition $z_0^j=(x_0,p_0)=z_0$ and $\lambda_0(x)$, $\lambda_1(x)$ as defined in \eqref{two_eigval_new1}. 

In addition to these three expressions for time-correlation, $\mathfrak{T}_{\mathrm{qm}}$ (the quantum mechanics correlation), $\mathfrak{T}_{\mathrm{mf}}$ (the mean-field approximation), and $\mathfrak{T}_{\mathrm{es}}$ (the classical excited state approximation), we also compute the approximation based only on the ground state contribution, $\mathfrak{T}_{\mathrm{gs}}$, obtained by setting the probability $q_1$ for the excited state equal to be zero, i.e.,
\begin{equation}\label{(ground-state_formula)}
\mathfrak{T}_{\mathrm{gs}}(\tau)=\int_{\mathbb{R}^2} p_\tau(z_0)\,p_0(z_0)\frac{e^{-\beta(\frac{|p_0|^2}{2}+\lambda_0(x_0))}}{\int_{\mathbb{R}^2} e^{-\beta(\frac{|p|^2}{2}+\lambda_0(x))}\mathrm{d}z}\mathrm{d}z_0\,,
\end{equation}
where $z_\tau=(x_\tau,p_\tau)$ solves the Hamiltonian dynamics with the potential $\lambda_0(x)$
\[ 
\begin{aligned}
\dot{x}_\tau&=p_\tau\,,\\
\dot{p}_\tau&=-\nabla\lambda_0(x_\tau)\,,
\end{aligned}
\]
with initial state $z_0=(x_0,p_0)$.

As discussed in \eqref{equilibrium_tau_0_compare}, the approximations with mean-field dynamics or excited state dynamics 
at an initial time $\tau=0$ are always identical, i.e.,
$\mathfrak{T}_{\mathrm{mf}}(0) = \mathfrak{T}_{\mathrm{es}}(0)$. 
For position-related observables, the classical ground state dynamics approximation $\mathfrak{T}_{\mathrm{gs}}(0)$ will be,
in general, different from $\mathfrak{T}_{\mathrm{mf}}(0)$ and $\mathfrak{T}_{\mathrm{es}}(0)$, 
which is consistent with our preceding observations on equilibrium density function in 
Section~\ref{subsection_5_1}, and is confirmed in the upcoming Figure~\ref{fig:corr_case_A5_xcorr}.

The numerical approximations of  $\mathfrak{T}_{\mathrm{mf}}(\tau)$, $\mathfrak{T}_{\mathrm{es}}(\tau)$, and $\mathfrak{T}_{\mathrm{gs}}(\tau)$ are computed with the Verlet method 
in combination with Simpson's formula.
We use $\mathcal{S}_{\mathrm{mf}}(\tau)$, $\mathcal{S}_{\mathrm{es}}(\tau)$, $\mathcal{S}_{\mathrm{gs}}(\tau)$, and $\mathcal{S}_{\mathrm{qm}}(\tau)$,
to denote the numerical approximations of $\mathfrak{T}_{\mathrm{mf}}(\tau)$, $\mathfrak{T}_{\mathrm{es}}(\tau)$,
$\mathfrak{T}_{\mathrm{gs}}(\tau)$, and $\mathfrak{T}_{\mathrm{qm}}(\tau)$, respectively.

In practise ground state molecular dynamics in the canonical ensemble is relatively well developed for realistic molecular systems and successful 
mean-field approximations appear in centroid and ring-polymer molecular dynamics \cite{marx_hutter}. 
Direct computations of excited state dynamics for realistic systems 
seem less attractive, due to
the challenge to efficiently compute excited electron states, cf. \cite{cances_nonlinear}. 
Here we compare these three alternatives for
a simple model problem in the hope of giving some
information also on realistic systems.

\subsubsection{Numerical results}
Following the discussion on the variances $\epsilon_1^2$ and $\epsilon_2^2$ with equations \eqref{epsilon_2_square_model}, \eqref{epsilon_1_square_model} and \eqref{eps_22}, we survey five different cases with varying parameter settings of $c$, $\delta$, and inverse temperature $\beta$, and make a comparison between the performances of different molecular dynamics approximations in each case. A summary of the parameters in each case is given in the Table~\ref{table_case_A_to_E_new_written_caption}.

    
    

\begin{figure}[!h]
\centering
\begin{subfigure}[b]{0.48\textwidth}
        \centering
        \includegraphics[width=\textwidth]{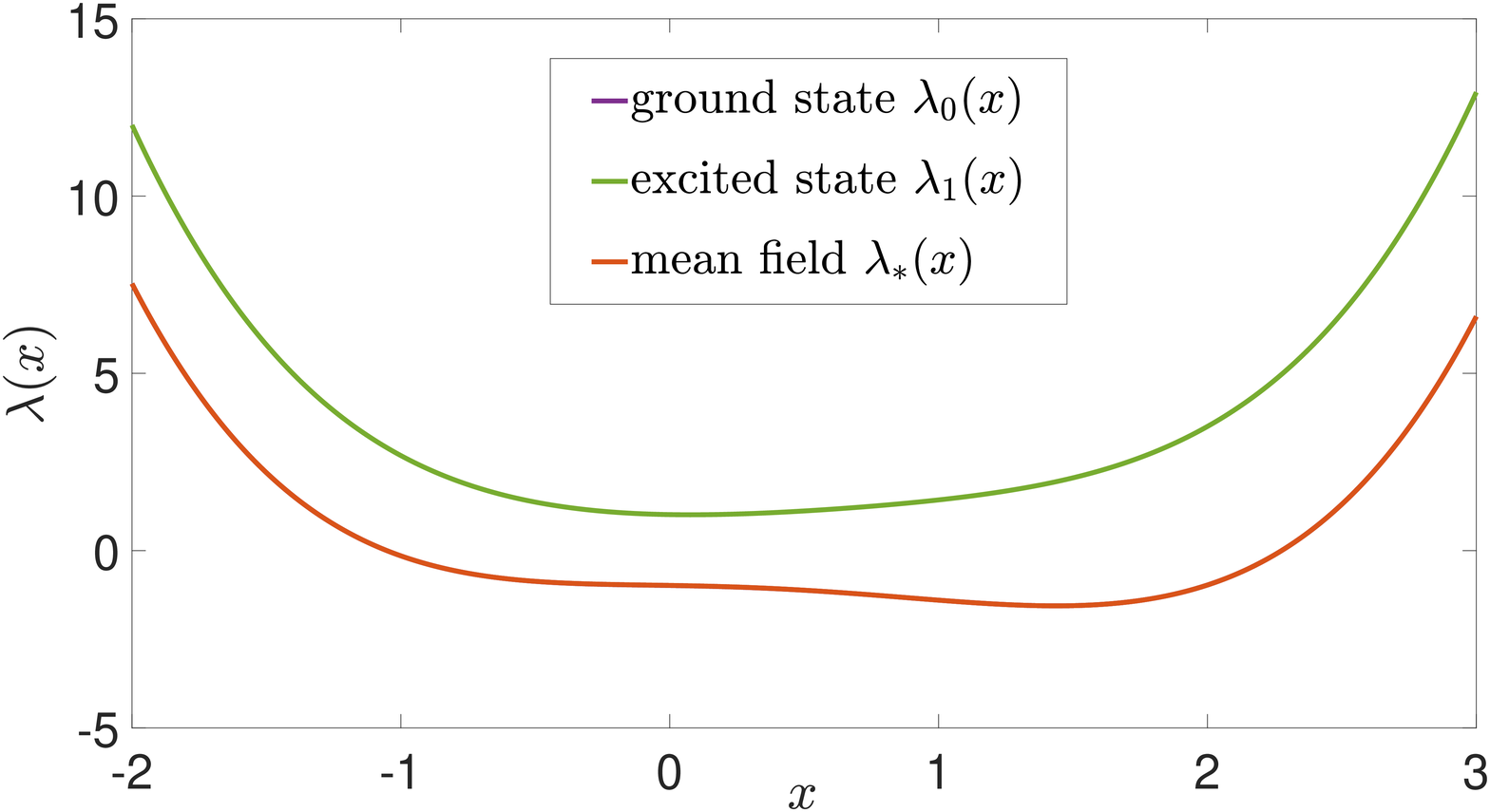}
        \caption{Case A: low temperature, large difference
        \\between eigenvalues, large gap }
        \label{fig:eigval_plot_case_A}
\end{subfigure}
\begin{subfigure}[b]{0.48\textwidth}
        \centering    
        \includegraphics[width=\textwidth]{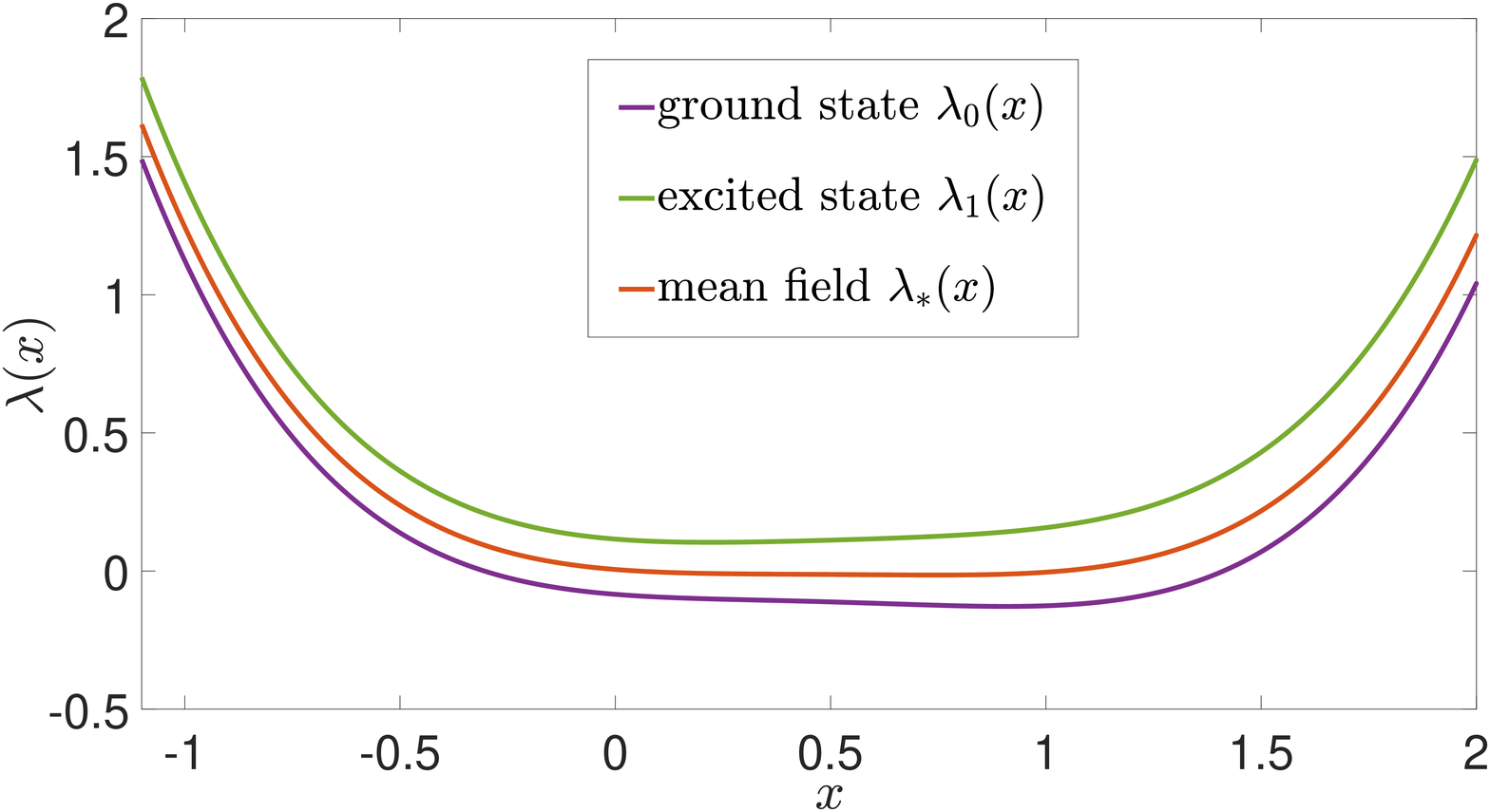}
        \caption{Case B: high temperature, small difference \\between eigenvalues, small gap}
        \label{fig:eigval_plot_case_B}
\end{subfigure}

\begin{subfigure}[b]{0.48\textwidth}
        \centering
        \includegraphics[width=\textwidth]{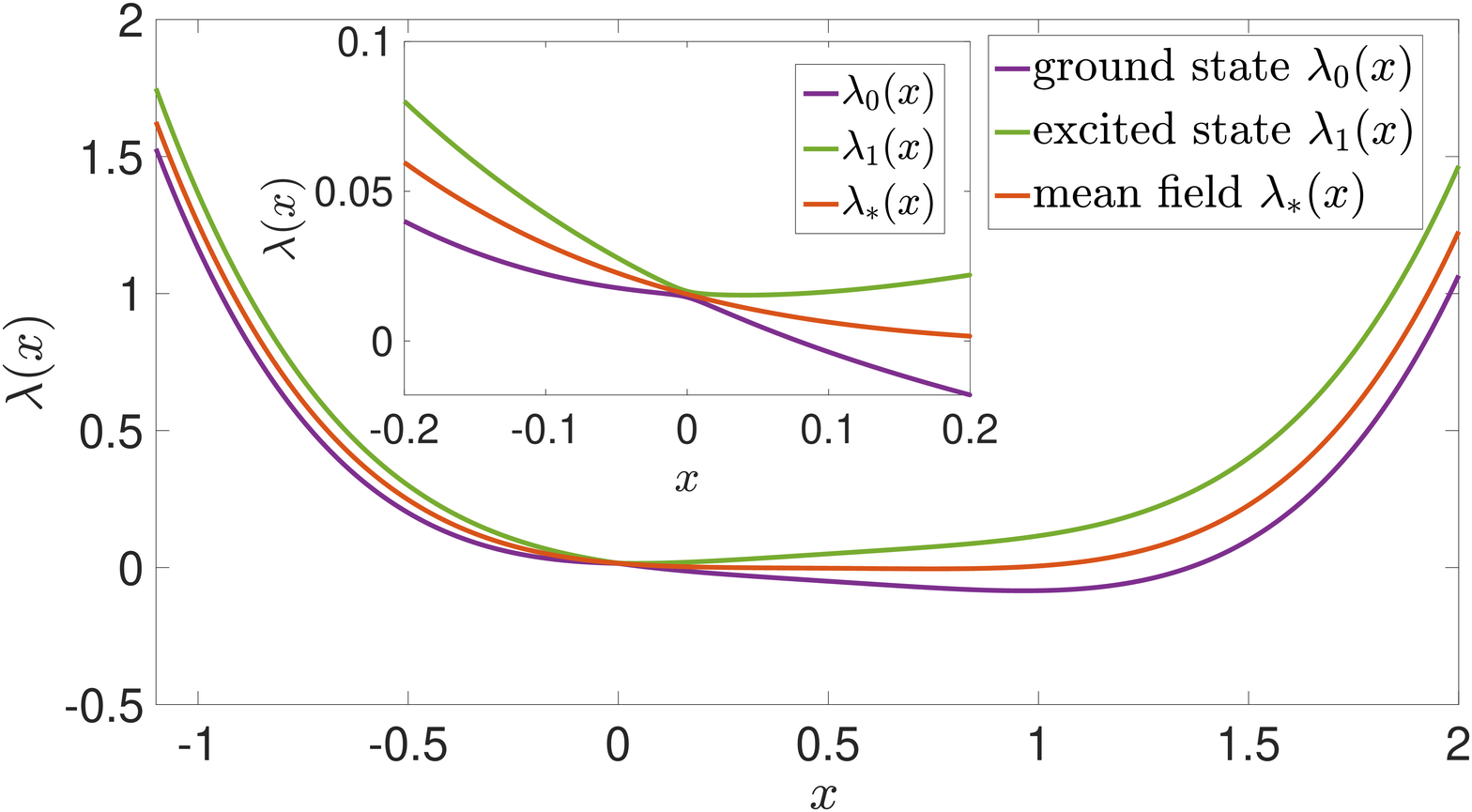}
        \caption{Case C: high temperature, small difference \\between eigenvalues, smallest gap}
        \label{fig:eigval_plot_case_C}
\end{subfigure}
\begin{subfigure}[b]{0.48\textwidth}
        \centering
        \includegraphics[width=\textwidth]{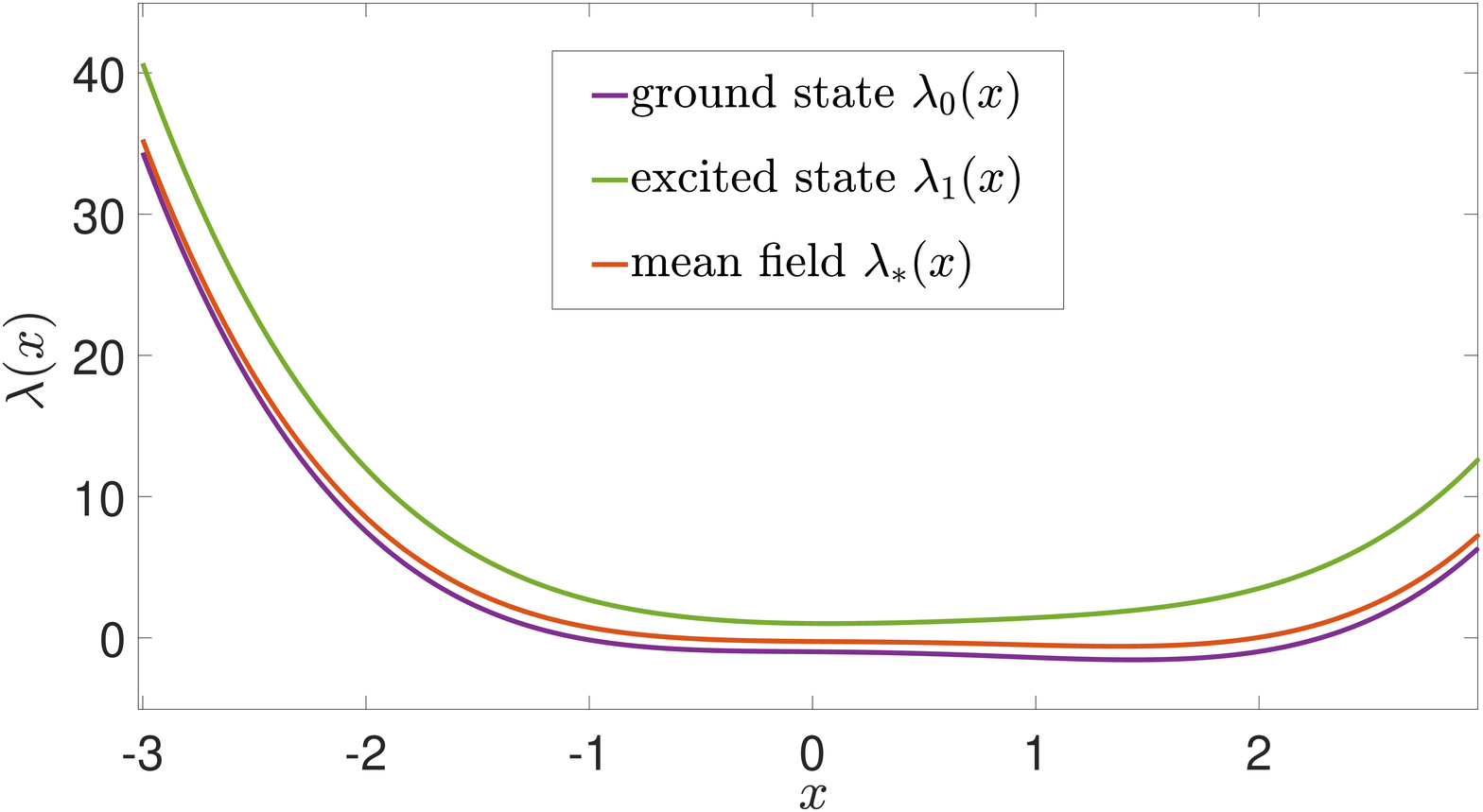}
        \caption{Case D: high temperature, large difference \\between eigenvalues, medium gap}
        \label{fig:eigval_plot_case_D}
\end{subfigure}

\caption{Eigenvalues of the matrix-valued potential $V(x)$ for test cases A to D. In subfigure (a), the mean-field potential $\lambda_\ast(x)$ (the red curve) is quite close to the ground state $\lambda_0(x)$ (the violet curve). }
\label{fig:eigvals-plots}
\end{figure}

\noindent\hypertarget{case_A}{\textbf{Case A}}\textbf{: \textit{Low temperature with large eigenvalue gap}}, $\beta=3.3$, $c=1$, $\delta=1$, $\epsilon_1^2=9.95\times 10^{-4}$, $\epsilon_2^2=1.25\times 10^{-4}$, the probability for the excited state $q_1=0.0002$ is almost negligible.

Figure~\ref{fig:eigval_plot_case_A} presents the eigenvalues $\lambda_0(x)$, $\lambda_1(x)$ and the mean-field potential function $\lambda_\ast(x)$ as defined in \eqref{h_simple} for Case A. With the parameters $c=1$ and $\delta=1$, the system has a large eigenvalue gap. Particularly in this low temperature setting, the mean-field potential $\lambda_\ast(x)$ is almost identical to the ground-state eigenvalue $\lambda_0(x)$. Since the probability for the excited state is very small ($q_1=0.0002$), the three molecular dynamics approximations $\mathfrak{T}_{\mathrm{mf}}(\tau)$, $\mathfrak{T}_{\mathrm{es}}(\tau)$, and $\mathfrak{T}_{\mathrm{gs}}(\tau)$ are similar.

\begin{figure}[h]
     \centering
     \begin{subfigure}[b]{0.49\textwidth}
         \centering
         \includegraphics[width=\textwidth]{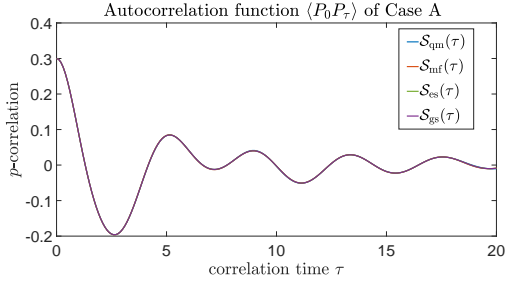}
         \caption{Case A: Auto-correlation curves}
         \label{fig:corr_case_A2}
     \end{subfigure}
     \hfill
     \begin{subfigure}[b]{0.50\textwidth}
         \centering
         \includegraphics[width=\textwidth]{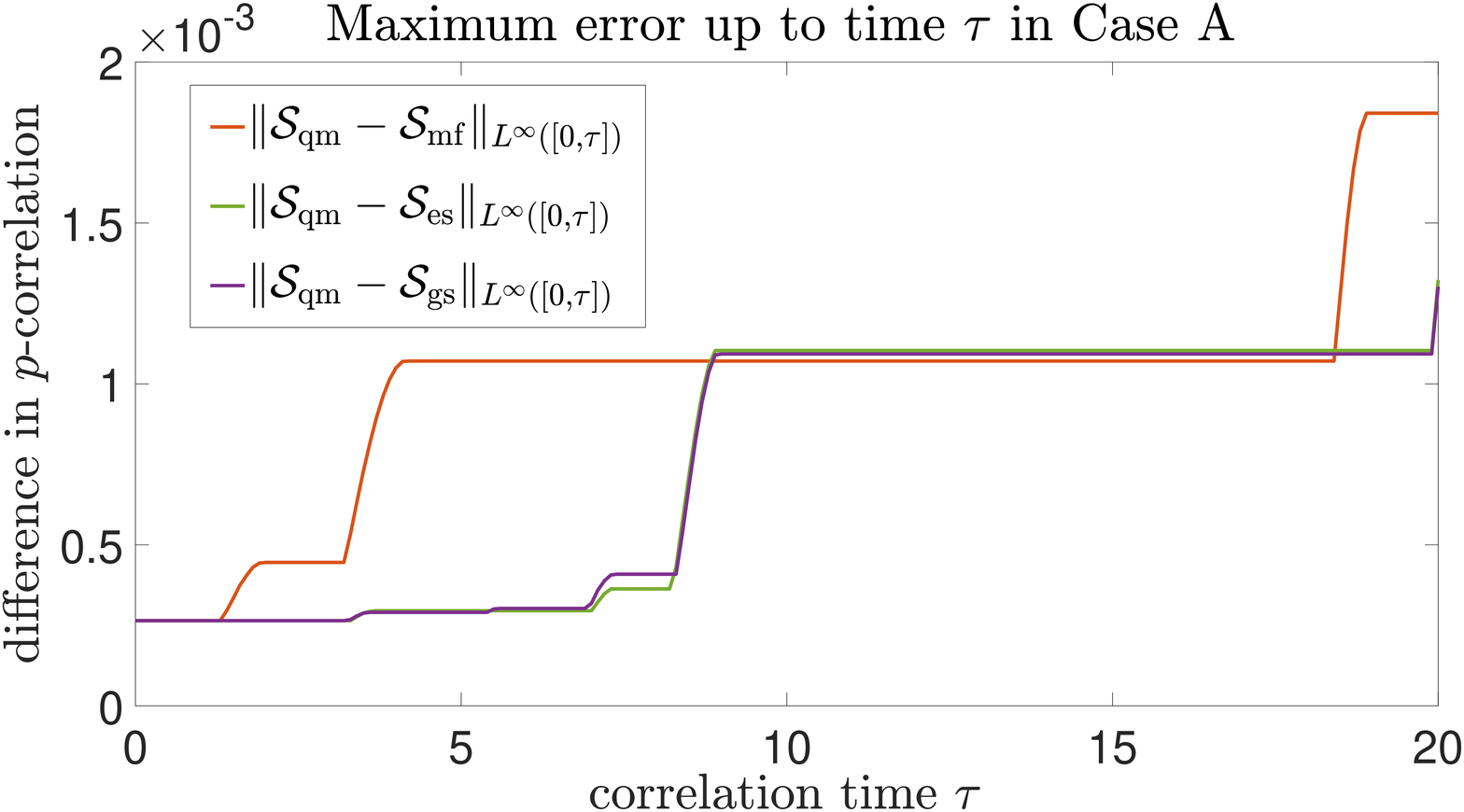}
         \caption{Case A: Maximum error curves up to time $\tau$}
         \label{fig:err_case_A2}
     \end{subfigure}
     \caption{Case A: (a) Auto-correlation function $\langle P_0P_\tau\rangle$ computed by quantum-mechanics formula, $\mathcal{S}_{\mathrm{qm}}$, with $M=1000$, 
     and by three molecular dynamics formulae. (b) The corresponding maximum errors up to time $\tau$, $\|\mathcal{S}_{\mathrm{qm}}-\mathcal{S}_{\mathrm{md}}\|_{L^\infty([0,\tau])}$.}
\end{figure}



In Figure~\ref{fig:corr_case_A2}, the quantum mechanics correlation function curve $\mathcal{S}_{\mathrm{qm}}(\tau)$ with mass ratio $M=1000$ 
is plotted as a function of correlation time $\tau$, together with the three molecular dynamics approximations $\mathcal{S}_{\mathrm{mf}}(\tau)$, $\mathcal{S}_{\mathrm{es}}(\tau)$, and $\mathcal{S}_{\mathrm{gs}}(\tau)$. The three molecular dynamics correlation function curves are almost on top of each other, as shown in Figure~\ref{fig:err_case_A2} with similarly small errors.
This case gives an example where all the three molecular dynamics work analogously, since $q_1\ll 1$, and we note that the error terms $\epsilon_1^2$, $\epsilon_2^2$ together with $1/M$ are all very small.

\noindent\hypertarget{case_B1}{\textbf{Case B}}\textbf{: \textit{High temperature with small difference between eigenvalues}}, $\beta=1$, $c=0.1$, $\delta=1$, $\epsilon_1^2=1.9\times 10^{-2}$, $\epsilon_2^2=3.5\times 10^{-3}$, the probability for the excited state $q_1=0.43$.

\begin{figure}[h]
     \centering
     \begin{subfigure}[b]{0.49\textwidth}
         \centering
         \includegraphics[width=\textwidth]{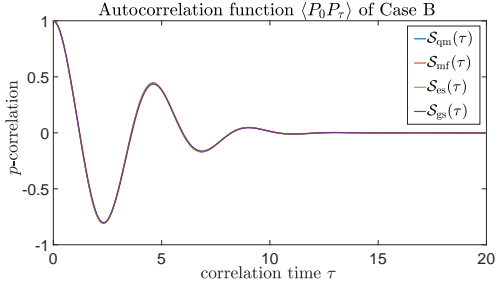}
         \caption{Case B: Auto-correlation curves}
         \label{fig:corr_case_B1}
     \end{subfigure}
     \hfill
     \begin{subfigure}[b]{0.49\textwidth}
         \centering
         \includegraphics[width=\textwidth]{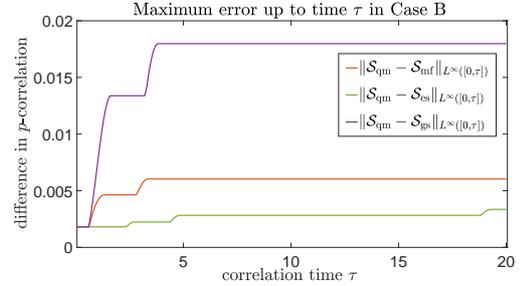}
         \caption{Case B: Maximum error curves up to time $\tau$}
         \label{fig:err_case_B1}
     \end{subfigure}
     \caption{Case B: (a) Auto-correlation function $\langle P_0P_\tau\rangle$ computed by quantum-mechanics formula, $\mathcal{S}_{\mathrm{qm}}$, with $M=100$, and by three molecular dynamics formulae. (b) The corresponding maximum errors up to time $\tau$ $\|\mathcal{S}_{\mathrm{qm}}-\mathcal{S}_{\mathrm{md}}\|_{L^\infty([0,\tau])}$.}
\end{figure}

In Figure~\ref{fig:eigval_plot_case_B} we observe that with the parameter setting of case B, the mean-field potential $\lambda_\ast(x)$ lies in between the ground state eigenvalue $\lambda_0(x)$ and the excited state eigenvalue $\lambda_1(x)$, indicating that by incorporating the effect of the excited state, the mean-field approximation $\mathfrak{T}_{\mathrm{mf}}(\tau)$ can make a difference from simply using the ground state molecular dynamics.

The improved accuracy of the mean-field molecular dynamics is verified by the Figure~\ref{fig:err_case_B1}, in which we  observe a smaller error of mean-field molecular dynamics approximation $\|\mathcal{S}_{\mathrm{qm}}-\mathcal{S}_{\mathrm{mf}}\|_{L^\infty([0,\tau])}$ (the red curve) compared with the molecular dynamics using only the ground state $\|\mathcal{S}_{\mathrm{qm}}-\mathcal{S}_{\mathrm{gs}}\|_{L^\infty([0,\tau])}$ (the violet curve). The excited state molecular dynamics $\mathcal{S}_{\mathrm{es}}(\tau)$ has the smallest error, manifesting an effective combination of the information from both the ground and the excited eigenstates.


\noindent\hypertarget{case_C}{\textbf{Case C}}\textbf{: \textit{High temperature, small difference between eigenvalues with avoided crossing}}, $\beta=1$, $c=0.1$, $\delta=0.01$, $\epsilon_1^2=9.1\times 10^{-3}$, $\epsilon_2^2=9.9\times 10^{-3}$, the probability for the excited state $q_1=0.46$.

The Case C has a similar parameter setting as the preceding \hyperlink{case_B1}{Case B}, with the only difference of a smaller parameter $\delta=0.01$. The small parameter $\delta$ leads to a small eigenvalue gap at $x=0$, i.e., the two eigenvalues $\lambda_0(x)$ and $\lambda_1(x)$ almost intersect at this point, as can be seen in the Figure~\ref{fig:eigval_plot_case_C}. Compared with \hyperlink{case_B1}{Case B}, the small eigenvalue gap also makes the probability for the excited state $q_1$ increase from 0.43 to 0.46 in Case C, with the same inverse temperature $\beta=1$.

\begin{figure}[!htb]
     \centering
     \begin{subfigure}[b]{0.49\textwidth}
         \centering
         \includegraphics[width=\textwidth]{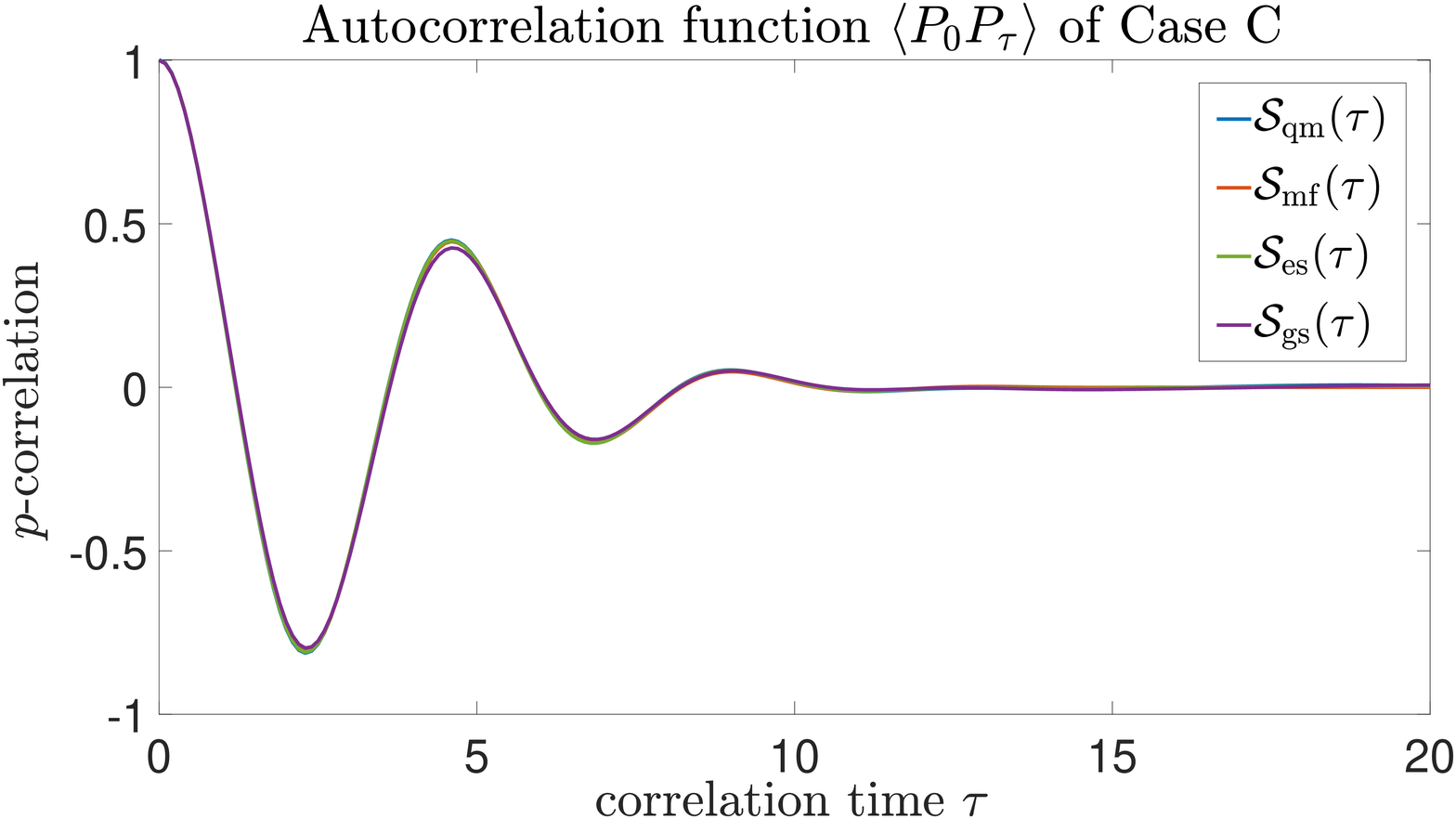}
         \caption{Case C: Auto-correlation curves}
         \label{fig:corr_case_A3}
     \end{subfigure}
     \hfill
     \begin{subfigure}[b]{0.49\textwidth}
         \centering
         \includegraphics[width=\textwidth]{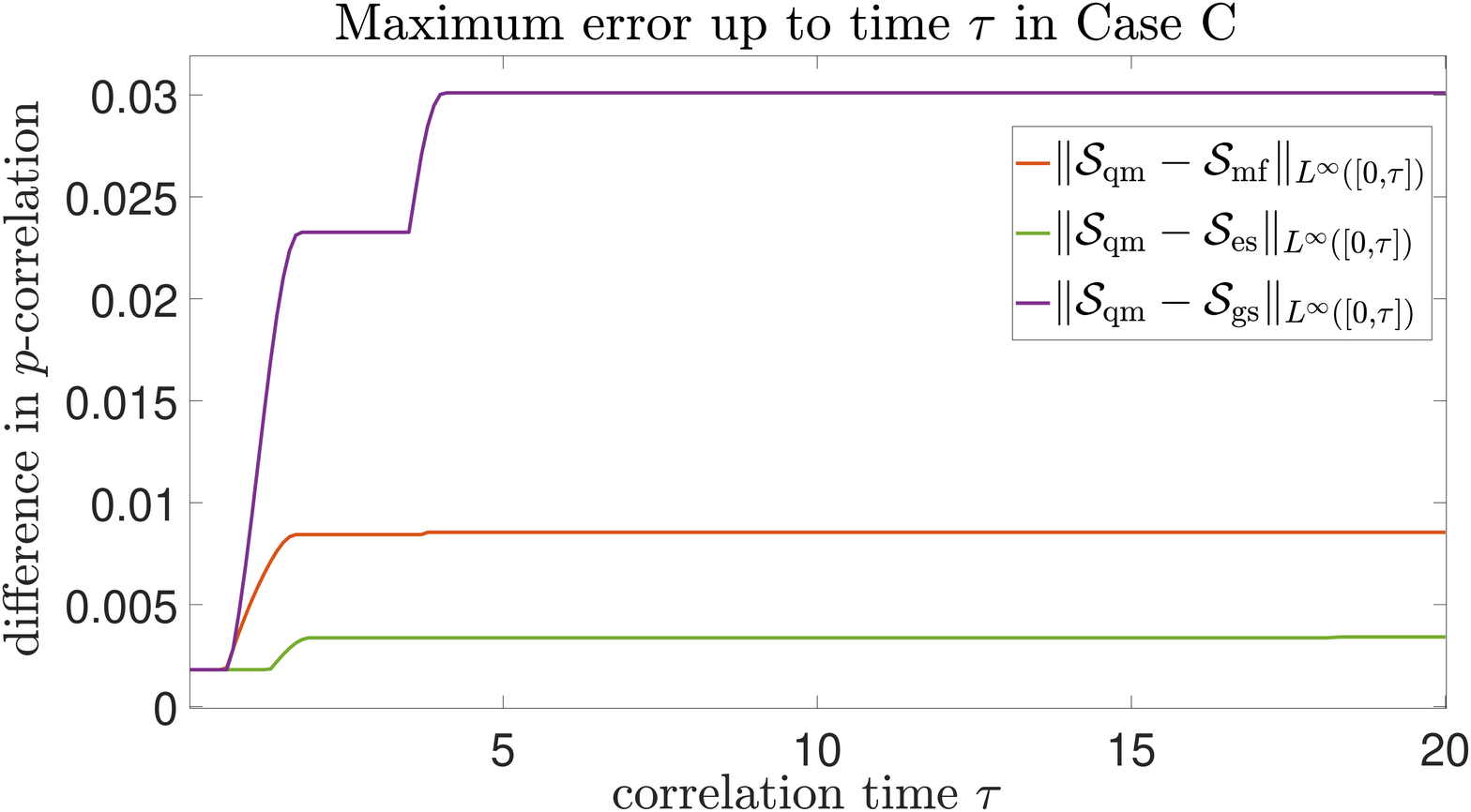}
         \caption{Case C: Maximum error curves up to time $\tau$}
         \label{fig:err_case_A3}
     \end{subfigure}
     \caption{Case C: (a) Auto-correlation function $\langle P_0P_\tau\rangle$ computed by quantum-mechanics formula, $\mathcal{S}_{\mathrm{qm}}$, with $M=100$, and by three molecular dynamics formulae. (b) The corresponding maximum errors up to time $\tau$ $\|\mathcal{S}_{\mathrm{qm}}-\mathcal{S}_{\mathrm{md}}\|_{L^\infty([0,\tau])}$.}
\end{figure}

The approximate $p$-auto-correlation function curves with their corresponding maximum errors up to time $\tau$ are plotted in the Figures \ref{fig:corr_case_A3} and \ref{fig:err_case_A3}, respectively. These two figures are quite similar to their corresponding plots in \hyperlink{case_B1}{Case B}, where the excited state approximation $\mathcal{S}_{\mathrm{es}}$ has the smallest error, and the mean-field approximation $\mathcal{S}_{\mathrm{mf}}$ achieves an improved accuracy compared to the ground state approximation $\mathcal{S}_{\mathrm{gs}}$. 

The similar approximation error of the three molecular dynamics in \hyperlink{case_B1}{Case B} and \hyperlink{case_C}{Case C} can be understood as a result of the relatively small difference between the two eigenvalues $\lambda_0$ and $\lambda_1$. For both cases the small parameter $c = 0.1$ leads to small $\epsilon_1$ and $\epsilon_2$ values, as summarized in  Table~\ref{table_case_A_to_E_new_written_caption}.



\noindent\hypertarget{case_D}{\textbf{Case D}}\textbf{: \textit{High temperature, large difference between eigenvalues with large gap}}, $\beta=0.28$, $c=1$, $\delta=1$, $\epsilon_1^2=2.01$, $\epsilon_2^2=0.42$, the probability for the excited state $q_1=0.30$.

\begin{figure}[!htb]
     \centering
     \begin{subfigure}[b]{0.49\textwidth}
         \centering
         \includegraphics[width=\textwidth]{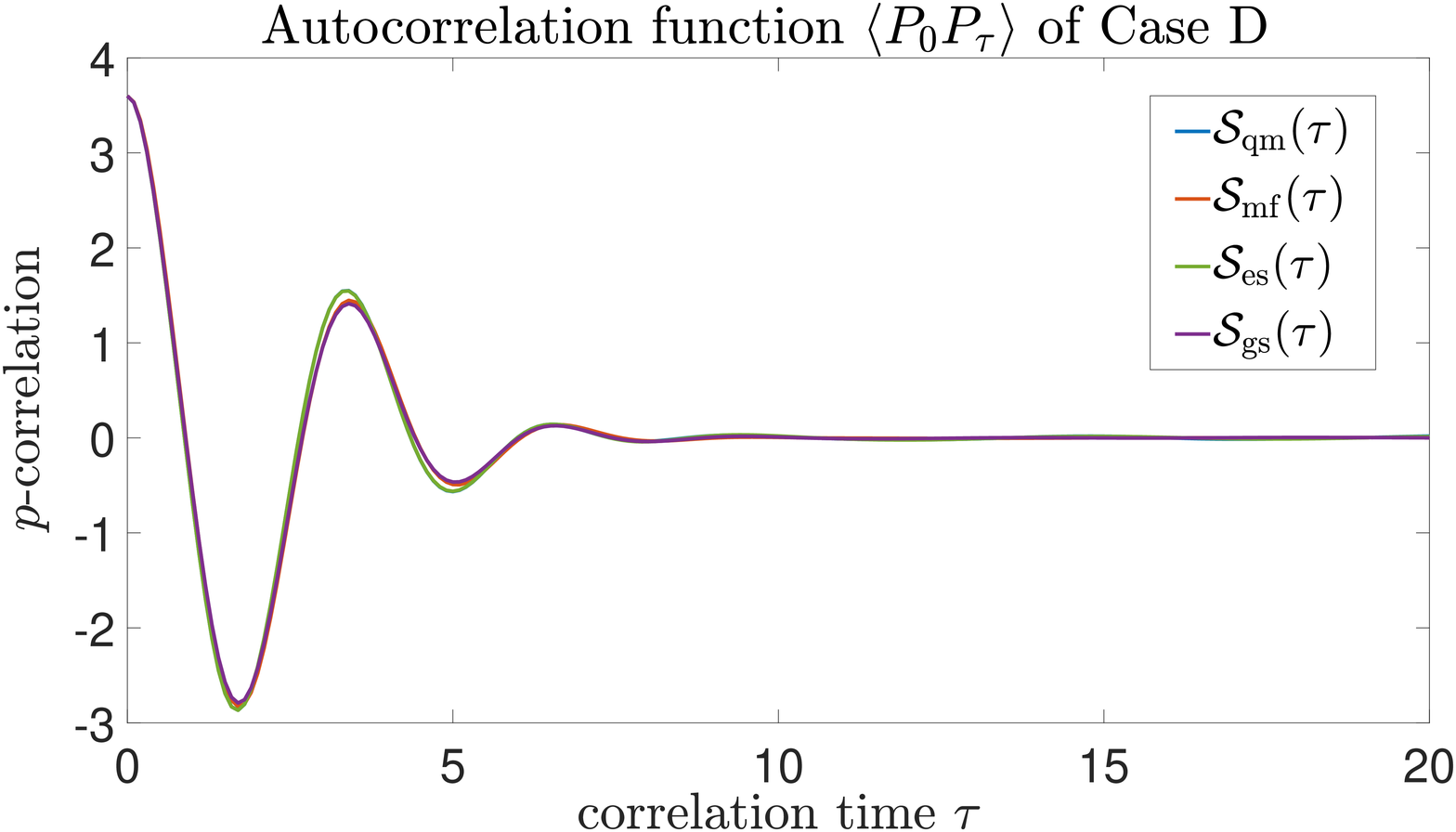}
         \caption{Case D: Auto-correlation curves}
         \label{fig:corr_case_A6}
     \end{subfigure}
     \hfill
     \begin{subfigure}[b]{0.49\textwidth}
         \centering
         \includegraphics[width=\textwidth]{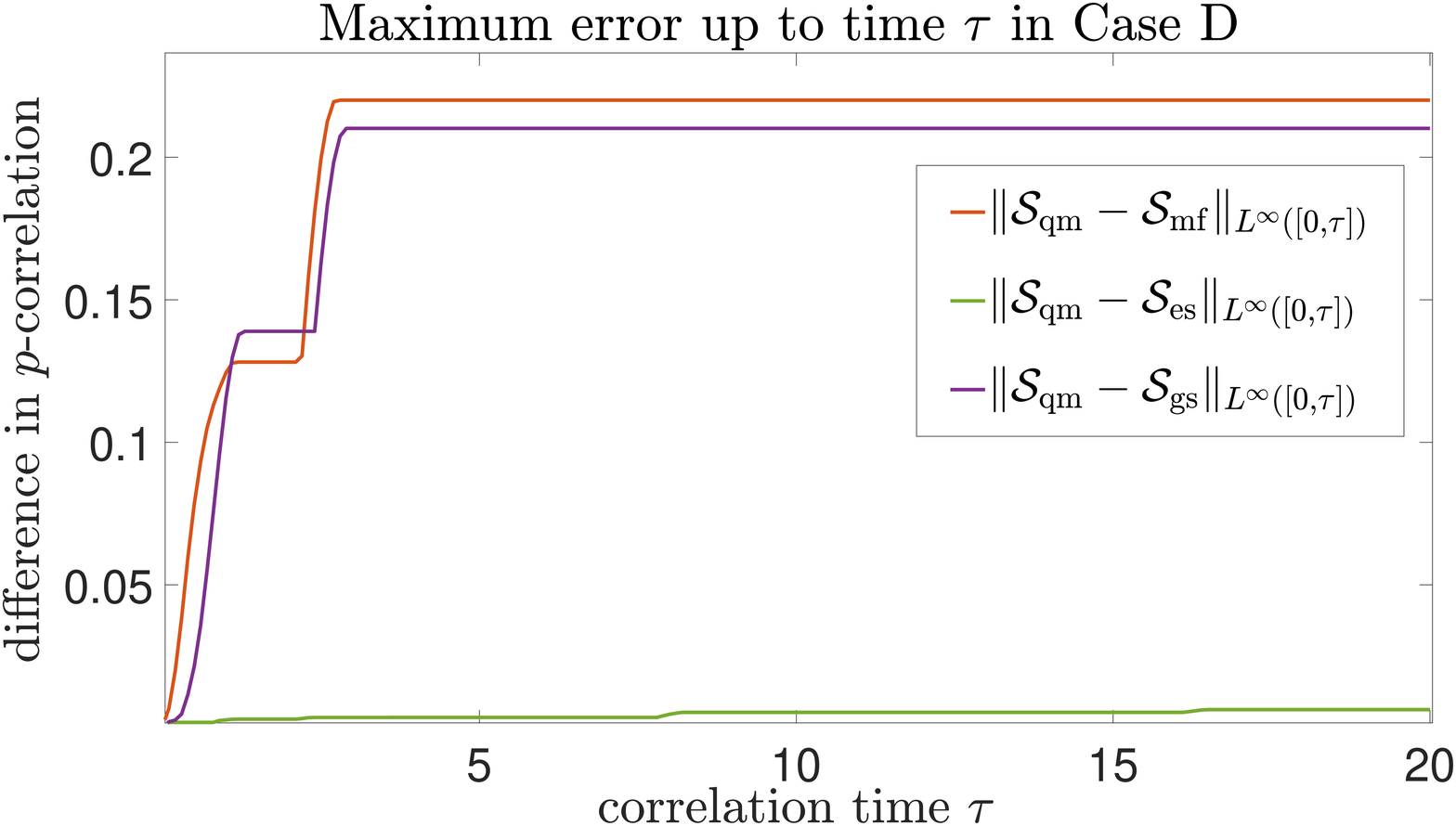}
         \caption{Case D: Maximum error curves up to time $\tau$}
         \label{fig:err_case_A6}
     \end{subfigure}
     \caption{Case D: (a) Auto-correlation function $\langle P_0P_\tau\rangle$ computed by quantum-mechanics formula, $\mathcal{S}_{\mathrm{qm}}$, with $M=100$, and by three molecular dynamics formulae. (b) The corresponding maximum errors up to time $\tau$ $\|\mathcal{S}_{\mathrm{qm}}-\mathcal{S}_{\mathrm{md}}\|_{L^\infty([0,\tau])}$.}
\end{figure}

For this case, we observe from Figure~\ref{fig:eigval_plot_case_D} that although the mean-field potential $\lambda_\ast(x)$ is still in between the ground state $\lambda_0(x)$ and the excited state $\lambda_1(x)$, the distance between $\lambda_\ast(x)$ and $\lambda_1(x)$ is much larger than that in the previous \hyperlink{case_B1}{Case B}. Also $\lambda_\ast(x)$ is much closer to  the ground state $\lambda_0(x)$ than to the excited state $\lambda_1(x)$. The parameter $\beta=0.28$ implies a relatively high temperature, with a considerable contribution from the excited state. Hence we cannot expect the mean-field molecular dynamics to be much better than ground state molecular dynamics. This is also verified by Figure~\ref{fig:err_case_A6} which shows that the error of mean-field molecular dynamics is of the same order as that of ground state molecular dynamics, while the excited state molecular dynamics remains accurate.

The mean-field and ground state molecular dynamics correlations include two approximations: replacing the matrix-valued potential $V(x)$ by a scalar potential $\lambda_\ast(x)$ or $\lambda_0(x)$ respectively, and replacing quantum dynamics with classical dynamics,
\begin{equation}
\begin{aligned}
\mathcal{S}_{\mathrm{qm}}(\tau)-\mathcal{S}_{\mathrm{mf}}(\tau)&=\big(\mathcal{S}_{\mathrm{qm}}(\tau)-\mathcal{S}_{\mathrm{qm},\lambda_\ast}(\tau)\big) + \big( \mathcal{S}_{\mathrm{qm},\lambda_\ast}(\tau)-\mathcal{S}_{\mathrm{mf}}(\tau)\big)\\
 \mathcal{S}_{\mathrm{qm}}(\tau)-\mathcal{S}_{\mathrm{gs}}(\tau)&=\big(\mathcal{S}_{\mathrm{qm}}(\tau)-\mathcal{S}_{\mathrm{qm},\lambda_0}(\tau)\big) + \big( \mathcal{S}_{\mathrm{qm},\lambda_0}(\tau)-\mathcal{S}_{\mathrm{gs}}(\tau)\big)
\end{aligned}\label{qm_scalar_potential_error_decompose}
\end{equation}
where $\mathcal{S}_{\mathrm{qm},\lambda_\ast}$ and $\mathcal{S}_{\mathrm{qm},\lambda_0}$ denote the approximation of auto-correlation function computed with quantum dynamics but using scalar-valued potentials $\lambda_\ast(x)$ and $\lambda_0(x)$, respectively. In the right hand side of \eqref{qm_scalar_potential_error_decompose}, the first terms $\big(\mathcal{S}_{\mathrm{qm}}(\tau)-\mathcal{S}_{\mathrm{qm},\lambda_\ast}(\tau)\big)$ and $\big(\mathcal{S}_{\mathrm{qm}}(\tau)-\mathcal{S}_{\mathrm{qm},\lambda_0}(\tau)\big)$ correspond to the potential approximations in quantum dynamics, while the second terms $\big( \mathcal{S}_{\mathrm{qm},\lambda_\ast}(\tau)-\mathcal{S}_{\mathrm{mf}}(\tau)\big)$ and $\big( \mathcal{S}_{\mathrm{qm},\lambda_0}(\tau)-\mathcal{S}_{\mathrm{gs}}(\tau)\big)$ are related to classical approximations of quantum dynamics using scalar potentials.

\begin{figure}[!htb]
     \centering
     \begin{subfigure}[b]{0.49\textwidth}
         \centering
         \includegraphics[width=\textwidth]{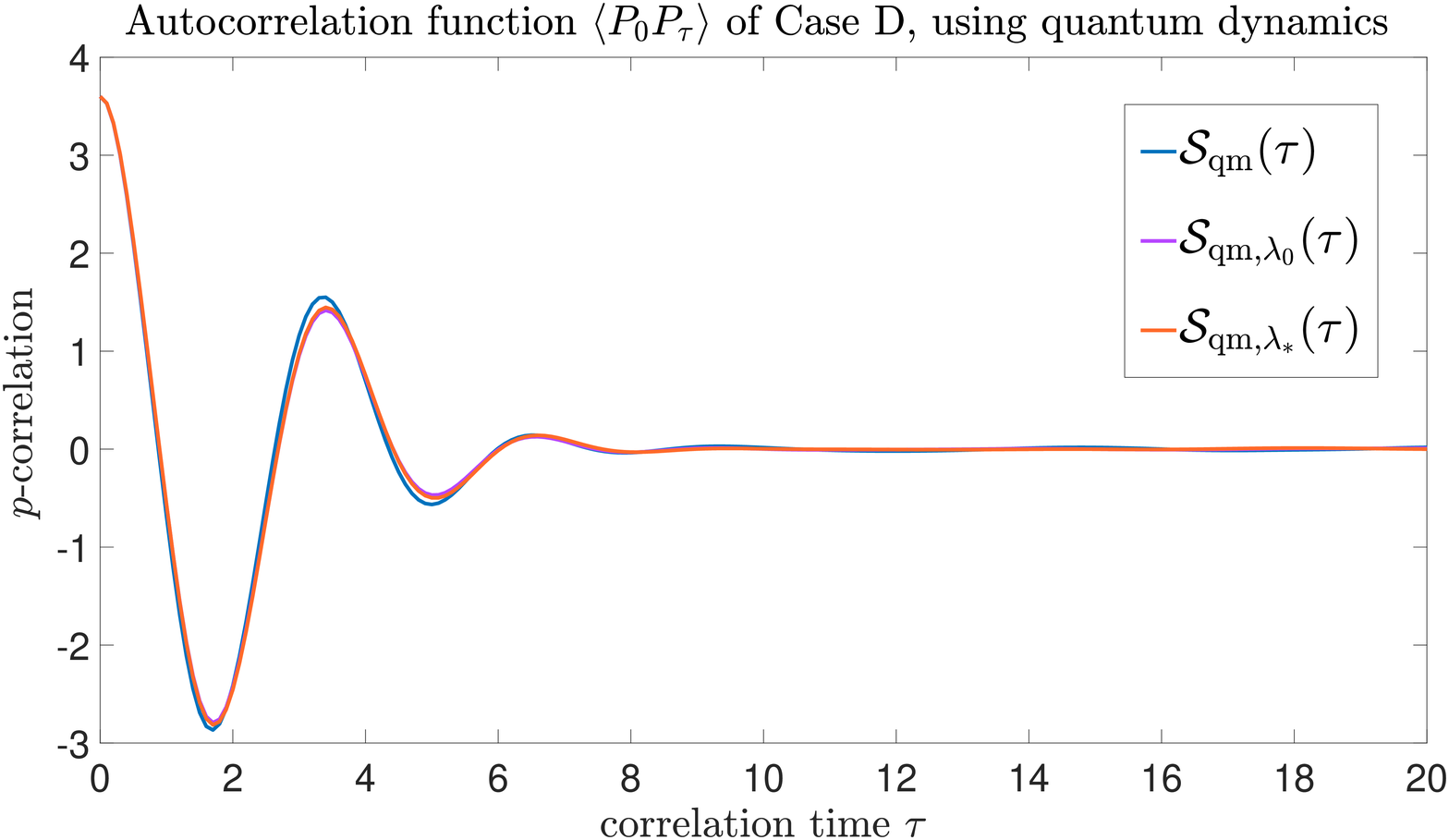}
         \caption{Case D: Auto-correlation curves}
         \label{fig:corr_qm_scalar_case_A6}
     \end{subfigure}
     \hfill
     \begin{subfigure}[b]{0.49\textwidth}
         \centering
         \includegraphics[width=\textwidth]{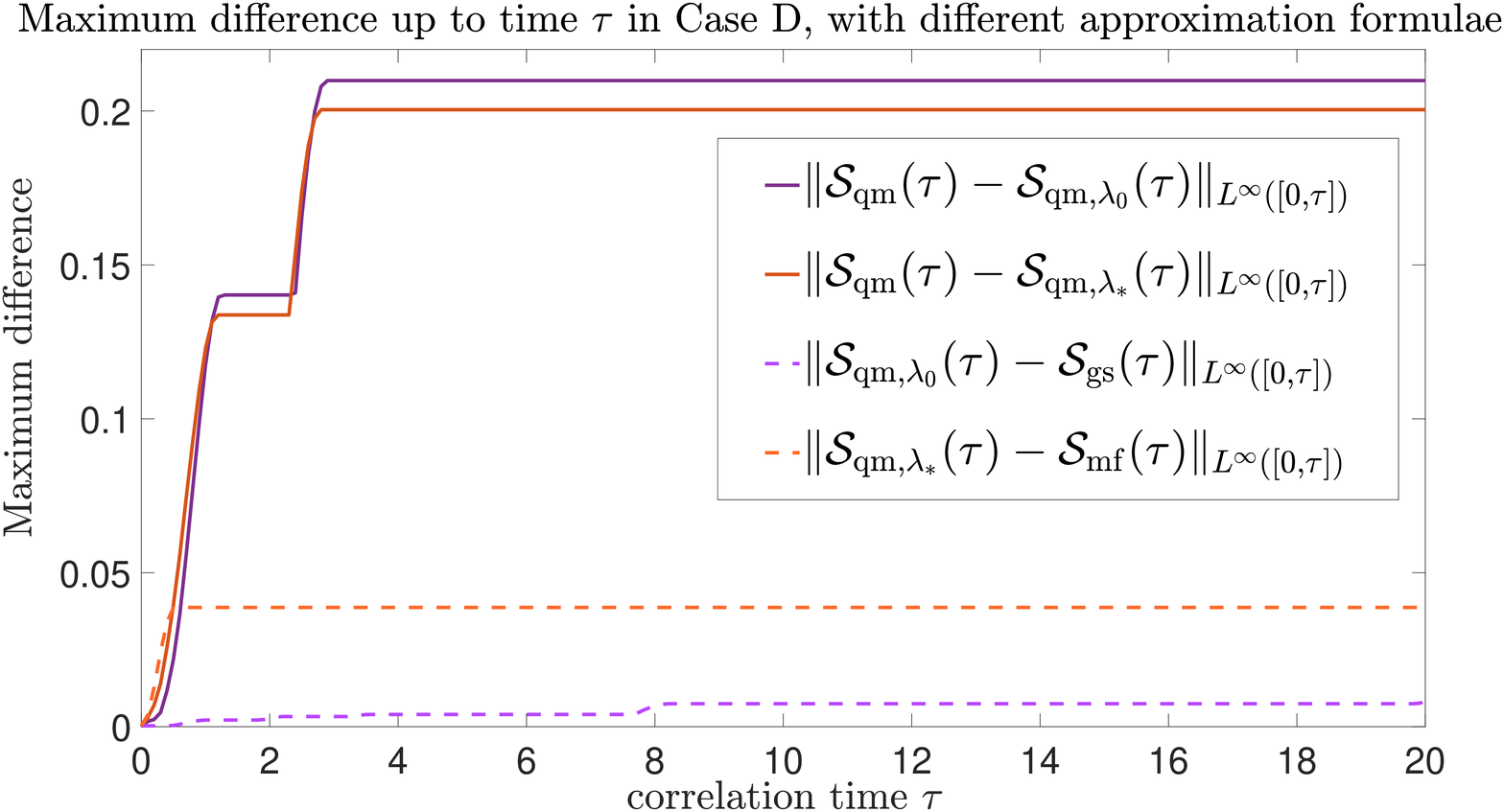}
         \caption{Case D: Maximum difference up to time $\tau$}
         \label{fig:err_qm_scalar_qm_case_A6}
     \end{subfigure}
     \caption{(a) Auto-correlation function $\langle P_0P_\tau\rangle$ curves $\mathcal{S}_{\mathrm{qm}}$, $\mathcal{S}_{\mathrm{qm},\lambda_0}$, and $\mathcal{S}_{\mathrm{qm},\lambda_\ast}$ computed by using matrix valued potential $V(x)$, and using scalar-valued potential $\lambda_0(x)$ or $\lambda_\ast(x)$ in the quantum mechanics formula \eqref{T_qm_def} in Case D. (b) Maximum error up to time $\tau$ in the $p$-auto-correlation curves computed with quantum mechanics formula using scalar-valued potential $\lambda_\ast(x)$ or $\lambda_0(x)$, by comparing them with the correlation computed from quantum mechanics formula using matrix-valued potential $V(x)$, and with their corresponding molecular dynamics approximations.
 }
\end{figure}

To investigate these two error contributions we compute the correlation function $\mathcal{S}_{\mathrm{qm},\lambda_\ast}(\tau)$ and $\mathcal{S}_{\mathrm{qm},\lambda_0}(\tau)$ for \hyperlink{case_D}{Case D}, using the scalar-valued potential $\lambda_\ast(x)$ or $\lambda_0(x)$ to replace the potential matrix $V(x)$ in the quantum dynamics. 
The corresponding auto-correlation curves and their maximum error up to time $\tau$ are shown in 
Figures~\ref{fig:corr_qm_scalar_case_A6} and \ref{fig:err_qm_scalar_qm_case_A6}. 


From Figure~\ref{fig:err_qm_scalar_qm_case_A6}, we clearly see that the errors $\|\mathcal{S}_{\mathrm{qm}}-\mathcal{S}_{\mathrm{qm},\lambda_\ast}\|_{L^\infty([0,\tau])}$ and $\|\mathcal{S}_{\mathrm{qm}}-\mathcal{S}_{\mathrm{qm},\lambda_0}\|_{L^\infty([0,\tau])}$ caused by substituting the potential matrix $V(x)$ with the 
scalar-valued potentials $\lambda_\ast(x)$ or $\lambda_0(x)$ is of the same order as the total errors $\|\mathcal{S}_{\mathrm{qm}}-\mathcal{S}_{\mathrm{mf}}\|_{L^\infty([0,\tau])}$ and $\|\mathcal{S}_{\mathrm{qm}}-\mathcal{S}_{\mathrm{gs}}\|_{L^\infty([0,\tau])}$ in Figure~\ref{fig:err_case_A6}. 
Hence we conclude that the main source of error in this case is the simplification by replacing the potential matrix $V(x)$ with a scalar-valued potential, and not the approximation of scalar potential quantum mechanics correlation by classical molecular dynamics.

We also vary the mass ratio $M$ between the heavy particle 
and the light particle, in order to study the corresponding behaviour of 
the approximation errors $\|\mathcal{S}_{\mathrm{qm}}(\tau)-\mathcal{S}_{\mathrm{qm},\lambda_\ast}(\tau)\|_{L^\infty([0,\tau])}$ and
$\|\mathcal{S}_{\mathrm{qm}}(\tau)-\mathcal{S}_{\mathrm{qm},\lambda_0}(\tau)\|_{L^\infty([0,\tau])}$ up to time $\tau=20$ for the mean-field molecular dynamics and ground state molecular dynamics in \hyperlink{case_D}{Case D}. As can be seen from the second and fourth column of
Table~\ref{table_change_M_case_D}, the main error caused by substitution of the potential matrix $V$ with the scalar valued potential $\lambda_\ast$ or $\lambda_0$ varies slightly as the $M$ value changes. 


\begin{table}[h]
\begin{center}
\begin{tabular}{ |c|c|c|c|c|c| }
\hline
 $M$ & $\|\mathcal{S}_{\mathrm{qm}}-\mathcal{S}_{\mathrm{qm},\lambda_\ast}\|$ & $\| \mathcal{S}_{\mathrm{qm,\lambda_\ast}} - \mathcal{S}_{\mathrm{mf}}\|$
 & $\|\mathcal{S}_{\mathrm{qm}}-\mathcal{S}_{\mathrm{qm},\lambda_0}\|$ & $\| \mathcal{S}_{\mathrm{qm,\lambda_0}} - \mathcal{S}_{\mathrm{gs}}\|$ & $\|\mathcal{S}_{\mathrm{qm}}-\mathcal{S}_{\mathrm{es}}\|$\\
 \hline
 100 & 0.2004 & 0.0387 & 0.2099 & 0.0081 & 0.0062 \\ 
 \hline
 50 & 0.1985 & 0.0384 & 0.2080 & 0.0201 & 0.0220 \\ 
 \hline
 20 & 0.1944 & 0.0375 & 0.2033 & 0.0377 & 0.0330 \\ 
 \hline
\end{tabular} 
\end{center}
\caption{Case D: Dependence of the error on the mass ratio $M$ at time $\tau=20$.}
\label{table_change_M_case_D}
\end{table}

\noindent\hypertarget{case_E}{\textbf{Case E}}\textbf{: \textit{High temperature, large difference between eigenvalues with avoided crossing}}, with $\beta=1$, $c=1$, $\delta=0.1$, $\epsilon_1^2=0.29$, $\epsilon_2^2=0.50$, the probability for the excited state $q_1=0.16$.

\begin{figure}[!htb]
     \centering
     \begin{subfigure}[b]{0.50\textwidth}
         \centering
         \includegraphics[width=\textwidth]{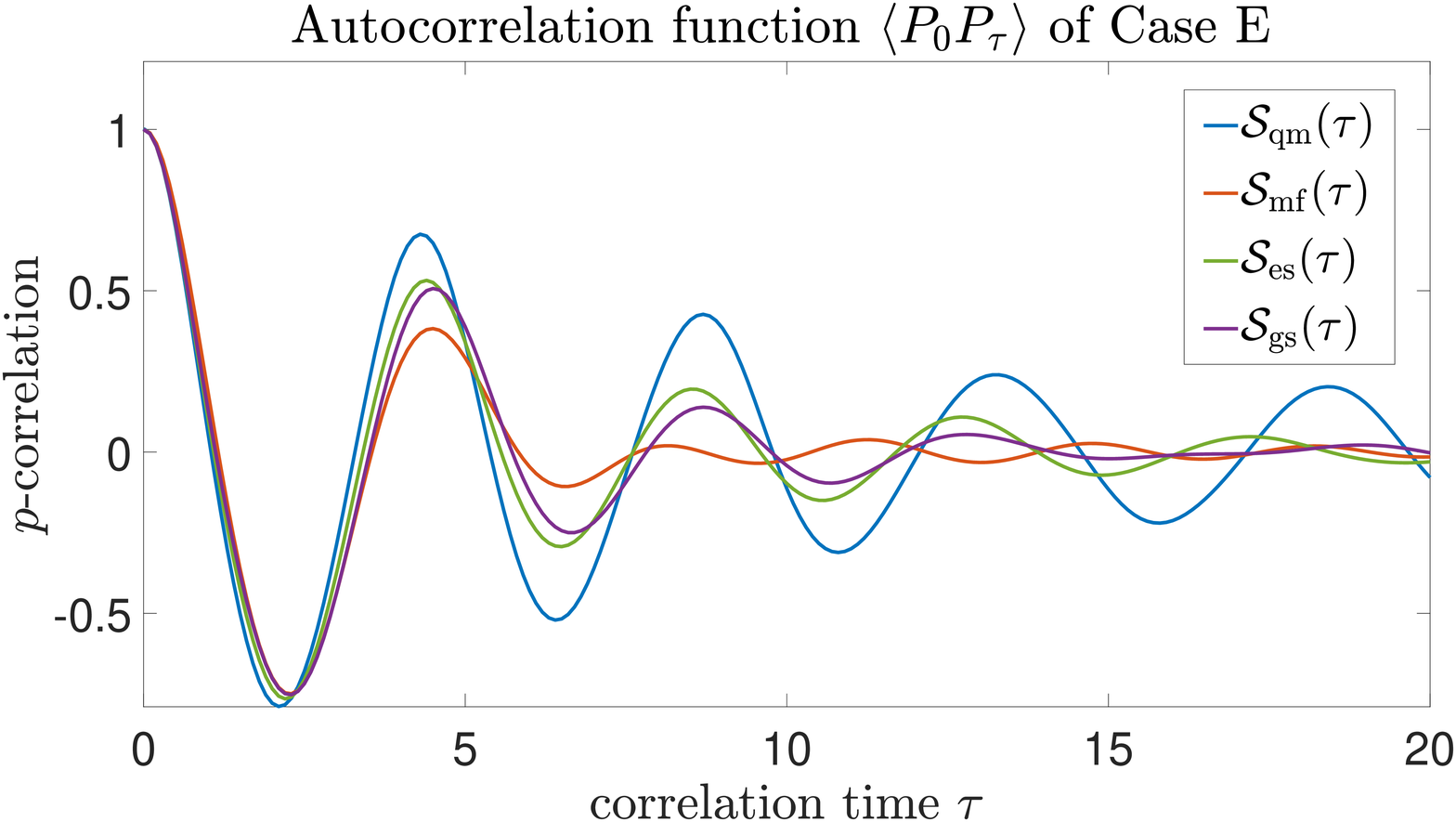}
         \caption{Case E: Auto-correlation curves}
         \label{fig:corr_case_A5}
     \end{subfigure}
     \hfill
     \begin{subfigure}[b]{0.49\textwidth}
         \centering
         \includegraphics[width=\textwidth]{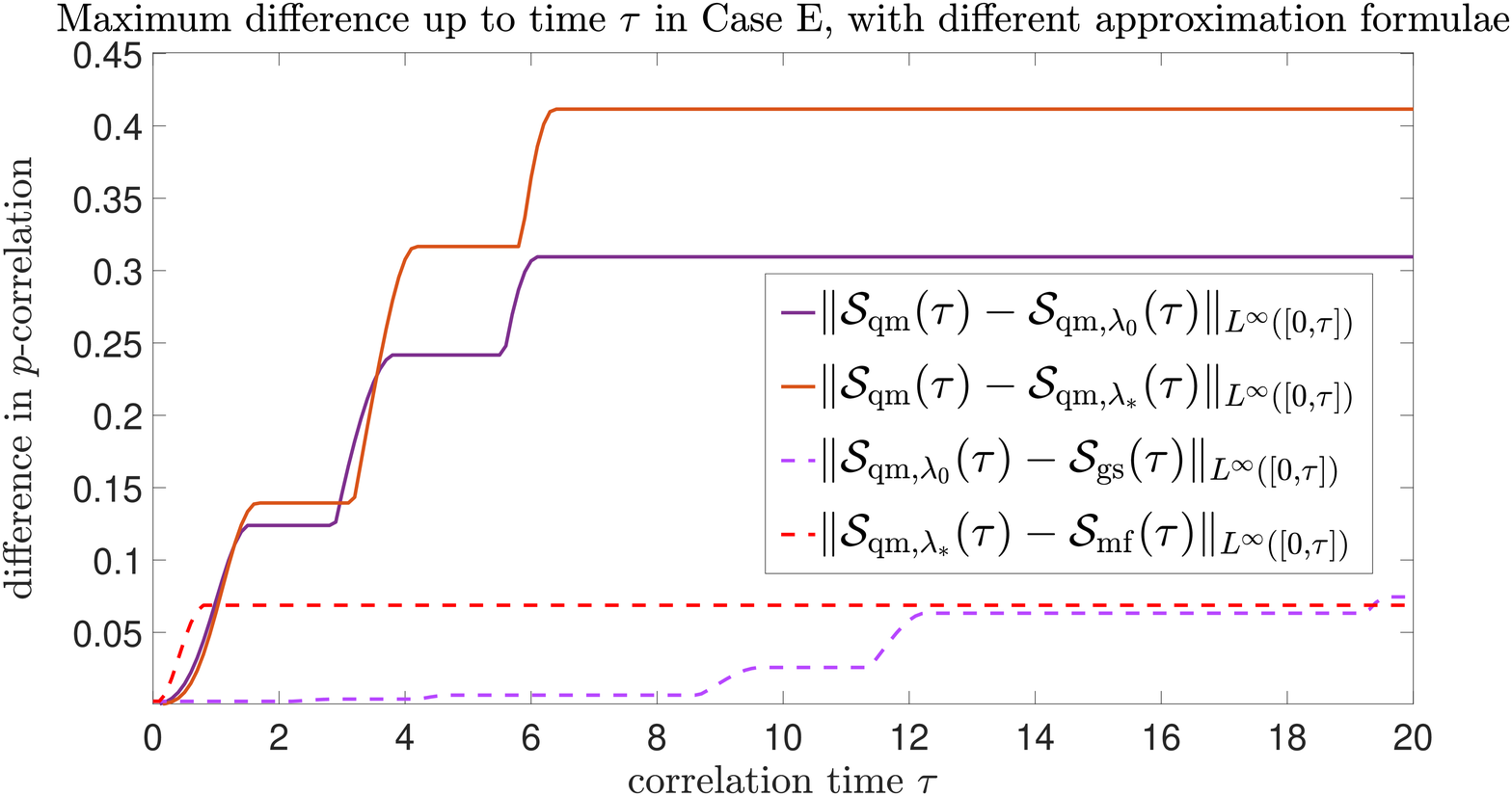}
         \caption{Case E: Maximum error curves up to time $\tau$}
         \label{fig:err_case_A5}
     \end{subfigure}
     \caption{(a) Auto-correlation function $\langle P_0P_\tau\rangle$ curves $\mathcal{S}_{\mathrm{qm}}$, $\mathcal{S}_{\mathrm{qm},\lambda_0}$, and $\mathcal{S}_{\mathrm{qm},\lambda_\ast}$ computed by using matrix valued potential $V(x)$ and using scalar-valued potential $\lambda_0(x)$ or $\lambda_\ast(x)$ in the quantum mechanics formula \eqref{T_qm_def} in Case E. (b) Maximum error up to time $\tau$ in the $p$-auto-correlation curves computed with quantum mechanics formula using scalar-valued potential $\lambda_\ast(x)$ or $\lambda_0(x)$, by comparing them with the correlation computed from quantum mechanics formula using matrix-valued potential $V(x)$, and with their corresponding molecular dynamics approximations.
 }
\end{figure}

This case has the same parameters as in Section~\ref{subsection_5_1}. In Figure~\ref{fig:eigval_plot_A51} we observe a pattern of the two eigenvalues $\lambda_0(x)$ and $\lambda_1(x)$ related to avoided crossing of potential surfaces. Our numerical results suggest that for this case all three molecular dynamics are only accurate for short time, as can be seen in Figure~\ref{fig:corr_case_A5}. 
Compared to \hyperlink{case_C}{Case C} and \hyperlink{case_D}{Case D}, where the excited state dynamics is accurate, the diminished eigenvalue regularity at the avoided crossing may explain the loss of accuracy in \hypertarget{case_E}{Case E}.

\begin{figure}[!htb]
    \begin{center}
  \includegraphics[height=4.5cm]{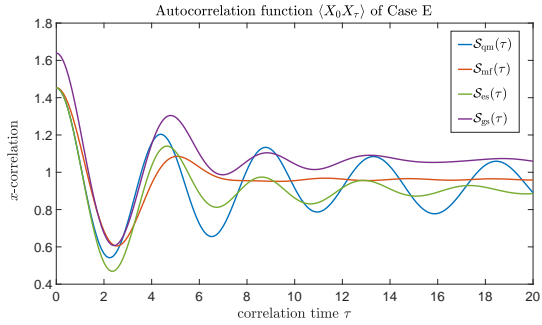}
  \caption{Case E: Auto-correlation function $\langle X_0X_\tau\rangle$ computed by quantum-mechanics formula with $M=100$ and by three molecular dynamics formulae.}\label{fig:corr_case_A5_xcorr}
  \end{center}
\end{figure}

Apart from the momentum auto-correlation function, we also computed in \hyperlink{case_E}{Case E} the correlation function between position observables $\widehat{x}_0$ and $\widehat{x}_\tau$, as plotted in Figure~\ref{fig:corr_case_A5_xcorr}. We observe that  for  short time range (e.g. $0\leq \tau \leq 0.1$), the error of the ground state molecular dynamics is  larger than the error of the mean-field molecular dynamics, which is consistent with the result for equilibrium observables in Section~\ref{subsection_5_1}, in which the ground state molecular dynamics has  larger error in approximating the density function $\mu_{\mathrm{qm}}(x)$ than the mean-field formula. Therefore the mean-field molecular dynamics can improve short-time approximation of position auto-correlation function compared to the ground state molecular dynamics.

We also changed the mass ratio $M$ in this case for the momentum auto-correlation function, from $M=100$ to $M=50$ and to smaller value $M=20$. When $M$ becomes smaller, we can expect the error of molecular dynamics approximation becomes larger, since the error includes the $\mathcal{O}(M^{-1})$ term. For \hyperlink{case_E}{Case E}, since we are only interested in the short time approximation, the time-dependent error term is not much larger than the $\mathcal{O}(M^{-1})$ term. Hence the effect of varying the mass ratio $M$ will be considerable. The dependence of the $L^\infty$-error in momentum auto-correlation approximation on the mass ratio $M$ is summarized in Table~\ref{table_change_M}, from which we observe an improved accuracy in all the three molecular dynamics approximations with an increased $M$ value.

\begin{table}[h]
\begin{center}
\begin{tabular}{ |c|c|c|c|c| }
\hline
$\tau$ & $M$ & $\|\mathcal{S}_{\mathrm{qm}}-\mathcal{S}_{\mathrm{mf}}\|_{L^\infty([0,\tau])}$ & $\|\mathcal{S}_{\mathrm{qm}}-\mathcal{S}_{\mathrm{es}}\|_{L^\infty([0,\tau])}$ & $\|\mathcal{S}_{\mathrm{qm}}-\mathcal{S}_{\mathrm{gs}}\|_{L^\infty([0,\tau])}$\\
 \hline
 1 & 20 & 0.1375 & 0.0816 & 0.0564 \\ 
 \hline
 1 & 50 & 0.1324 & 0.0766 & 0.0513 \\ 
 \hline
 1 & 100 & 0.1270 & 0.0712 & 0.0460 \\  
 \hline
 2 & 20 & 0.1657 & 0.1515 & 0.1043 \\ 
 \hline
 2 & 50 & 0.1533 & 0.1356 & 0.0884 \\ 
 \hline
 2 & 100 & 0.1425 & 0.1226 & 0.0759 \\  
 \hline
\end{tabular} 
\end{center}
\caption{Case E: Dependence of the error on the mass ratio $M$ at different correlation times $\tau$.}
\label{table_change_M}
\end{table}

\subsection{Conclusion from numerical comparisons}
From the study with equilibrium observables in Section~\ref{subsection_5_1}, we see that by considering the contributions of excited states the classical mean-field approximation of quantum mechanics density at equilibrium achieves a substantial improvement from the approximation which uses only the information from the ground state. The error of the mean-field approximation will decrease as the mass ratio $M$ increases, following the $\mathcal{O}(M^{-1})$ relation.

For the time-dependent observables, specifically by studying the momentum auto-correlation function, we know from \hyperlink{case_A}{Case A} that for a low temperature setting with a large eigenvalue gap, where the probability for an excited state is small, all three molecular dynamics with the mean-field approximation, excited state approximation, or ground state approximation work similarly well.

From \hyperlink{case_B1}{Case B} and \hyperlink{case_D}{Case D}, we observe that the error of the mean-field approximation decreases as the difference between two eigenvalues diminishes (i.e. parameter $c$ becomes small). Furthermore, for \hyperlink{case_B1}{Case B} with a small difference between two eigenvalues, the mean-field approximation improves the accuracy of molecular dynamics compared to using the ground state only.

With a small eigenvalue difference, even including the avoided crossing in \hyperlink{case_C}{Case C} the result is similar to \hyperlink{case_B1}{Case B}, that is the mean-field approximation is still more accurate than the ground state approximation.

From \hyperlink{case_D}{Case D} we know that when the system temperature is high and the difference between two eigenvalues is not small, the excited state approximation outperforms both the mean-field and the ground state molecular dynamics. 

From \hyperlink{case_E}{Case E} we see that when the difference between the eigenvalues are sufficiently large and when the potential matrix includes avoided crossings, all three molecular dynamics approximations are accurate for a short time range only.

The small error terms $\epsilon_1^2$ and $\epsilon_2^2$, in Table \ref{table_case_A_to_E_new_written_caption}, for the Cases A, B and C
is consistent with the actual error being small for the mean-field approximation, while in Cases D and E
where the mean-field approximation error is large these error terms are in fact large. Therefore
the experiments indicate that the error estimate could be useful to estimate the mean-field approximation error
also for realistic problems when the quantum observable is not computable.

Figures \ref{fig:err_qm_scalar_qm_case_A6} and \ref{fig:err_case_A5} together with Table~\ref{table_change_M_case_D} 
show that in  \hyperlink{case_D}{Case D} and \hyperlink{case_E}{Case E}, where mean-field and ground state approximations are not accurate, the error of the mean-field and ground state molecular dynamics are dominated by the matrix valued potential replaced by a scalar potential on the quantum level, since the classical approximation error of the quantum dynamics for the corresponding scalar potential is clearly smaller.

\section*{Acknowledgment}
This research was supported by
Swedish Research Council grant 2019-03725 and KAUST grant OSR-2019-CRG8-4033.3. 
The work of P.P. was supported in part by the U.S. Army Research Office Grant W911NF-19-1-0243.

\bibliographystyle{abbrv}

\begin{thebibliography}{99}
\bibitem{amour} Amour L., Jager L. and  Nourrigat J.,
\textit{The Weyl symbol of Schr\"odinger semigroups},
Ann. Henri Poincare \textbf{16} (2015), 1479–1488.

\bibitem{robert}
Bouzounia A. and Robert D., \emph{Uniform semiclassical estimates for the
  propagation of quantum observables}, Duke Math. J. \textbf{111} (2002),
  223--252.

\bibitem{cances_nonlinear}
 {Canc{\`{e}}s E., Le Bris C. and Lions P-L.},
\textit {Molecular simulation and related topics: some open mathematical problems},
 Nonlinearity \textbf{21} (2008), T165-T176.


\bibitem{dornheim} Dornheim T., Groth S.,  Filinov A.V. and Bonitz M., \textit{Path integral Monte Carlo simulation
of degenerate electrons: Permutation-cycle
properties}, J. Chem. Phys. \textbf{151} 014108 (2019).

\bibitem{evans}
Evans L. C., \emph{Partial differential equation}, American Mathematical
  Society, Providence, RI, 1998.
  
  \bibitem{feynman}  Feynman R.F., {\em Statistical Mechanics: A Set of Lectures} (Westview Press 1998).
  
  \bibitem{Figalli-Paul}Figalli A., Ligabò M. and Paul T., {\em  Semiclassical Limit for Mixed States with Singular and Rough Potentials,} Indiana University Mathematics Journal, \textbf{61} (2012) 193--222.
 


  \bibitem{ford_kac2} Ford G. W. and Kac M., {\em On the quantum Langevin equation.}  J. Statist. Phys.  \textbf{46}  (1987), 
803--810.


\bibitem{Golse-Paul}François G. and Paul T.,
{\em Semiclassical evolution with low regularity,}
Journal de Mathématiques Pures et Appliquées,
\textbf{151} (2021) 257--311.
  
\bibitem{gronwall} Gronwall, T. H., \textit{ Note on the derivatives with respect to a parameter of the solutions of a system of differential equations}, Ann. of Math. \textbf{20}  (1919) 292–296.

  \bibitem{ring-poly}Habershon S., Manolopoulos D.E., Markland T.E. and Miller T.F. 3rd., 
\textit{Ring-polymer molecular dynamics: quantum effects in chemical dynamics from classical trajectories in an extended phase space.}
Annu Rev Phys Chem., \textbf{64} (2013) 387-413.

\bibitem{Hansen}Hansen J.-P. and McDonald I.R., {\em Theory of Simple Liquids}, Academic Press (1985).

\bibitem{HS}Hoel H. and Szepessy A., {\em Classical Langevin dynamics derived from quantum mechanics},
Discrete \& Continuous Dynamical Systems - B, \textbf{25} (2020) 4001-4038.
  
  \bibitem{KPSS} Kammonen A., Plech\'a\v{c} P., Sandberg M. and Szepessy A., {\em Canonical quantum observables for molecular systems approximated by ab initio molecular dynamics},  Ann. Henri Poincare \textbf{19} (2018), 2727-2781.
  
  \bibitem{karatzas} Karatzas I. and Shreve S.E., {\em Brownian Motion and Stochastic Calculus} Springer Verlag (1998).
  
  \bibitem{md_book} Leimkuhler B. and Matthews C., {\em Molecular Dynamics: With Deterministic and Stochastic Numerical Methods,} Springer, Berlin (2015).

\bibitem{Lions-Paul}
Lions P.-L. and Paul T., {\em Sur les mesures de Wigner}, Revista Matemática Iberoamericana, \textbf{9.3} (1993) 553-618.

  \bibitem{martinez_book}
Martinez A., \emph{An introduction to semiclassical and microlocal analysis}, Springer Verlag, 2002.

\bibitem{marx_hutter}Marx D. and Hutter J., \textit{ Ab Initio Molecular Dynamics: Basic theory and advanced methods},
Cambridge University Press (2009).

\bibitem{morandi}Morandi G., Napoli F. and Ercolessi E., {\em Statistical Mechanics: An Intermediate Course}, World Scientific Publishing (2001).

\bibitem{papanicolaou} Papanicolaou G. and Hersh R., \textit{Non-Commuting Random Evolutions, and an Operator-Valued Feynman-Kac Formula},
Comm. Pure Appl. Math., \textbf{XXV} (1972) 337-367.

\bibitem{tuckerman_compare} Pérez A., Tuckerman M.E., M\"user M.H.,  \textit{A comparative study of the centroid and ring-polymer molecular dynamics methods for approximating quantum time correlation functions from path integrals,} J Chem Phys. 2009 May \textbf{14};130(18):184105. 

\bibitem{Teufel} Stiepan H.-M. and Teufel S., \textit{ Semiclassical approximations for Hamiltonians with operator-valued symbols}, Comm. Math. Phys. \textbf{320} (3) (2013) 821-849.

\bibitem{teller} 
Teller E., \emph{The crossing of potential surfaces}, J. Phys. Chem. \textbf{41} (1937), 109--116.

\bibitem{wigner} Wigner E., \textit{ On the Quantum Correction For Thermodynamic Equilibrium},
Phys. Rev. \textbf{40}  (1932) 749-759.

\bibitem{zworski}Zworski  M., \textit{ Semiclassical Analysis},  Providence, RI, American Mathematical Society (2012).



\end{thebibliography}


\appendix
\section{Numerical Implementations}
\subsection{Finite difference approximation of the equilibrium density}\label{subsection_5_3_1}
To approximate the quantum mechanics density formula \eqref{mu_qm}, we use a fourth-order finite difference approximation of the Laplacian in the Hamiltonian operator $\widehat{H}$ (5.1), with the formula
\[
f''(x)=\frac{-f(x-2h)+16f(x-h)-30f(x)+16f(x+h)-f(x+2h)}{12h^2}+\mathcal{O}(h^4).
\]
This discretization is performed on the computational domain $\Omega=[-6,6]$ with a uniform mesh $x_k=-6+k\Delta x$, for $k=0,1,\cdots,K$, with $\Delta x=\frac{12}{K}$. The choice of this computational domain is obtained by checking that the quantum mechanics density $\mu_{\mathrm{qm}}(x)$ approximately vanishes on the boundary of this domain, so that the homogeneous Dirichlet boundary condition can be assumed. By applying this discretization, we approximate the eigenvalue problem
\[ \widehat{H}\Phi_n=E_n\Phi_n, \]
with the following algebraic eigenvalue problem
\begin{equation}\label{fdm_eigval_problem}
    H_d\phi_n=e_n\phi_n,
\end{equation}
where the $2(K+1)\times 2(K+1)$ matrix $H_d$ is given by
\begin{equation}
\begin{split}
&H_d:=\frac{1}{2M\cdot 12(\Delta x)^2}\times\\
&\left[\begin{smallmatrix}
h_{11,0}+30 & h_{12,0} & -16 & 0 & 1 &   &   &   &   &   &  \\
h_{21,0} & h_{22,0}+30 & 0 & -16 & 0 & 1 &   &   &   &   &  \\
-16 & 0 & h_{11,1}+30 & h_{12,1} & -16 & 0 & 1 &   &   &   &  \\
0 & -16 & h_{21,1} & h_{22,1}+30 & 0 & -16 & 0 & 1 &   &   &  \\
1 & 0 & -16 & 0 & h_{11,2}+30 & h_{12,2} & -16 & 0 & 1 &   &  \\
  & 1 & 0 & -16 & h_{21,2} & h_{22,2}+30 & 0 & -16 & 0 & 1 &  \\
  &   & \ddots & \ddots & \ddots & \ddots & \ddots & \ddots & \ddots & \ddots & \ddots \\
  &   &   & 1 & 0 & -16 & 0 & h_{11,K-1}+30 & h_{12,K-1} & -16 & 0 \\
  &   &   &   & 1 & 0 & -16 & h_{21,K-1} & h_{22,K-1}+30 & 0 & -16 \\
  &   &   &   &   & 1 & 0 & -16 & 0 & h_{11,K}+30 & h_{12,K} \\
  &   &   &   &   &   & 1 & 0 & -16 & h_{21,K} & h_{22,K}+30 \\
\end{smallmatrix}\right]
\end{split}\label{Hd-matrix}
\end{equation}
and the eigenfunctions $\Phi_n$ are approximated with the $2(K+1)$-length vector $\phi_n$, as
\[ \phi_n=[\phi_{n,0,1},\phi_{n,0,2},\phi_{n,1,1},\phi_{n,1,2},\cdots,\phi_{n,K,1},\phi_{n,K,2}]^T. \]
Here, the entry terms $h_{ij,k}$ of the $H_d$ matrix are given by
\[ h_{ij,k}=2M(12\Delta x^2)V_{ij}(x_k),\quad\textup{for }i,j=1,2,\textup{ and }k=0,1,\cdots,K. \]
In practice, based on the finite difference scheme \eqref{fdm_eigval_problem}, we approximate the quantum mechanics density $\mu_{\mathrm{qm}}$ in \eqref{mu_qm} by
\begin{equation}\label{(5.4)}
\frac{\sum_n\Big( |\phi_{n,k,1}|^2+|\phi_{n,k,2}|^2 \Big)e^{-\beta e_n}}{\sum_k\sum_n\Big( |\phi_{n,k,1}|^2+|\phi_{n,k,2}|^2 \Big)e^{-\beta e_n}\Delta x},\ \textup{ for }\ k=0,1,\cdots, K,
\end{equation}
and the eigenvalues $e_n$ and eigenvectors $\phi_n$ here are obtained by using the \texttt{Matlab} function \texttt{eig}.

\subsection{Numerical solution of the mean-field Hamiltonian system}\label{subsection_5_3_2}

In the mean-field trace formula \eqref{mean-field_md_formula} for time correlation function $\mathcal{S}_{\mathrm{mf}}(\tau)$, we need to solve the Hamiltonian system
\[ 
\begin{aligned}
\dot{x}_t&=\nabla_p h(x_t,p_t),\\
\dot{p}_t&=-\nabla_x h(x_t,p_t),\\
\end{aligned}
\]
to obtain the state variable $x_\tau$ at time $t=\tau$.
Specifically, given the initial state $(x_0,p_0)$ at time $t_0=0$, we apply the velocity Verlet method, where for each discrete time point $t_n:=t_0+n\Delta t$, the dynamics of state variables $(x,p)$ in \eqref{h} is approximated by
\[ 
\begin{aligned}
p_{n+\frac{1}{2}}&=p_n+\frac{\Delta t}{2}\cdot\big( -\nabla_x h(x_n,p_n) \big),\\
x_{n+1} &= x_n + \Delta t\cdot p_{n+\frac{1}{2}},\\
p_{n+1}&=p_{n+\frac{1}{2}} + \frac{\Delta t}{2}\cdot\big( -\nabla_x h(x_{n+1},p_{n+\frac{1}{2}}) \big).
\end{aligned}
\]
The integrals in molecular dynamics formulas \eqref{mean-field_md_formula}, \eqref{excited_state_formula}, and \eqref{(ground-state_formula)} are computed with a fourth-order composite Simpson's method, with a discretized mesh $x_l=x_0+l\Delta x$, $p_l=p_0+l\Delta p$, for $l=0,1,\cdots L$ on the phase space $(x,p)$ and the computational domain $\Omega=[x_0,x_L]\times[p_0,p_L]$ is taken to be sufficiently large, with $\Delta x=(x_L-x_0)/L$, $\Delta p = (p_L-p_0)/L$.

\end{document}